\documentclass[10pt,reqno]{amsart}
\usepackage{amssymb,amsmath,amsfonts,amsthm,enumerate,enumitem,stmaryrd,tensor,mathtools,dsfont,upgreek,bbm,mathrsfs,nameref,xcolor,bm,tikz,stmaryrd,todonotes,graphicx,thmtools}
\usepackage[framemethod=TikZ]{mdframed}
\usetikzlibrary{arrows}
\usepackage[left=0.7 in, right=0.7 in, top=0.7 in, bottom=0.7 in]{geometry}

\usepackage{anyfontsize}
\usepackage{scrextend}
\usepackage{hyperref}
\usepackage{caption}
\usepackage{subcaption}
\changefontsizes{10pt}

\numberwithin{equation}{subsection}

\setlist[enumerate]{leftmargin=0.3in}
\setlist[itemize]{leftmargin=0.3in}
\definecolor{popblue}{RGB}{55,115,255}
\definecolor{lightbl}{RGB}{155,205,255}
\definecolor{depthbl}{RGB}{145,215,255}
\definecolor{fancyre}{RGB}{225,55,115}
\definecolor{lightgr}{RGB}{230,255,230}
\definecolor{darkgre}{RGB}{25,105,25}
\definecolor{darkblu}{RGB}{15,75,185}
\definecolor{mellowy}{RGB}{225,225,35}

\renewcommand{\tilde}[1]{\widetilde{#1}}
\renewcommand{\Bar}{\overline}

\renewcommand{\S}{\mathbb{S}}
\renewcommand{\v}{\mathsf{v}}
\renewcommand{\O}{\mathbb{O}}
\newcommand{\R}{\mathbb{R}}
\newcommand{\N}{\mathbb{N}}
\newcommand{\Z}{\mathbb{Z}}

\newcommand{\Y}{\mathbb{Y}}

\renewcommand{\P}{\mathbb{P}}

\newcommand{\imp}{\;\Rightarrow\;}
\newcommand{\m}{\mathrm}

\newcommand{\lv}{\lVert}
\newcommand{\rv}{\rVert}

\newcommand{\al}{\alpha}
\newcommand{\be}{\beta}
\newcommand{\es}{\varnothing}

\newcommand{\ep}{\varepsilon}
\newcommand{\f}{\frac}
\newcommand{\sig}{\sigma}
\newcommand{\gam}{\gamma}
\newcommand{\del}{\delta}

\newcommand{\pd}{\partial}

\newcommand{\grad}{\nabla}
\newcommand{\bpm}{\begin{pmatrix}}
\newcommand{\epm}{\end{pmatrix}}

\newcommand{\loc}{\m{loc}}

\renewcommand{\le}{\leqslant}
\renewcommand{\ge}{\geqslant}

\newcommand{\norm}[1]{\left\lv#1\right\rv}
\newcommand{\bnorm}[1]{\Big\lv#1\Big\rv}

\newcommand{\tnorm}[1]{\lv#1\rv}

\newcommand{\bp}[1]{\Big(#1\Big)}

\newcommand{\tp}[1]{(#1)}
\newcommand{\tfloor}[1]{\lfloor #1\rfloor}

\newcommand{\babs}[1]{\Big|#1\Big|}

\newcommand{\tabs}[1]{|#1|}

\newcommand{\bsb}[1]{\Big[{#1}\Big]}
\newcommand{\ssb}[1]{\big[{#1}\big]}
\newcommand{\tsb}[1]{[{#1}]}

\newcommand{\bcb}[1]{\Big\{{#1}\Big\}}
\newcommand{\tcb}[1]{\{{#1}\}}

\providecommand{\tbr}[1]{\langle #1 \rangle}

\def\XXint#1#2#3{{\setbox0=\hbox{$#1{#2#3}{\int}$ }
\vcenter{\hbox{$#2#3$ }}\kern-.6\wd0}}

\renewcommand{\bf}[1]{\mathbf{#1}}

\DeclareMathOperator{\supp}{supp}

\declaretheoremstyle[
    headfont=\bfseries\rmfamily\color{popblue}, bodyfont=\normalfont, 
    mdframed={
        linewidth=0.8pt,
        linecolor=lightbl, backgroundcolor=lightbl!5,
        skipabove=4pt,
    }
]{blueStyle}

\declaretheoremstyle[
    headfont=\bfseries\rmfamily\color{popblue}, 
    bodyfont=\itshape,
]{boringStyle}

\declaretheoremstyle[
    headfont=\bfseries\rmfamily\color{darkgre}, bodyfont=\normalfont, 
    mdframed={
        linewidth=0.8pt,
        linecolor=darkgre!35, backgroundcolor=lightgr!7,
        skipabove=4pt,
    }
]{greenStyle}

\declaretheoremstyle[
    headfont=\bfseries\rmfamily\color{fancyre}, bodyfont=\normalfont, 
    mdframed={
        linewidth=0.8pt,
        linecolor=fancyre!35, backgroundcolor=mellowy!6,
        skipabove=4pt,
    }
]{redStyle}

\newtheorem{prop}{\color{popblue}{Proposition}}[section]
\newtheorem{thm}[prop]{\color{popblue}{Theorem}}
\newtheorem{defn}[prop]{\color{popblue}{Definition}}
\newtheorem{lem}[prop]{\color{popblue}{Lemma}}
\newtheorem{coro}[prop]{\color{popblue}{Corollary}}
\newtheorem{rmk}[prop]{\color{popblue}{Remark}}

\newtheorem{innercustomthm}{\color{popblue}{Theorem}}
\newenvironment{customthm}[1]
{\renewcommand\theinnercustomthm{#1}\innercustomthm}
{\endinnercustomthm}

\author{R\u{a}zvan-Octavian Radu}
\address{
  Department of Mathematics\\
  Imperial College London\\
  180 Queen's Gate\\
  London SW7 2AZ, United Kingdom
}
\email[R.-O. Radu]{r.radu2@imperial.ac.uk}
\thanks{R.-O. Radu was partially supported by a Charlotte Elizabeth Procter Fellowship and a Chapman Fellowship.}
\author{Noah Stevenson}
\address{
  Forschungsinstitut f\"{u}r Mathematik\\
  ETH Z\"{u}rich\\
  R\"{a}mistrasse 101\\
  8092 Z\"{u}rich, Switzerland
}
\email[N. Stevenson]{noah.stevenson@math.ethz.ch}
\thanks{N. Stevenson was partially supported by an NSF Graduate Research Fellowship and a Hermann-Weyl-Instructorship.}

\DeclareFontFamily{U}{cbgreek}{}
\DeclareFontShape{U}{cbgreek}{m}{n}{
  <-6>    grmn0500
  <6-7>   grmn0600
  <7-8>   grmn0700
  <8-9>   grmn0800
  <9-10>  grmn0900
  <10-12> grmn1000
  <12-17> grmn1200
  <17->   grmn1728
}{}
\DeclareFontShape{U}{cbgreek}{bx}{n}{
  <-6>    grxn0500
  <6-7>   grxn0600
  <7-8>   grxn0700
  <8-9>   grxn0800
  <9-10>  grxn0900
  <10-12> grxn1000
  <12-17> grxn1200
  <17->   grxn1728
}{}

\DeclareRobustCommand{\qoppa}{%
  \text{\usefont{U}{cbgreek}{\normalorbold}{n}\symbol{19}}%
}

\makeatletter
\newcommand{\normalorbold}{%
  \ifnum\pdf@strcmp{\math@version}{bold}=\z@ bx\else m\fi
}
\makeatother

\title[Smooth and swirling steady vortex rings]{Smooth and swirling steady vortex rings near~the~Hill--Norbury~family}

\begin{document}
\begin{abstract}
  An investigation of local aspects of the space of axisymmetric traveling-wave solutions to the three-dimensional Euler equations near Hill's zero-flux spherical vortex and Norbury's small-flux vortex rings is the endeavor undertaken herein. The aforementioned family consists of idealized, swirl-free, and nonsmooth objects, their azimuthal relative vorticity being prescribed by a scaled characteristic function of the compact vortex core.
  
  We confront the desingularization question, both within the swirl-free solution subspace and upon admitting swirl, and prove that each member of the Hill--Norbury family of sufficiently small flux is the limit of traveling-wave solutions to the Euler equations that are smooth to infinite order and have compact toroidal vorticity support. In the latter swirling regime, however, a curious flexibility is revealed: total vorticity may diverge along sequences converging in a natural desingularization topology.

  Our construction, founded upon a five-dimensional reformulation of the Stokes stream function, is an amalgamation of delicate thin-domain analysis, singular-parameter Newton-type nonlinear inversion, and topological fixed-point arguments. To the best of the authors' knowledge, this constitutes the first verification that the small-flux Hill--Norbury family is nowhere isolated from the collection of smooth vortex rings.
\end{abstract}
\maketitle
\section{Introduction}\label{SO_introduction}


\subsection{Axisymmetric traveling waves and Hicks' equation}\label{SS_axi_TW_and_Hicks}

We study an incompressible, inviscid, and homogeneous fluid occupying all of three-dimensional Euclidean space $\R^3$. Its dynamics are described by a velocity field $\bf{u}:\R^+\times\R^3\to\R^3$ and a scalar pressure $\bf{p}:\R^+\times\R^3\to\R$ obeying Euler's equations
\begin{equation}\label{Eulers_Equations_Dynamic_Form}
  \pd_0\bf{u} + \bf{u}\cdot\grad\bf{u} + \grad\bf{p} = 0\text{ and }\grad\cdot\bf{u} = 0\text{ in }\R^+\times\R^3,
\end{equation}
where $\pd_0$ denotes differentiation with respect to the $\R^+ = (0,\infty)$ time factor and $\nabla = \tp{\pd_1,\pd_2,\pd_3}$ is the spatial gradient in $\R^3$.

Our interest lies in the construction and study of finite-energy traveling waves propagating at a nontrivial constant velocity. We may select coordinates such that these waves move in the positive $e_3$-direction, and thus for some speed $W>0$ we seek stationary profiles $u:\R^3\to\R^3$ and $p:\R^3\to\R$ such that
\begin{equation}\label{_TRAVELING_ANSATZ_}
  \bf{u}\tp{t,x} = u\tp{x - tWe_3}\text{ and }\bf{p}\tp{t,x} = p(x - tWe_3)\text{ for all }\tp{t,x}\in\R^+\times\R^3.
\end{equation}
Under the ansatz~\eqref{_TRAVELING_ANSATZ_}, the dynamic system~\eqref{Eulers_Equations_Dynamic_Form} is satisfied provided that
\begin{equation}\label{_TRAVELING_EULER_EQUATIONS_INTRO_FORM_}
  \tp{u - W e_3}\cdot\grad u + \grad p = 0 \text{ and }\grad\cdot u = 0\text{ in }\R^3.
\end{equation}

Next, let us impose axisymmetry about the direction of traveling wave propagation. Denote the standard $\R^3$-cylindrical coordinates by $r\in\R^+$, $\theta\in[0,2\pi)$, and $z\in\R$; these satisfy for $x = \tp{x_1,x_2,x_3}$ in $\R^3$ away from the symmetry axis:
\begin{equation}\label{cylindrical_coordinates_conventions}
  x_1 = r\cos\tp{\theta},\;x_2 = r\sin\tp{\theta},\text{ and }x_3 = z.
\end{equation}
The corresponding cylindrical orthonormal frame is
\begin{equation}\label{cylindrical_coordinates_orthonormal_frame}
  e_r = \cos\tp{\theta}e_1 + \sin\tp{\theta}e_2,\;e_\theta = -\sin\tp{\theta}e_1 + \cos\tp{\theta}e_2,\text{ and }e_z = e_3.
\end{equation}
An axisymmetric scalar function and the cylindrical components of an axisymmetric vector field are independent of $\theta$ and are therefore determined entirely by the values on the meridional plane
\begin{equation}\label{intro_meridional_plane}
  \Pi = \tp{0,\infty}\times\R.
\end{equation}

We make the further ansatz that the stationary profiles $u$ and $p$, which are to satisfy~\eqref{_TRAVELING_EULER_EQUATIONS_INTRO_FORM_}, are not only axisymmetric, but also determined via a Stokes stream function $\psi:\Pi\to\R$, a head derivative function $\Gamma:\R\to\R$, a circulation function $S:\R\to\R$, and a flux parameter $k\ge0$ through the identities
\begin{equation}\label{stokes_stream_function_ansatz}
  u = -\f{\pd_z\psi}{r}e_r + \f{\pd_r\psi}{r}e_z + \f{S(\Psi)}{r}e_\theta\text{ and }p = - H(\Psi) - \f12\tabs{u - We_3}^2
\end{equation}
where
\begin{equation}\label{defn_rel_stream_and_head}
  \Psi = \psi - \f{W}{2}r^2 - k\text{ and }H(s) = \int_0^s\Gamma\tp{\sig}\;\m{d}\sig.
\end{equation}

Then, $\grad\cdot u = 0$ is automatically enforced by~\eqref{stokes_stream_function_ansatz}; moreover, a direct calculation shows that the momentum equations in the traveling Euler system~\eqref{_TRAVELING_EULER_EQUATIONS_INTRO_FORM_} are satisfied, provided that
\begin{equation}\label{Hicks_equation_introduction_form}
  -L\psi = r^2\Gamma\tp{\Psi} + \tp{SS'}\tp{\Psi}\text{ in }\Pi\text{ where }L = r\pd_r\bp{\f{1}{r}\pd_r} + \pd_z^2.
\end{equation}
Equation~\eqref{Hicks_equation_introduction_form} goes by several names in the literature: for example, one sees Hicks' equation~\cite{Hicks1898}, the Bragg--Hawthorne equation~\cite{MR36103}, and the Grad--Shafranov equation~\cite{GradRubin1958,Shafranov1966}. In this paper we employ the first of these designations. As we are considering only finite-energy vortex rings, we supplement Hicks' equation~\eqref{Hicks_equation_introduction_form} with the normalization and decay conditions
\begin{equation}\label{_NORMALIZATION_AND_DECAY_CONDITIONS_}
  \psi = 0\text{ on }\pd\Pi,\;\psi(r,z) = \psi(r,-z),\text{ and }\psi,\f{\pd_z\psi}{r},\f{\pd_r\psi}{r}\to0\text{ as }r^2+z^2\to\infty.
\end{equation}
The second of these conditions centers the profile by requiring $\psi$ to be $z$-even.

Let us now comment on the geometric meaning of the quantities introduced in~\eqref{stokes_stream_function_ansatz} and~\eqref{defn_rel_stream_and_head}, which is more clearly seen when using the relative velocity
\begin{equation}
  u - We_3 = -\f{\pd_z\Psi}{r}e_r + \f{\pd_r\Psi}{r}e_z + \f{S\tp{\Psi}}{r}e_\theta,
\end{equation}
which satisfies $\tp{u - W e_3}\cdot\grad\Psi = 0$; thus, the level sets of $\Psi$ in $\Pi$ generate the stream surfaces of the relative velocity. Moreover, if $\mathcal{C}_{r,z}$ denotes the azimuthal circle $\tcb{\tp{r\cos\tp{\theta},r\sin\tp{\theta},z}\;:\;0\le\theta<2\pi}$ then the identities
\begin{equation}
  p + \f{1}{2}\tabs{u - W e_3}^2 = - H\tp{\Psi\tp{r,z}}\text{ and }\int_{\mathcal{C}_{r,z}}u\cdot\m{d}\ell = 2\pi S\tp{\Psi\tp{r,z}}
\end{equation}
show that on the stream surface $\tcb{\Psi = s}$, the head and circulation are identically constant and are given by $-H(s)$ and $2\pi S(s)$, respectively.

One computes that the vorticity under the ansatz~\eqref{stokes_stream_function_ansatz} is given by
\begin{equation}\label{_VORTICITY_LOOKY_HERE_}
  \grad\times u = \bp{r\Gamma\tp{\Psi} + \f{1}{r}\tp{SS'}\tp{\Psi}}e_\theta + \f{S'\tp{\Psi}}{r}\tp{-\pd_z\Psi e_r + \pd_r\Psi e_z}.
\end{equation}
Wanting to study vortex rings with internal circulation, we impose the following further restrictions on $\Gamma$ and $S$, namely
\begin{equation}\label{_VORTEX_RING_ANSATZ_}
  S(s) = \Gamma(s) = 0\text{ for }s\le0,\;\Gamma\tp{s}>0\text{ for }s>0,\text{ and }\supp S\Subset\R^+.
\end{equation}
We then see from~\eqref{_VORTICITY_LOOKY_HERE_} and~\eqref{_VORTEX_RING_ANSATZ_} that the bounded set $\tcb{\Psi>0}\subset\Pi$ is the meridional vortex core, in that 
\begin{equation}\label{_MER_REV_IS_SUP_VORT_}
  \supp\tp{\grad\times u} = \Bar{\tcb{x\in\R^3\;:\;\Psi(\tp{x_1^2+x_2^2}^{1/2},x_3)>0}},
\end{equation}
where the overset bar denotes the topological closure.

We now comment on the meaning of the parameter $k\ge0$, which has units of volume flux and admits a direct interpretation. The normalization~\eqref{_NORMALIZATION_AND_DECAY_CONDITIONS_} gives $\Psi\tp{0,z} = -k\le 0$. For the vortex rings considered below, the meridional vortex core intersects the equatorial plane $z=0$. Let
\begin{equation}\label{_the_size_of_the_hole_}
  r_- = \sup\tcb{r\ge0\;:\;\Psi(\rho,0)\le0\text{ for all }0\le\rho\le r}.
\end{equation}
For these configurations, $r_-<\infty$, and $r_->0$ whenever $k>0$. For the disk
\begin{equation}\label{_together_with_this_disk_}
  D_- = \tcb{x\in\R^3\;:\;x_3 = 0,\;x_1^2+x_2^2\le r_-^2}
\end{equation}
we have
\begin{equation}
  \int_{D_-}\tp{u - W e_3}\cdot e_3\;\m{d}\mathcal{H}^2 = 2\pi\int_0^{r_-}\pd_r\Psi\tp{\rho,0}\;\m{d}\rho = 2\pi k,
\end{equation}
where $\mathcal{H}^d$ denotes the $d$-dimensional Hausdorff measure; here we used $\Psi(r_-,0) = 0$ and $\Psi(0,0) = -k$. Thus $2\pi k$ is precisely the flux of the relative velocity field $u - We_3$ through the `hole' $D_-$ at the center of the vortex ring. See Figure~\ref{fig_norbury_flux} for a relevant visualization.

\begin{center}
  \begin{minipage}{0.7\textwidth}
    \centering
    \includegraphics[width=\textwidth]{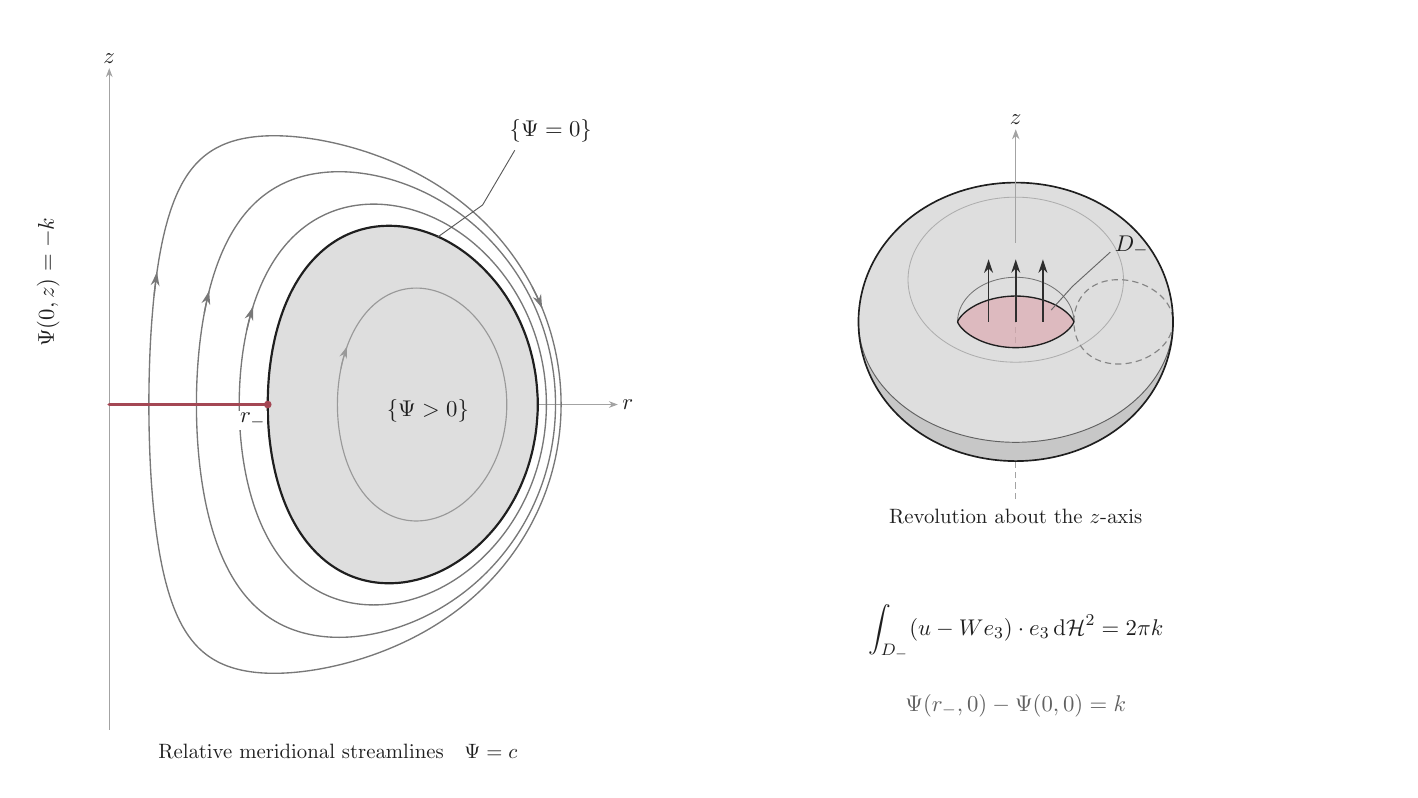}
  \end{minipage}
\end{center}

{
\captionsetup{hypcap=false,width=0.9\textwidth}
\captionof{figure}{
  An example of a vortex ring arising from Hicks' equation~\eqref{Hicks_equation_introduction_form}. On the left, we depict the meridional vortex core $\tcb{\Psi>0}$, with $\Psi$ defined in~\eqref{defn_rel_stream_and_head}, together with several relative streamlines $\tcb{\Psi=c}$ and, in red, the segment connecting $(0,0)$ to $(r_-,0)$, where $r_-$ is defined in~\eqref{_the_size_of_the_hole_}. On the right, we depict the corresponding three-dimensional support of the vorticity, as in~\eqref{_MER_REV_IS_SUP_VORT_}, together with the disk $D_-$, defined in~\eqref{_together_with_this_disk_}, spanning its central hole. The arrows indicate the relative velocity $u-We_3$, whose flux through $D_-$ is $2\pi k$.}
\label{fig_norbury_flux}
}

\subsection{The solutions of Hill and Norbury, nondimensionalization}\label{SS_sol_HN_and_nondim}

The preceding discussion assumes sufficient regularity to interpret Hicks' equation~\eqref{Hicks_equation_introduction_form} classically. While we shall only construct smooth solutions to the Euler and Hicks equations, we are in a certain sense working perturbatively around certain distinguished idealized and non-smooth vortex rings in which the vorticity law, or head derivative function, is discontinuous. As such, we shall require the weak formulation of the vortex ring problem developed in the works of Fraenkel and Berger~\cite{MR422916}, Amick and Fraenkel~\cite{MR816615}, and Amick and Turner~\cite{MR929976}.

Given $\psi,\tilde{\psi}\in C^\infty_{\m{c}}\tp{\Pi}$, we define the bilinear kinetic energy pairing and the associated norm via
\begin{equation}
  \tbr{\psi,\tilde{\psi}}_{\mathtt{o}} = \int_{\Pi}\f{1}{r}\tp{\pd_r\psi\pd_r\tilde{\psi} + \pd_z\psi\pd_z\tilde{\psi}}\;\m{d}r\;\m{d}z\text{ and }\tnorm{\psi}_{\mathtt{o}} = \tsb{\tbr{\psi,\psi}_{\mathtt{o}}}^{1/2}.
\end{equation}
The completion of $C^\infty_{\m{c}}\tp{\Pi}$ under the norm $\tnorm{\cdot}_{\mathtt{o}}$ is then denoted $H(\Pi)$. The weak formulation of the generalized vortex ring problem is as follows. Given functions $S$ and $\Gamma$ satisfying~\eqref{_VORTEX_RING_ANSATZ_} and parameters $k\ge0$ and $W>0$, find a function $\psi\in H(\Pi)$ such that
\begin{equation}\label{_GENERALIZED_VORTEX_RING_PROBLEM_}
  \tbr{\psi,\tilde{\psi}}_{\mathtt{o}} = \int_{\Pi}\bp{\f{1}{r}\tp{SS'}\bp{\psi - \f{W}{2}r^2 - k} + r\Gamma\bp{\psi - \f{W}{2}r^2 - k}}\tilde{\psi}\;\m{d}r\;\m{d}z\text{ for all }\tilde{\psi}\in C^\infty_{\m{c}}\tp{\Pi}.
\end{equation}

The idealized problem is obtained by specializing to the case
\begin{equation}\label{_IDEALIZED_VORTEX_RING_PROBLEM_}
  S = 0\text{ and }\Gamma = \lambda\mathds{1}_{\tp{0,\infty}}
\end{equation}
for a vortex strength parameter $\lambda>0$. In this idealized case, the Stokes stream function $\psi\in H(\Pi)$ satisfies, for all $\tilde{\psi}\in C^\infty_{\m{c}}\tp{\Pi}$,
\begin{equation}\label{_IDEALIZED_VORTEX_RING_PROBLEM_EQUATION_}
  \tbr{\psi,\tilde{\psi}}_{\mathtt{o}} = \lambda \int_{A\tp{\psi,k}}\tilde{\psi} r\;\m{d}r\;\m{d}z
\end{equation}
where the set $A(\psi,k)$ is the vortex core, defined as
\begin{equation}
  A\tp{\psi,k} = \bcb{\tp{r,z}\in\Pi\;:\;\psi(r,z) > \f{W}{2}r^2 + k}.
\end{equation}
The reconstructed fluid vorticity is then
\begin{equation}\label{_RECON_FLUID_VORT_}
  \grad\times u = \lambda r\mathds{1}_{A\tp{\psi,k}}e_\theta,
\end{equation}
in the sense of distributions. In this idealized formulation, we emphasize that the vortex core $A\tp{\psi,k}$ is not prescribed in advance; rather, its boundary is determined simultaneously with $\psi$.

The earliest known example of a steady vortex `ring' is the spherical vortex discovered in the late nineteenth century by Micaiah Hill~\cite{Hill1894}. We denote its stream function by $\psi^{\m{H}}$ here. Hill's spherical vortex is an explicit $z$-even solution to the idealized problem~\eqref{_IDEALIZED_VORTEX_RING_PROBLEM_EQUATION_} in the zero flux case $k = 0$; in fact, Hill's solution is known to be the unique $z$-even weak solution of zero flux, see Amick and Fraenkel~\cite{MR816615}. The solid of revolution generated by the closure of its meridional vortex core is a ball of radius
\begin{equation}\label{_HILLS_VORTEX_RADIUS_}
  R_{\m{H}} = \bp{\f{15 W}{2\lambda}}^{1/2}\text{ so that }A(\psi^{\m{H}},0) = \tcb{\tp{r,z}\in\Pi\;:\;r^2 + z^2 < R_{\m{H}}^2}.
\end{equation} 
One has
\begin{equation}\label{_HILLS_VORTEX_STREAM_FUNCTION_}
  \psi^{\m{H}}\tp{r,z} = \f{\lambda r^2}{30}\begin{cases}\tp{5R_{\m{H}}^2 - 3\tp{r^2 + z^2}}&\text{ if }r^2 + z^2\le R_{\m{H}}^2,\\
    2R_{\m{H}}^5\tp{r^2 + z^2}^{-3/2}&\text{if }r^2 + z^2\ge R_{\m{H}}^2,
  \end{cases}
  \text{ for all }\tp{r,z}\in\Pi.
\end{equation}

Well over half a century after Hill's discovery, John Norbury~\cite{MR302044} considered the question of constructing nearby solutions to the idealized problem~\eqref{_IDEALIZED_VORTEX_RING_PROBLEM_EQUATION_} for small positive values of the flux parameter $k$. He was ultimately successful and constructed a non-explicit family of $z$-even
weak solutions
\begin{equation}
  \tcb{\psi_k^{\m{N}}}_{0<k<k_{\m{N}}}\subset H(\Pi),
\end{equation}
for some small $k_{\m{N}}>0$. In striking contrast with Hill's solution, for
$k>0$ the solids of revolution generated by the closures of the meridional vortex cores $\Bar{A(\psi_k^{\m{N}},k)}$ are topological solid tori not intersecting the $e_3$-axis. Their inner radii are of order $\tp{k/W}^{1/2}$. Norbury's solutions converge to Hill's spherical vortex in the vanishing-flux limit $k\to0$.

The parameters $W$ and $\lambda$ merely fix the velocity and length scales of the idealized problem. We now endeavor to nondimensionalize by first letting $\m{L} = R_{\m{H}}$ and $\m{U} = W/2$ denote our choice of length and velocity scales. For the duration of this discussion, we let an overset tilde denote dimensionless quantities. Set $t = \m{L}\m{U}^{-1}\tilde{t}$, $x = \m{L}\tilde{x}$, $s = \m{U}\m{L}^2\tilde{s}$,
\begin{multline}\label{non_dimensionalized_functions}
  \bf{u}\tp{t,x} = \m{U}\tilde{\bf{u}}\tp{\tilde{t},\tilde{x}},\;\bf{p}\tp{t,x} = \m{U}^2\tilde{\bf{p}}\tp{\tilde{t},\tilde{x}},\;\psi(r,z) = \m{U}\m{L}^2\tilde{\psi}\tp{\tilde{r},\tilde{z}},\;\Psi(r,z) = \m{U}\m{L}^2\tilde{\Psi}\tp{\tilde{r},\tilde{z}},\\
  k = \m{U}\m{L}^2\tilde{k},\;H(s) = \m{U}^2\tilde{H}\tp{\tilde{s}},\;S(s) = \m{U}\m{L}\tilde{S}\tp{\tilde{s}},\;\text{and}\;\Gamma(s) = \m{U}\m{L}^{-2}\tilde{\Gamma}\tp{\tilde{s}}.
\end{multline}
With these choices, Euler's equations~\eqref{Eulers_Equations_Dynamic_Form} and~\eqref{_TRAVELING_EULER_EQUATIONS_INTRO_FORM_}, the traveling wave and stream function ans\"{a}tze~\eqref{_TRAVELING_ANSATZ_} and~\eqref{stokes_stream_function_ansatz}, Hicks' equation~\eqref{Hicks_equation_introduction_form}, and the weak formulations~\eqref{_GENERALIZED_VORTEX_RING_PROBLEM_} and~\eqref{_IDEALIZED_VORTEX_RING_PROBLEM_EQUATION_} retain precisely the same forms in the dimensionless variables. Moreover,
\begin{equation}\label{key_transformations_of_the_nondimensionalization}
  \tilde{W} = 2,\;\tilde{\lambda} = \m{L}^2\m{U}^{-1}\lambda = 15,\;\tilde{R_{\m{H}}} = 1,\;\tilde{k} = \tp{\m{U}\m{L}^2}^{-1}k = 4\lambda\tp{15W^2}^{-1}k.
\end{equation}
We henceforth work exclusively in these units and suppress the overset tildes. We note in particular that the nondimensionalized idealized vortex ring problem is to find $\psi\in H(\Pi)\setminus\tcb{0}$ such that
  \begin{equation}\label{_NON_DIM_NORBURY_WEAK_PROBLEM_}
    \tbr{\psi,\tilde{\psi}}_{\mathtt{o}} = 15\int_{A\tp{\psi,k}}\tilde{\psi} r \;\m{d}r\;\m{d}z\text{ for all }\tilde{\psi}\in C^\infty_{\m{c}}\tp{\Pi}\text{ with }A\tp{\psi,k} = \tcb{\tp{r,z}\in\Pi\;:\;\psi(r,z)>r^2 + k}.
  \end{equation}

\begin{defn}[The Hill--Norbury family]\label{defn on the Hill Norbury family}

  Fix $\pmb{\mathtt{k}}>0$ sufficiently small so that Norbury's centered branch~\cite{MR302044} of weak solutions to the nondimensional weak formulation of the idealized steady vortex ring problem~\eqref{_NON_DIM_NORBURY_WEAK_PROBLEM_} is defined for every $0<k\le\pmb{\mathtt{k}}$. We set $\psi^0 = \psi^{\m{H}}$ and $\psi^k = \psi^{\m{N}}_k$ for $0<k\le\pmb{\mathtt{k}}$, and denote by $u^k\in\tp{L^2\cap\m{LL}}\tp{\R^3;\R^3}$ (see~\eqref{definition_of_the_log_lip_norm}) the corresponding finite-energy, swirl-free velocity field
  \begin{equation}\label{Hill_Norbury_Velocity}
    u^k = -\f{\pd_z\psi^k}{r}e_r + \f{\pd_r\psi^k}{r}e_z.
  \end{equation}
  Each $u^k$ is a weak traveling-wave solution of Euler's equations with speed $2$ and, in the sense of distributions,
  \begin{equation}\label{_UNFORTUNATELY_THIS_GETS_REFERENCED_}
    \grad\times u^k = 15r\mathds{1}_{A\tp{\psi^k,k}}e_\theta.
  \end{equation}

  We refer to $\tcb{u^k}_{0\le k\le\pmb{\mathtt{k}}}$ as the portion of the Hill--Norbury family considered in this paper. See Figure~\ref{fig_Hill_Norbury_family} for a numerical depiction of some of these vortex rings, including rings beyond the perturbative regime.
\end{defn}

\begin{center}
  \begin{minipage}{0.7\textwidth}
    \centering
    \includegraphics[width=\textwidth]{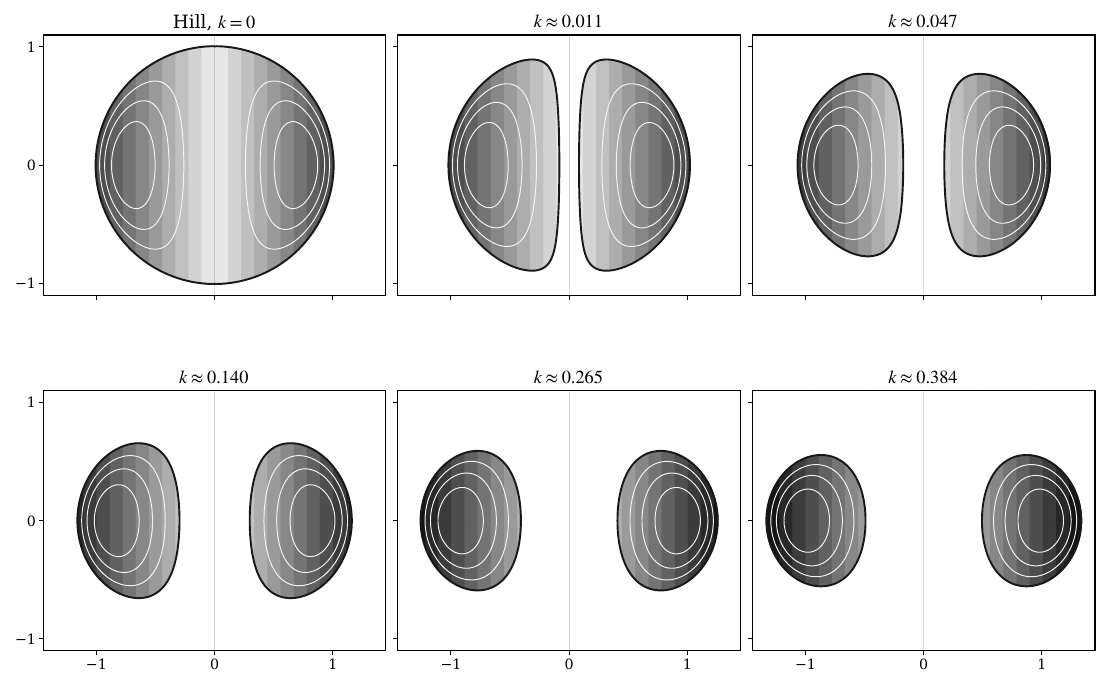}
  \end{minipage}
\end{center}

{
\captionsetup{hypcap=false,width=0.9\textwidth}
\captionof{figure}{Depicted here, for various values of $k\ge0$, are members of the numerically continued nondimensionalized Hill--Norbury family, represented by solutions $\psi$ to the problem~\eqref{_NON_DIM_NORBURY_WEAK_PROBLEM_}. The $z$-coordinate varies vertically and the signed radial coordinate $r$ varies horizontally. The black contour is the boundary of the meridional vortex core $A(\psi,k)$, which has been reflected across the symmetry axis. The discretely binned grayscale records $\tabs{\grad\times u}$, with white being zero and increasing magnitude with increasing darkness. The overlaid contours are level curves of the relative stream function $\Psi$, defined in terms of $\psi$ via~\eqref{defn_rel_stream_and_head}. The first panel is Hill's spherical vortex, whose explicit expression is recorded in equation~\eqref{_HILLS_VORTEX_STREAM_FUNCTION_}. The remaining panels sample the Norbury family for increasing values of $k$; we emphasize that only an initial small-flux portion of this family is discussed by our theorems.}
\label{fig_Hill_Norbury_family}
}

\subsection{Literature review and niche identification}\label{SS_lit_review}
Having now fixed in Definition~\ref{defn on the Hill Norbury family} the classical family that serves as the target of our later constructions, we now place it within the broader theory of steady vortex rings and elucidate the questions addressed by this paper.

Following Helmholtz's foundational study of vortex motion~\cite{Helmholtz1858}, the early theory of steady vortex rings developed around two rather different limiting configurations. Hill~\cite{Hill1894} discovered an explicit finite-energy traveling wave whose vorticity fills a ball, furnishing the basic large-core model for the subject. At the opposite extreme, Fraenkel~\cite{Fraenkel1970,Fraenkel1972} constructed rings of vanishing cross-section and determined their asymptotic structure. Norbury~\cite{MR302044} approached the problem from Hill's end, constructing a local family of toroidal idealized rings that collapses onto the spherical vortex as the flux tends to zero; his subsequent numerical study~\cite{Norbury1973} traced this family well beyond the perturbative regime. Thus, by the early 1970s, the two extreme shapes of the classical theory were visible, together with strong numerical evidence that they belong to a common family.

A second line of work sought vortex rings through general existence methods rather than continuation from a distinguished explicit solution. Fraenkel and Berger~\cite{MR422916} initiated a broad variational theory for swirl-free rings, while Ni~\cite{MR583638} replaced their constrained minimization by a mountain-pass construction, working for a slightly more restricted class of head functions. Of particular importance here, Ni was also the first to use the transformation that identifies the axisymmetric stream-function operator with the Laplacian acting on cylindrically symmetric functions in $\R^5$ in the vortex ring problem. Related nonlinear and variational constructions were developed by Berger and Fraenkel~\cite{MR589430} and by Friedman and Turkington~\cite{MR628444}, the latter building on Benjamin's vorticity-space formulation~\cite{MR671099}; Turkington~\cite{MR977488} subsequently treated rings with nonzero swirl and helicity. We also mention the later filament-concentration constructions of de Valeriola and Van Schaftingen~\cite{MR3101789} and Cao, Wan, and Zhan~\cite{MR4152780}.

The idealized Hill--Norbury problem possesses a separate rigidity and continuation theory. Amick and Fraenkel proved that Hill's spherical vortex is the unique centered solution of the zero-flux problem~\cite{MR816615} and later established conditional uniqueness of the small-flux Norbury rings~\cite{MR918795}. Amick and Turner~\cite{MR929976} continued Norbury's local branch to an unbounded, closed, and connected global set. Together with the asymptotic and numerical studies of Fraenkel~\cite{Fraenkel1972} and Norbury~\cite{Norbury1973}, the Amick--Turner branch supports the classical vortex ring interpolation conjecture, which predicts a continuous family joining Hill's vortex to the thin-core rings. We shall not address this global conjecture here, although the authors are optimistic that machinery developed below might prove relevant in the future.

The questions considered in this paper are instead local to the Hill--Norbury family. These idealized rings have a discontinuous vorticity law, and their stream functions are of class $C^{1,1}$ but not $C^2$ across the vortex boundary. Norbury's fixed-point construction and the Amick--Turner continuation are strongly tied to this free-boundary structure. Conversely, the aforementioned general variational theories permit smoother vorticity laws, but do little to locate the resulting solutions near any prescribed member of the Hill--Norbury family. It is therefore not known from the existing theory whether these classical rings are isolated from smooth steady vortex rings.

This leads to two problems. First, can one desingularize each small-flux Hill--Norbury ring within the swirl-free class while retaining convergence of the velocity in finite energy and in $C^\alpha$ for every $0\le\alpha<1$? Second, what becomes possible once swirl is admitted? In particular, can one remain arbitrarily close to the Hill--Norbury velocity fields while independently controlling quantities such as the helicity and total meridional vorticity? The first question concerns non-isolation inside a rigid symmetry class; the second concerns the freedom in the vorticity that is compatible with strong convergence of the velocity.

These two questions have both physical and mathematical significance. Showing that the Hill--Norbury family belongs to the closure of the smooth traveling-wave solution space demonstrates that the classical rings are not merely artifacts of a discontinuous constitutive law. At the same time, the problem gives a concrete setting in which to study the local geometry of the steady Euler solution set.

A further motivation comes from convex integration, a family of methods originating in Nash's work on isometric embeddings~\cite{MR0065993} and entering fluid mechanics with the work of De Lellis and Sz\'ekelyhidi~\cite{MR2600877,MR3090182}. Recently, Bru\`e, Colombo, and Kumar~\cite{BrueColomboKumar2026} constructed nonunique two-dimensional Euler flows using building blocks based on the Lamb--Chaplygin dipole, while Nguyen and Wang~\cite{NguyenWang2026} used localized moving Hill vortices in a convex-integration construction for the three-dimensional Navier--Stokes equations. Such schemes rely on the iterative addition of highly oscillatory perturbations with carefully controlled frequency content and therefore require sufficiently smooth building blocks. In the works above, this smoothness is obtained by mollification, at the cost of losing the exact traveling-wave character of the underlying Euler solutions. The exact smooth desingularization developed in this paper avoids this tradeoff and may therefore provide building blocks with sharper error estimates, potentially leading to stronger flexibility results.

\subsection{Main results}\label{SS_main_results}

Before stating our main results, we fix one final piece of notation. For $d\in\N^+$ and $\es\neq U\subseteq\R^d$, the log-Lipschitz space $\m{LL}\tp{U}$ consists of continuous and bounded functions $f:\Bar{U}\to\R$ for which
\begin{equation}\label{definition_of_the_log_lip_norm}
  \tnorm{f}_{\m{LL}\tp{U}} = \tnorm{f}_{C^0_b\tp{\Bar{U}}} + \sup\bcb{\f{\tabs{f(x) - f(y)}}{|x - y|\tp{1 + \tabs{\log\tabs{x - y}}}}\;:\;x,y\in U,\;0<|x-y|\le 1}<\infty.
\end{equation}
For vector-valued functions we use the same notation and definition, with the absolute value replaced by the Euclidean norm. We shall also consider the related space
\begin{equation}\label{definition_of_LL1_space}
  \m{LL}^1\tp{U} = \tcb{f\in C^1_{\m{b}}\tp{\Bar{U}}\;:\;\grad f\in\m{LL}\tp{U;\R^d}}\text{ with the norm }\tnorm{f}_{\m{LL}^1\tp{U}} = \tnorm{f}_{C^1_{\m{b}}\tp{\Bar{U}}} + \tnorm{\grad f}_{\m{LL}\tp{U}}.
\end{equation}

The first main theorem addresses the first question raised above.
\begin{customthm}{I}[Swirl-free desingularization, proved in Section~\ref{SS_proofs_main_results}]\label{main_thm_1_swirl_free}
  There exist constants $0<k^\star\le\pmb{\mathtt{k}}$ and $0<C<\infty$ such that for all $0\le k\le k^\star$ there exists a sequence of axisymmetric functions
  \begin{equation}\label{_CITE_ME_PLEASE_}
    \tcb{u_n}_{n\in\N}\subset\tp{L^2\cap C^\infty}\tp{\R^3;\R^3}\text{ with }\sup_{n\in\N}\tnorm{u_n}_{\tp{L^2\cap\m{LL}}\tp{\R^3}}\le C
  \end{equation}
  such that the following hold.
  \begin{enumerate}
    \item Upon letting $u^k\in\tp{L^2\cap\m{LL}}\tp{\R^3;\R^3}$ denote the Hill--Norbury velocity field of flux $k$ (see Definition~\ref{defn on the Hill Norbury family}), for every $0<\al<1$ the limit
    \begin{equation}\label{_CITE_ME_TOO_PLEASE_}
      \tnorm{u_n - u^k}_{L^2\tp{\R^3}} + \tnorm{u_n - u^k}_{C^\al_{\m{b}}\tp{\R^3}}\to0\text{ as }n\to\infty
    \end{equation}
    holds.
    \item For all $n\in\N$, we have $e_\theta\cdot u_n = 0$ away from the symmetry axis and there exists $p_n\in C^\infty\tp{\R^3}$ such that the traveling incompressible Euler equations
    \begin{equation}
      \tp{u_n - 2 e_3}\cdot\grad u_n + \grad p_n = 0\text{ and }\grad\cdot u_n = 0\text{ in }\R^3
    \end{equation}
    are satisfied.
    \item For all $n\in\N$, the radial component of $u_n$ is odd in $x_3$, while its azimuthal and axial components are even in $x_3$; moreover, $e_r\cdot u_n(x)>0$ whenever $x_1^2+x_2^2>0$ and $x_3>0$. The set $\supp\tp{\grad\times u_n}$ is contained within $\tcb{x\in\R^3\;:\;|x|<2,\;x_1^2+x_2^2\neq0}$ and is a topological solid torus with smooth boundary.
  \end{enumerate}
\end{customthm}

A couple of remarks are in order. The overarching takeaway of Theorem~\ref{main_thm_1_swirl_free} is that the portion of the Hill--Norbury family with sufficiently small flux, including Hill's solution itself, is nowhere isolated from the collection of infinitely smooth, axisymmetric, and swirl-free steady vortex ring solutions to Euler's equations.

The desingularization sequences $\tcb{u_n}_{n\in\N}$ are generated by sequences of smooth Stokes stream functions $\tcb{\psi_n}_{n\in\N}$, exactly through identity~\eqref{stokes_stream_function_ansatz}, that solve the Hicks equation~\eqref{Hicks_equation_introduction_form} with identically zero circulation and, as suggested by the idealized Heaviside vorticity law~\eqref{_IDEALIZED_VORTEX_RING_PROBLEM_}, head derivative functions $\Gamma_n$ that are smooth approximations of $15\mathds{1}_{\tp{0,\infty}}$.

\begin{figure}[ht]
    \centering

    \includegraphics[width=0.7\textwidth]
        {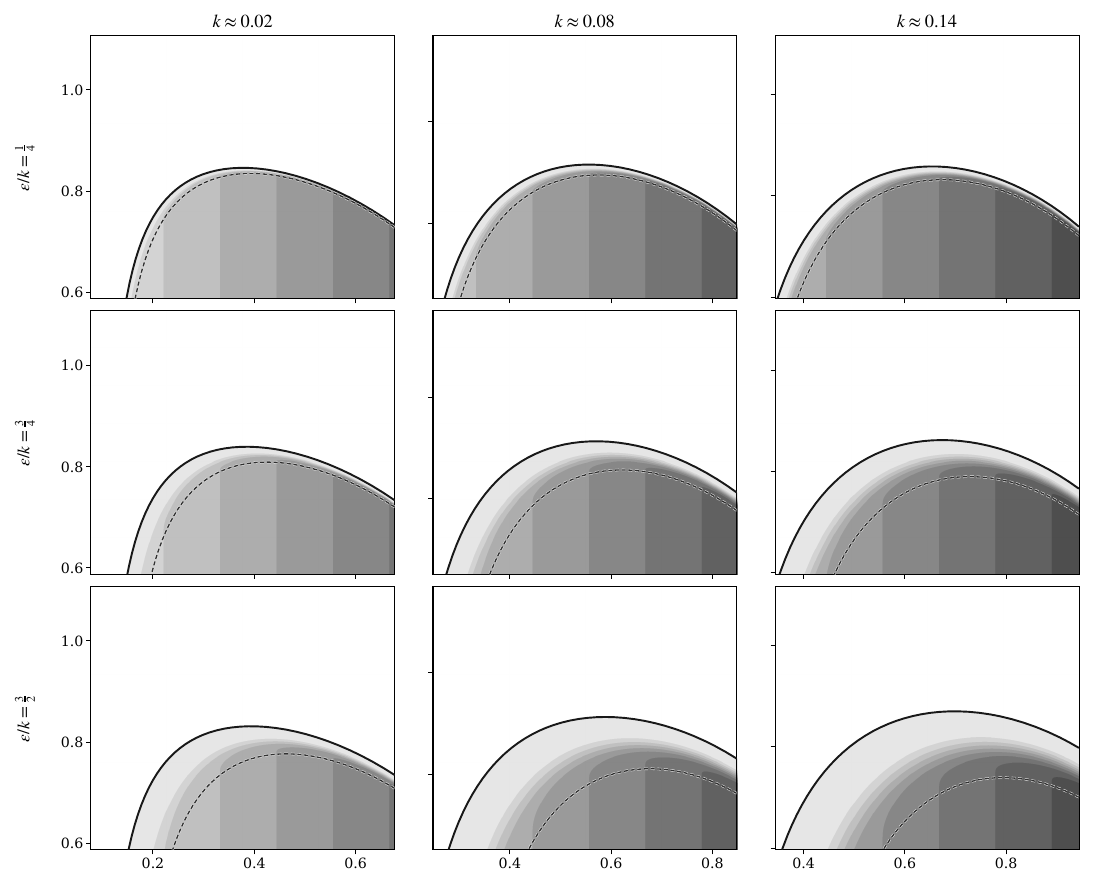}

    \caption{Depicted here, for various values of the flux $k>0$ and
    regularization ratio $\ep/k$, are numerical solutions of the smooth
    swirl-free problem described in~\eqref{thm_I_specialized_hicks}.
    Each panel gives a magnified view of the upper portion of the meridional
    vortex core, with the $z$-coordinate varying vertically and the radial
    coordinate $r$ varying horizontally. The solid black contour is the
    boundary $\tcb{\Psi=0}$ of the meridional vortex core, while the dashed
    contour is the level set $\tcb{\Psi=\ep}$. The discretely binned grayscale
    records $\tabs{\grad\times u}$, with white being zero and increasing
    magnitude with increasing darkness. Thus, between the solid and dashed
    contours the smooth head derivative
    $\Gamma(s)=15\gam(s/\ep)$ transitions from zero to its full value $15$,
    while on the interior side of the dashed contour one has
    $\Gamma(\Psi)=15$. The columns illustrate increasing values of $k$, while
    the rows illustrate increasing values of $\ep/k$. In particular, as the
    regularization scale decreases relative to the flux, the transition layer
    collapses toward the vortex boundary and the discontinuous vorticity law
    of the corresponding idealized Hill--Norbury ring is recovered.}

    \label{fig_Hill_Norbury_desingularization}
\end{figure}

More precisely, there exists a fixed transition profile $\gam\in C^\infty\tp{\R}$ with $\gam(s) = 0$ for $s\le 0$ and $\gam(s) = 1$ for $s\ge 1$ and a small parameter $\del_\star>0$ such that upon setting
\begin{equation}\label{form_of_the_smoothed_out_heaviside_and_flux}
  \Gamma_n\tp{s} = 15\gam\tp{2^ns/\del_\star}\text{ and }k_n = \begin{cases}
    k&\text{if }k>0,\\
    \del_\star/2^n&\text{if }k=0
  \end{cases}
\end{equation}
we have
\begin{equation}\label{thm_I_specialized_hicks}
  -L\psi_n = r^2\Gamma_n\tp{\psi_n - r^2 - k_n}\text{ in }\Pi\text{ and }u_n = -\f{\pd_z\psi_n}{r}e_r + \f{\pd_r\psi_n}{r}e_z.
\end{equation}
The effect of the regularization parameter on the vorticity profile and the geometry of the vortex core is illustrated numerically in
Figure~\ref{fig_Hill_Norbury_desingularization}.

Notice, in particular, that each element along the sequences constructed for $0<k\le k^\star$ has precisely the same flux value as the target Norbury ring. On the other hand, when $k = 0$, and so the target is Hill's solution, the members along our desingularization sequences have the strictly positive but decreasing flux values $n\mapsto\del_\star/2^n$. This leads to a notable geometric feature of our desingularization at $k = 0$. Every member of the approximating sequence has vorticity support equal to a smooth topological solid torus separated from the symmetry axis, whereas the vorticity support of the limiting Hill vortex is a ball. Thus the convergence in $L^2\cap C^\alpha_{\m b}$ for every $0<\alpha<1$ in~\eqref{_CITE_ME_TOO_PLEASE_}, even together with the uniform log-Lipschitz bound~\eqref{_CITE_ME_PLEASE_}, permits a change in the topology of the vortex core in the limit.

Our next main theorem considers the second question raised above on what happens around the Hill--Norbury family when nontrivial swirl is permitted.
\begin{customthm}{II}[Flexible desingularization with swirl, proved in Section~\ref{SS_proofs_main_results}]\label{main_thm_2_with_swirl}
  Let $C$ and $k^\star$ be as in Theorem~\ref{main_thm_1_swirl_free}. For all $0\le k\le k^\star$ and all choices of sequences 
  \begin{equation}\label{_CONSTRAINTS_ON_HELI_TOT_MER_VORT_SEQ_}
      \tcb{h_n}_{n\in\N}\subset[-1/C,1/C],\;\tcb{m_n}_{n\in\N}\subset[1,\infty)\text{ for which }h_n\to0\text{ as }n\to\infty,
  \end{equation}
  there exists a sequence of axisymmetric functions
  \begin{equation}
    \tcb{u_n}_{n\in\N}\subset\tp{L^2\cap C^\infty}\tp{\R^3;\R^3}\text{ with }\sup_{n\in\N}\tnorm{u_n}_{\tp{L^2\cap\m{LL}}\tp{\R^3}}\le C
  \end{equation}
  such that the following hold.
  \begin{enumerate}
    \item Upon letting $u^k\in\tp{L^2\cap\m{LL}}\tp{\R^3;\R^3}$ denote the Hill--Norbury velocity field of flux $k$ (see Definition~\ref{defn on the Hill Norbury family}), for every $0<\al<1$ the limit
    \begin{equation}\label{NATURAL_DESING_TOPO}
      \tnorm{u_n - u^k}_{L^2\tp{\R^3}} + \tnorm{u_n - u^k}_{C_{\m{b}}^\al\tp{\R^3}}\to0\text{ as }n\to\infty
    \end{equation}
    holds.
    \item For all $n\in\N$, we have $e_\theta\cdot u_n \not\equiv 0$ and there exists $p_n\in C^\infty\tp{\R^3}$ such that the traveling incompressible Euler equations
    \begin{equation}
      \tp{u_n - 2 e_3}\cdot\grad u_n + \grad p_n = 0\text{ and }\grad\cdot u_n = 0\text{ in }\R^3
    \end{equation}
    are satisfied. Moreover, the inclusions
    \begin{equation}\label{_LOCATION_OF_THE_SWIRL_AWAY_FROM_RING_BOUNDARY_}
      \supp\tp{e_\theta\cdot u_n}\subseteq\tcb{x\in\R^3\;:\;C^{-2}\le x_1^2 + x_2^2\text{ and }|x|\le1 - C^{-1}}
    \end{equation}
    hold.
    \item For all $n\in\N$, the radial component of $u_n$ is odd in $x_3$, while its azimuthal and axial components are even in $x_3$; moreover, $e_r\cdot u_n(x)>0$ whenever $x_1^2+x_2^2>0$ and $x_3>0$. The set $\supp\tp{\grad\times u_n}$ is contained within $\tcb{x\in\R^3\;:\;|x|<2,\;x_1^2+x_2^2\neq0}$ and is a topological solid torus with smooth boundary.
    \item The helicity and total meridional vorticity of $u_n$ are, for all $n\in\N$, given by
    \begin{equation}\label{_THANKS_FOR_LOOKING_AT_THIS_LABEL_}
      h_n = \int_{\R^3}u_n\cdot\grad\times u_n\text{ and }m_n = \int_{\R^3}\tabs{\tp{I - e_\theta\otimes e_\theta}\tp{\grad\times u_n}},
    \end{equation}
    respectively.
  \end{enumerate}
\end{customthm}

Several comments on the second main theorem are warranted. The upshot is that once we are authorized to desingularize outside of the class of swirl-free traveling wave solutions, we are afforded very little control of the vorticity simply by being nearby the Hill--Norbury family. This family of idealized vortex rings possesses discontinuous yet bounded vorticity which, from~\eqref{_UNFORTUNATELY_THIS_GETS_REFERENCED_}, is completely aligned with the azimuthal direction. However, Theorem~\ref{main_thm_2_with_swirl} allows us to specify the total meridional vorticity sequence $\tcb{m_n}_{n\in\N}\subset[1,\infty)$ completely arbitrarily; in particular, by using the latter identity of the fourth item, selecting $m_n\to\infty$ forces
\begin{equation}
  \int_{\R^3}\tabs{\tp{I - e_\theta\otimes e_\theta}\tp{\grad\times u_n}}\to\infty\text{ as }n\to\infty,
\end{equation}
but, for the limiting Hill--Norbury velocity $u^k$, which is approached by $\tcb{u_n}_{n\in\N}$ in the natural desingularization topology~\eqref{NATURAL_DESING_TOPO}, we have $\tp{I - e_\theta\otimes e_\theta}\tp{\grad\times u^k} \equiv 0$.

The arbitrary specification of the helicity and total meridional vorticity sequences $\tcb{h_n}_{n\in\N}$ and $\tcb{m_n}_{n\in\N}$ satisfying~\eqref{_CONSTRAINTS_ON_HELI_TOT_MER_VORT_SEQ_} along the swirling desingularization sequence $\tcb{u_n}_{n\in\N}$ is achieved through an admittedly pathological choice of smooth circulation functions $\tcb{S_n}_{n\in\N}$. More precisely, the desingularization sequence is again obtained via a sequence of Stokes stream functions solving Hicks' equation; we have
\begin{multline}\label{HICKS_FOR_THEOREM_2}
  -L\psi_n = r^2\Gamma_n\tp{\psi_n - r^2 - k_n} + \tp{S_nS_n'}\tp{\psi_n - r^2 - k_n}\text{ in }\Pi\\\text{and }u_n = -\f{\pd_z\psi_n}{r}e_r + \f{\pd_r\psi_n}{r}e_z + \f{S_n\tp{\psi_n - r^2 - k_n}}{r}e_\theta
\end{multline}
where $\Gamma_n$ and $k_n$ are as in~\eqref{form_of_the_smoothed_out_heaviside_and_flux}. Roughly speaking, the functions $\tcb{S_n}_{n\in\N}$ are selected to be families of increasingly oscillatory `combs' with decreasing amplitudes that assign circulation values to stream surfaces upholding the conspiracy of the generic prescription~\eqref{_THANKS_FOR_LOOKING_AT_THIS_LABEL_} without disrupting the limit~\eqref{NATURAL_DESING_TOPO}; in fact $S_n\to0$ in the velocity-relevant log-Lipschitz and $C^\alpha$ norms and $S_nS_n'\to0$, even though the total variation carried by $S_n'$ can be prescribed to be arbitrarily large. We emphasize that this mechanism takes place on interior stream surfaces lying in a fixed compact region separated both from the symmetry axis and from the vortex boundary, as suggested by~\eqref{_LOCATION_OF_THE_SWIRL_AWAY_FROM_RING_BOUNDARY_}. For more precise information, see Section~\ref{SS_swirl_operator}.

\subsection{Key elements of the proof}\label{SS_on_elements_of_the_proof}

We conclude the introduction with an outline of the proof of Theorem~\ref{main_thm_1_swirl_free}. The conceptual starting point is the Chandrasekhar transformation, which converts the axisymmetric Hicks equation~\eqref{Hicks_equation_introduction_form} into a semilinear elliptic equation on $\R^5$ with the axisymmetric Laplacian $L$ on $\Pi$ mapping to the Cartesian Laplacian $\Delta$ on $\R^5$. See Section~\ref{SS_chandra_trans_Euler_recon} for more details on this transformation.

In this $\R^5$ formulation, Hill's spherical vortex becomes the purely radial function $\v$, written explicitly in equation~\eqref{Hills_Vortex} and characterized by the Newtonian-potential identity
\begin{equation}\label{_HILL_R5_II_}
  \v(x) = \f{15}{8\pi^2}\int_{\R^5}\f{1}{\tabs{x - y}^3}\mathds{1}_{\tp{0,\infty}}\tp{\v\tp{y} - 1}\;\m{d}y\text{ for }x\in\R^5.
\end{equation}
Replacing the Heaviside function $\mathds{1}_{\tp{0,\infty}}$ by the smooth transition profiles $\Gamma_\ep$, which are defined in~\eqref{defn_regularized_heaviside}, writing the unknown perturbatively as $\v+\phi$, and subtracting off~\eqref{_HILL_R5_II_} leads one to the nonlinear integral equation 
\begin{equation}\label{_NONLINEAR_INTEGRAL_EQUATION_}
  \phi\tp{x} = \f{15}{8\pi^2}\int_{\R^5}\f{1}{\tabs{x - y}^3}\tp{\Gamma_\ep\tp{\tabs{\P y}^2\tp{\v(y) + \phi(y) - 1} - k} - \mathds{1}_{\tp{0,\infty}}\tp{\v\tp{y} - 1}}\;\m{d}y\text{ for }x\in\R^5,
\end{equation}
where $\P:\R^5\to\R^5$ is $\P y = \tp{y_1,y_2,y_3,y_4,0}$. We then write~\eqref{_NONLINEAR_INTEGRAL_EQUATION_} operator-theoretically as
\begin{equation}
  F_\ep\tp{\phi} = N_{\ep,k}\tp{\phi} + f_\ep.
\end{equation}
Here the principal operator is $F_\ep$ (defined in~\eqref{defn_principal_part_operator}), the perturbative operator is $N_{\ep,k}$ (defined in~\eqref{defn_pert_op}) and it exclusively carries the flux dependence, and the small source term is $f_\ep$ (defined in~\eqref{_defn_of_the_source_term_}), which arises from the regularization of the vorticity law.

Perhaps the most central analytic difficulty is to invert the principal part nonlinear operator $F_\ep$ on an $\ep$-uniform neighborhood of the origin in $C^1_{\m{b}}$ with estimates that remain useful in the singular limit $\ep\to0$. As it turns out, this is far from a routine application of the implicit function theorem: the Fr\'echet differential of $F_\ep$ involves the derivative of $\Gamma_\ep$, which has size comparable to $1/\ep$ and is concentrating near the four-dimensional hypersurfaces $\tcb{\v+\phi=1}$.

We first identify the singular limit as $\ep\to0$ of the linearized operators by utilizing stable coordinates on these moving level sets. The limiting formal linearization is an integral operator on $\S^4$, which, somewhat miraculously, is explicitly diagonalized by spherical harmonics. Moreover, its only kernel is the first spherical harmonic space, corresponding to the translation modes. The imposed cylindrical and reflection symmetries pleasantly eliminate this kernel. A delicate compactness argument, crucially using the precise uniform level set estimates of Appendix~\ref{appendix_estimates_on_level_sets}, then promotes this limiting nondegeneracy to a quantitative closed-range estimate uniform in $\ep$, but in the weaker $C^0_{\m{b}}$ norm. This, in turn, combines with a Newton-type quantitative local surjectivity result (Theorem 3.13 in~\cite{VPD}), thus yielding inverse maps for $F_\ep$ on an $\ep$-independent neighborhood in $C^1_{\m{b}}$, even though the norms of their differentials are seemingly not uniformly bounded in $\ep$.

The remaining terms require a different strategy. In particular, $N_{\ep,k}$ contains the factor $|\P y|^2$, arising from the cylindrical geometry, and its degenerating behavior near the symmetry axis precludes a simple perturbative estimate. We decompose the relevant region dyadically according to $|\P y|$ and combine the coarea formula with the uniform density estimates for perturbed level hypersurfaces developed in Appendix~\ref{appendix_estimates_on_level_sets}. This procedure produces the not necessarily sharp, yet sufficient bound
\begin{equation}\label{__BB1__}
  \tnorm{N_{\ep,k}(\phi)}_{C^1_{\m{b}}}
  \lesssim (\ep+k)^{1/6}.
\end{equation}
A similar but simpler strategy applied to the source term gives us the estimate
\begin{equation}\label{__BB2__}
  \tnorm{f_\ep}_{C^1_{\m{b}}}\lesssim\ep\tp{1 + \tabs{\log\ep}}.
\end{equation}
After the principal part $F_\ep$ has been inverted, the full equation is therefore reduced to the fixed point equation:
\begin{equation}\label{__BB3__}
  \phi = \tp{F_\ep}^{-1}\tp{N_{\ep,k}\tp{\phi} + f_\ep}.
\end{equation}
The inverse estimate on $F_\ep$ when combined with~\eqref{__BB1__} and~\eqref{__BB2__} shows the right hand side of~\eqref{__BB3__} to induce a self-map on a closed ball in the function space. This self-map is also compact. Hence an invocation of Schauder's fixed point theorem produces the sought-after smooth solutions. These solutions are initially defined on $B_{\R^5}[0,3] = \tcb{y\in\R^5\;:\;|y|\le 3}$; however, they admit unique global extensions after a simple argument utilizing the Newtonian potential. Additionally, we get uniform weighted decay, log-Lipschitz bounds, positivity, reflection symmetry, and strict half-space monotonicity in the traveling direction.

It remains to identify the singular limit. When the target flux is zero, the quantitative estimates directly force the constructed solutions to converge to Hill's vortex. For a fixed sufficiently small positive flux $k$, compactness instead gives a candidate limit, and one must show that the smooth head derivatives $15\Gamma_{\ep_n}$ converge in the equation to the discontinuous law $15\mathds{1}_{(0,\infty)}$. The strict monotonicity obtained above provides the mechanism needed to control the shrinking transition layers and pass to this limit in a weak formulation. The limiting profile is therefore shown to solve the $\R^5$-manifestation of the idealized $k$-Norbury problem~\eqref{_NON_DIM_NORBURY_WEAK_PROBLEM_}. The conditional uniqueness theorem of Amick and Fraenkel~\cite{MR918795}, together with our uniform a priori bounds, identifies this limiting profile with Norbury's original solution. Finally, the Chandrasekhar transformation, the $\R^5\to\R^3$ Euler reconstruction, and the velocity estimates of Section~\ref{SO_euler_recon_and_vortex_geometry} return the solutions to the familiar axisymmetric $\R^3$ formulation, yield the convergence asserted in Theorem~\ref{main_thm_1_swirl_free}, and give the stated geometry of the vortex cores.

\section{Nonlinear operators and mapping properties}\label{SO_nonlinear_ops}


Our task in this section is to formulate the construction of smooth vortex rings near Hill's solution in terms of a nonlinear integral equation on $\R^5$. The resulting operator naturally decomposes into a singular principal part, a perturbative part associated with the Norbury parameter, and a collection of swirl operators. We introduce these objects in turn and establish the mapping estimates that will be used in the inversion and existence arguments below.

\subsection{Operator framework}\label{SS_operator_frame}

We now introduce our main cast of operators and functions. The $\R^5$ manifestation of Hill's spherical vortex, originally discovered by Hill~\cite{Hill1894} in the intimately related $\R^3$ formulation discussed in Sections~\ref{SS_sol_HN_and_nondim} and~\ref{SS_chandra_trans_Euler_recon}, is the function
\begin{equation}\label{Hills_Vortex}
  \v\in C^{1,1}_\loc\tp{\R^5}\text{ defined via }\v(y) = \begin{cases}
    \f52 - \f32\tabs{y}^2&\text{for }|y|\le1,\\
    \tabs{y}^{-3}&\text{for }|y|\ge 1.
  \end{cases}
\end{equation}
The subsequent lemma shows that Hill's vortex~\eqref{Hills_Vortex} is a solution to the integral equation~\eqref{Hills_Integral_Equation} in which the vorticity function is the Heaviside function. Our interest lies in constructing nearby solutions with smoothed out vorticity functions. As such, we introduce some relevant notation. Let us fix
\begin{equation}\label{reg_vort_func}
  \gam\in C^\infty\tp{\R;[0,1]}\text{ satisfying }\gam(t) = 0\text{ for }t\le0,\;\gam(t)=1\text{ for }t\ge1,\text{ and }\gam'(t)>0\text{ for }0<t<1.
\end{equation}
Then, for $0<\ep\le2$, we denote
\begin{equation}\label{defn_regularized_heaviside}
  \Gamma_\ep\in C^\infty\tp{\R}\text{ via }\Gamma_\ep\tp{t} = \gam\tp{t/\ep}\text{ for all }t\in\R.
\end{equation}
Finally, we shall frequently rely on the projection operator 
\begin{equation}\label{the_important_projection_operator}
  \P:\R^5\to\R^5\text{ defined as }\P = I - e_5\otimes e_5.
\end{equation}

\begin{lem}[Simple observations]\label{lem on simple observations}
  The following hold.
  \begin{enumerate}
    \item For all $x\in\R^5$ the equality
    \begin{equation}\label{Hills_Integral_Equation}
      \v(x) = \f{15}{8\pi^2}\int_{\R^5}\f{1}{\tabs{x - y}^3}\mathds{1}_{\tp{0,\infty}}\tp{\v\tp{y}-1}\;\m{d}\mathcal{H}^{5}\tp{y}
    \end{equation}
    holds, where $\v$ is defined in equation~\eqref{Hills_Vortex}.
    \item For all $0<\ep\le2$, $0\le k\le1$, and $\phi\in C^1_{\m{b}}\tp{B_{\R^5}[0,3]}$ satisfying $\tnorm{\phi}_{C^1_{\m{b}}\tp{B_{\R^5}[0,3]}}\le 1/2$, the inclusion
    \begin{equation}\label{A crucial inclusion}
      \tcb{y\in B_{\R^5}[0,3]\;:\;\Gamma_\ep\tp{\tabs{\P y}^2\tp{\v\tp{y} + \phi(y) - 1} - k} > 0}\subseteq B_{\R^5}[0,2]
    \end{equation}
    holds.
    \item For all $\phi\in C^1_{\m{b}}\tp{B_{\R^5}[0,3]}$ satisfying $\tnorm{\phi}_{C^1_{\m{b}}\tp{B_{\R^5}\tsb{0,3}}}\le1/54$ we have
    \begin{equation}\label{an_important_positivity}
      \min\tcb{\v\tp{y} + \phi\tp{y}\;:\;y\in B_{\R^5}\tsb{0,3}}\ge1/54.
    \end{equation}
  \end{enumerate}
\end{lem}
\begin{proof}
  From the definition~\eqref{Hills_Vortex}, we deduce 
  \begin{equation}
    \mathds{1}_{\tp{0,\infty}}\tp{\v\tp{y} - 1} = \mathds{1}_{B_{\R^5}(0,1)}\tp{y}
  \end{equation}
  for all $y\in\R^5$. Thus it is natural to define the function
  \begin{equation}\label{thing_that_should_be_Hills}
    \tilde{\v}\in C^{1}\tp{\R^5}\text{ via }\tilde{\v}\tp{x}  = \f{15}{8\pi^2}\int_{\R^5}\f{1}{\tabs{x - y}^3}\mathds{1}_{B_{\R^5}\tp{0,1}}\tp{y}\;\m{d}\mathcal{H}^5\tp{y}.
  \end{equation}
  The kernel $\f{1}{8\pi^2}\tabs{x-y}^{-3}$ is precisely the $\R^5$-Newtonian potential for $-\Delta$, while the factor $15$ accounts for the right hand side below. For any test function $f\in C^1_{\m{c}}\tp{\R^5}$, the identities
  \begin{equation}\label{The equality of Laplacians}
    \int_{\R^5}\grad\tilde{\v}\cdot\grad f = 15\int_{B(0,1)}f = \int_{\R^5}\grad\v\cdot\grad f.
  \end{equation}
  are then seen to hold.

  Identity~\eqref{The equality of Laplacians} shows that $\Delta\tp{\v - \tilde{\v}} = 0$. Then, by appealing to Liouville's theorem for harmonic functions (see, e.g., Theorem 2.16 in Folland~\cite{MR1357411}) and the observation that both $\v$ and $\tilde{\v}$ decay to zero at infinity, the equality $\v = \tilde{\v}$ is justified.

  Let us now prove the second item. Assume that $y\in B_{\R^5}[0,3]$. If $|y|\le 2$, there is nothing to show, so consider the case that $2<|y|\le 3$. Then, since $\tabs{\phi\tp{y}}\le1/2$, the inequality 
  \begin{equation}
    \v\tp{y} + \phi(y) - 1\le -7/8 + 1/2<0
  \end{equation}
  is true. In turn, 
  \begin{equation}
    |\P y|^2\tp{\v\tp{y} + \phi\tp{y} - 1}-k\le-k\le0\imp\Gamma_\ep\tp{\tabs{\P y}^2\tp{\v\tp{y} + \phi\tp{y} - 1} - k} = 0.
  \end{equation}
  Thus $y$ cannot belong to the set on the left hand side of~\eqref{A crucial inclusion}.

  The third item is a trivial consequence of the observation that $\v\tp{y}\ge1/27$ for all $y\in B_{\R^5}[0,3]$.
\end{proof}

We introduce a functional and operator framework aimed at perturbing Hill's solution~\eqref{Hills_Integral_Equation}.

\begin{defn}[Spaces and mappings]\label{defn on spaces and mappings}
  For $0<r\le1/2$, $0<\ep\le 1$, and $0\le k\le 1$ we make the following definitions.
  \begin{enumerate}
    \item The space of orthogonal transformations of $\R^5$ leaving the $e_5$-axis invariant is denoted by
    \begin{equation}\label{the_special_subset_of_orthogonal_group}
      \O = \tcb{T\in\mathcal{L}\tp{\R^5}\;:\;T^{\m{t}}T = TT^{\m{t}} = I\text{ and }\tp{Te_5 = e_5\text{ or }Te_5 = - e_5}}.
    \end{equation}
    \item Let $q\in\tcb{0,1}$. The space of $C^q$ functions in $B_{\R^5}[0,3]$ that are cylindrically symmetric and $x_5$-even is
    \begin{equation}
      \mathcal{B}^q = \tcb{\phi\in C^q_{\m{b}}\tp{B_{\R^5}[0,3]}\;:\;\forall\;T\in\O,\;\forall\;x\in B_{\R^5}[0,3],\;\phi(Tx) = \phi(x)},\;\tnorm{\cdot}_{\mathcal{B}^q} = \tnorm{\cdot}_{C^q_{\m{b}}\tp{B(0,3)}}.
    \end{equation}
    The origin centered balls in $\mathcal{B}^q$ are denoted simply as
    \begin{equation}
      \mathcal{B}^q[r] = \tcb{\phi\in\mathcal{B}^q\;:\;\tnorm{\phi}_{\mathcal{B}^q}\le r}.
    \end{equation}
    When $q = 1$ we shall write $\mathcal{B} = \mathcal{B}^1$.
    \item The principal operator is the map $K_\ep:\mathcal{B}[1/2]\to\mathcal{B}$ with the action
    \begin{equation}\label{defn_principal_operator}
      K_\ep\tp{\phi} = -\f{15}{8\pi^2}\int_{B_{\R^5}[0,3]}\f{1}{\tabs{\cdot - y}^3}\tp{\Gamma_\ep\tp{\v\tp{y} + \phi\tp{y} - 1} - \Gamma_\ep\tp{\v\tp{y} - 1}}\;\m{d}\mathcal{H}^5\tp{y}\text{ for }\phi\in\mathcal{B}[1/2],
    \end{equation}
    where we recall that $\v$ is Hill's solution defined in equation~\eqref{Hills_Vortex} and the regularized Heaviside functions $\Gamma_\ep$ are given in equation~\eqref{defn_regularized_heaviside}.
    \item The perturbative operator is the map $N_{\ep,k}:\mathcal{B}[1/2]\to\mathcal{B}$ with the action
    \begin{equation}\label{defn_pert_op}
      N_{\ep,k}\tp{\phi} = \f{15}{8\pi^2}\int_{B_{\R^5}[0,3]}\f{1}{\tabs{\cdot - y}^3}\tp{\Gamma_\ep\tp{\tabs{\P y}^2\tp{\v\tp{y} + \phi\tp{y} - 1} - k} - \Gamma_\ep\tp{\v\tp{y} + \phi\tp{y} - 1}}\;\m{d}\mathcal{H}^5\tp{y},
    \end{equation}
    again for $\phi\in\mathcal{B}[1/2]$, where we recall that $\P$ is the projection operator defined in~\eqref{the_important_projection_operator}.
    \item The principal part operator is the map $F_\ep:\mathcal{B}[1/2]\to\mathcal{B}$ with the action
    \begin{equation}\label{defn_principal_part_operator}
      F_\ep\tp{\phi} = \phi + K_\ep\tp{\phi},\;\phi\in\mathcal{B}[1/2].
    \end{equation}
   \end{enumerate}
\end{defn}

\begin{rmk}\label{remark on symmetry preservation}
  The maps $K_\ep$ and $N_{\ep,k}$ defined in~\eqref{defn_principal_operator} and~\eqref{defn_pert_op}, respectively, indeed map into $\mathcal{B}$ in the sense that the symmetry of being invariant under $\O$ (defined in~\eqref{the_special_subset_of_orthogonal_group}) is preserved under these nonlinear transformations. To see this, one needs only observe that $\v$ (defined in~\eqref{Hills_Vortex}) is radial and that for all $T\in\O$ and $x\in\R^5$ we have $T\P x = \P Tx$ and $\tabs{T x} = \tabs{x}$.
\end{rmk}

\begin{prop}[Basic smoothness and compactness]\label{prop on basic smoothness}
  There exists a constant $C\in\R^+$ such that for all $0<\ep\le 1$ and all $0\le k\le 1$ the operators $K_\ep,N_{\ep,k}:\mathcal{B}[1/2]\to\mathcal{B}$ are well-defined, smooth in the Fr\'echet sense, compact, and the following hold for all $\phi,\tilde{\phi}\in\mathcal{B}[1/2]$.
  \begin{enumerate}
    \item We have the inclusions $K_\ep\tp{\phi},N_{\ep,k}\tp{\phi}\in\m{LL}^1\tp{B_{\R^5}[0,3]}$ and the bound
    \begin{equation}\label{uniform LL1 esimate}
      \max\tcb{\tnorm{K_\ep\tp{\phi}}_{\m{LL}^1},\tnorm{N_{\ep,k}\tp{\phi}}_{\m{LL}^1}}\le C
    \end{equation}
    holds.
    \item The Fr\'echet differential $DK_\ep\tp{\phi}\in\mathcal{L}\tp{\mathcal{B}}$ extends to a linear mapping
    \begin{equation}\label{the_extended_frechet_differential}
      A_\ep\tp{\phi}\varphi = -\f{15}{8\pi^2}\int_{B_{\R^5}[0,3]}\f{1}{\tabs{\cdot - y}^3}\tp{\Gamma_\ep}'\tp{\v\tp{y} + \phi\tp{y} - 1}\varphi\tp{y}\;\m{d}\mathcal{H}^5\tp{y}.
    \end{equation}
    This extension acts between the spaces
    \begin{equation}
      A_\ep\tp{\phi}:C^0_{\m{b}}\tp{B_{\R^5}[0,3]}\to C^1_{\m{b}}\tp{B_{\R^5}[0,3]}\text{ and }A_\ep\tp{\phi}:C^1_{\m{b}}\tp{B_{\R^5}[0,3]}\to C^2_{\m{b}}\tp{B_{\R^5}[0,3]}
    \end{equation}
    and the bounds 
    \begin{equation}\label{the_bounds_right_here}
      \tnorm{A_\ep\tp{\phi}\varphi}_{C^1}\le\tp{C/\ep}\tnorm{\varphi}_{C^0},\;\tnorm{A_\ep\tp{\phi}\varphi}_{C^2}\le\tp{C/\ep^2}\tnorm{\varphi}_{C^1},\text{ and}\;\tnorm{\tp{A_\ep\tp{\phi} - A_\ep\tp{\tilde{\phi}}}\varphi}_{C^1}\le\tp{C/\ep^2}\tnorm{\phi - \tilde{\phi}}_{C^0}\tnorm{\varphi}_{C^0}
    \end{equation}
    hold.
  \end{enumerate}
\end{prop}
\begin{proof}
  The symmetry assertion required for the maps $K_\ep,N_{\ep,k}:\mathcal{B}[1/2]\to\mathcal{B}$ is the content of Remark~\ref{remark on symmetry preservation}, while their well-definedness as mappings into $C^1$ follows from the standard local estimates for Newtonian potentials of bounded, compactly supported functions. The verification of Fr\'echet smoothness is also a simple observation; these nonlinear operators are structurally compositions of a bounded linear operator on $C^1_{\m{b}}\tp{B_{\R^5}[0,3]}$ (convolution with the Newtonian potential) and outer composition with the smooth function $\Gamma_\ep$. Therefore, Fr\'echet smoothness is now a routine computation by using the converse to Taylor's theorem (see, e.g., Section 2.4B in Abraham, Marsden, and Ratiu~\cite{MR960687}). 

  The estimate~\eqref{uniform LL1 esimate} claimed within the first item, in particular, establishes that the operators $K_\ep$ and $N_{\ep,k}$ are compact, thanks to the Arzel\`a-Ascoli theorem.

  The establishment of~\eqref{uniform LL1 esimate} for $N_{\ep,k}$ is much the same as it is for $K_\ep$; thus, we shall exclusively focus on the latter case. It is plain from Young's inequality that 
  \begin{equation}
    \tnorm{K_\ep\tp{\phi}}_{C^1_{\m{b}}\tp{B[0,3]}}\le C\sup_{x\in B_{\R^5}[0,3]}\int_{B_{\R^5}[0,3]}\bp{\f{1}{\tabs{x - y}^3} + \f{1}{\tabs{x - y}^4}}\;\m{d}\mathcal{H}^5\tp{y}\le C.
  \end{equation}

  To estimate the log-Lipschitz norm on the gradient, we initially write for any $x,\tilde{x}\in B_{\R^5}[0,3]$.
  \begin{equation}\label{cited_suprisingly}
    \tabs{\grad K_\ep\tp{\phi}\tp{x} - \grad K_\ep\tp{\phi}\tp{\tilde{x}}}\le C\int_{B_{\R^5}[0,3]}\babs{\f{x - y}{\tabs{x - y}^5} - \f{\tilde{x} - y}{\tabs{\tilde{x} - y}^5}}\;\m{d}\mathcal{H}^5\tp{y}.
  \end{equation}
  Now we split the region of integration according to 
  \begin{equation}
    D(x,\tilde{x}) = \tcb{y\in B_{\R^5}[0,3]\;:\;\tabs{y - \tp{x + \tilde{x}}/2}\le 2\tabs{x - \tilde{x}}}.
  \end{equation}
  For $y\in D(x,\tilde{x})$ we do not exploit the difference in~\eqref{cited_suprisingly}; rather, we use the inclusion 
  \begin{equation}
    D\tp{x,\tilde{x}}\subset B(x,5\tabs{x - \tilde{x}}/2)\cap B(\tilde{x},5\tabs{x - \tilde{x}}/2)
  \end{equation}
  and radial-spherical coordinates to bound
  \begin{equation}
    \int_{B_{\R^5}[0,3]\cap D(x,\tilde{x})}\f{1}{\tabs{z - y}^4}\;\m{d}\mathcal{H}^{5}\tp{y}\le C\tabs{x - \tilde{x}}
  \end{equation}
  for $z\in\tcb{x,\tilde{x}}$.

  On the other hand, for $y\in B_{\R^5}[0,3]\setminus D(x,\tilde{x})$, we first employ the fundamental theorem of calculus to get the bound
  \begin{eqnarray}
    \babs{\f{x - y}{\tabs{x - y}^5} - \f{\tilde{x} - y}{\tabs{\tilde{x} - y}^5}}\le C\f{\tabs{x - \tilde{x}}}{\tabs{y - \tp{x + \tilde{x}}/2}^5}.
  \end{eqnarray}
  Then, we use radial-spherical coordinates again to estimate
  \begin{equation}
    \int_{B_{\R^5}[0,3]\setminus D(x,\tilde{x})}\f{1}{\tabs{y - \tp{x + \tilde{x}}/2}^5}\;\m{d}\mathcal{H}^5\tp{y}\le C\tp{1 + \tabs{\log\tabs{x - \tilde{x}}}}.
  \end{equation}
  Synthesizing the above casework shows that 
  \begin{equation}
    \tabs{\grad K_\ep\tp{\phi}\tp{x} - \grad K_\ep\tp{\phi}\tp{\tilde{x}}}\le C\tabs{x - \tilde{x}}\tp{1 + \tabs{\log\tabs{x - \tilde{x}}}}.
  \end{equation}
  The validity of inequality~\eqref{uniform LL1 esimate} is now apparent.

  We focus now on the second item; it is sufficient for us to only prove the bounds~\eqref{the_bounds_right_here}. The first of these is essentially a direct consequence of the pointwise estimate 
  \begin{equation}
    \tabs{\grad\tp{A_\ep\tp{\phi}\varphi}\tp{x}}\le\f{C}{\ep}\int_{B_{\R^5}[0,3]}\f{\tabs{\varphi\tp{y}}}{\tabs{x - y}^4}\;\m{d}\mathcal{H}^5\tp{y}\le\f{C}{\ep}\tnorm{\varphi}_{C^0}\text{ for }x\in B_{\R^5}[0,3]
  \end{equation}
  and Young's inequality. 

  The second bound asserted in~\eqref{the_bounds_right_here} requires integration by parts in order to stay away from a kernel singularity strength agreeing with the spatial dimension. The pointwise estimate
  \begin{equation}
    \tabs{\grad^2\tp{A_\ep\tp{\phi}\varphi}\tp{x}}\le\f{C}{\ep^2}\int_{B_{\R^5}[0,3]}\f{\tabs{\varphi(y)} + \tabs{\grad\varphi\tp{y}}}{\tabs{x - y}^4}\;\m{d}\mathcal{H}^5\tp{y}
  \end{equation}
  holds.

  The final estimate of~\eqref{the_bounds_right_here} is a consequence of the fundamental theorem of calculus and the pointwise bound
  \begin{multline}
    \tabs{\grad\tp{A_\ep\tp{\phi}\varphi - A_\ep\tp{\tilde{\phi}}\varphi}\tp{x}}\\\le C\int_{B_{\R^5}[0,3]}\f{1}{\tabs{x - y}^4}\int_0^1\tabs{\tp{\Gamma_\ep}''\tp{\v\tp{y} - 1 + \tau\phi\tp{y} + \tp{1 - \tau}\tilde{\phi}\tp{y}}}\;\m{d}\tau\tabs{\phi(y) - \tilde{\phi}\tp{y}}\tabs{\varphi\tp{y}}\;\m{d}\mathcal{H}^5\tp{y}.
  \end{multline}
\end{proof}

To conclude this subsection, we initiate the study of the vanishing limit of parameters in the maps $A_\ep\tp{\phi}$ that are defined in equation~\eqref{the_extended_frechet_differential}. The stable local level set coordinates constructed in Appendix~\ref{SS_level_set_coordinates} are an essential part of this analysis.

\begin{prop}[Limit Identification]\label{prop on limit identification}
  Let $\tcb{\ep_n}_{n\in\N}\subset(0,1]$, $\tcb{\phi_n}_{n\in\N}\subset\mathcal{B}[1/2]$, $\varphi\in C^0_{\m{b}}\tp{B_{\R^5}[0,3]}$, and $\tcb{\varphi_n}_{n\in\N}\subset C^0_{\m{b}}\tp{B_{\R^5}[0,3]}$ be such that the limits 
  \begin{equation}\label{convergence hypotheses}
    \ep_n\to0,\;\tnorm{\phi_n}_{\mathcal{B}}\to0,\;\tnorm{\varphi_n - \varphi}_{C^0}\to0\text{ as }n\to\infty.
  \end{equation}
  hold. Then 
  \begin{equation}\label{the limit that we would like to show}
    \int_{B_{\R^5}[0,3]}\tp{\Gamma_{\ep_n}}'\tp{\v\tp{y} + \phi_n\tp{y} - 1}\varphi_n\tp{y}\;\m{d}\mathcal{H}^5\tp{y} \to \f13\int_{\S^4}\varphi\tp{y}\;\m{d}\mathcal{H}^4\tp{y}\text{ as }n\to\infty,
  \end{equation}
  where $\S^4$ denotes the unit sphere in $\R^5$.
\end{prop}
\begin{proof}
  We begin by making a pair of observations. First:
  \begin{equation}\label{observation_1_}
    \tabs{\grad\v\tp{y}} = 3\text{ for all }y\in\S^4,\;\v\tp{y} - 1>0\text{ for all }|y|<1,\;\v\tp{y} -1 < 0\text{ for all }|y|>1.
  \end{equation}
  Here we recall that Hill's solution $\v$ is defined in equation~\eqref{Hills_Vortex}. Second: for any $n\in\N$
  \begin{equation}\label{observation_2_}
    \tcb{y\in B_{\R^5}[0,3]\;:\;\tp{\Gamma_{\ep_n}}'\tp{\v\tp{y} + \phi_n\tp{y} - 1}\neq 0}\subseteq\tcb{y\in B_{\R^5}[0,3]\;:\;\tabs{\v\tp{y} - 1}\le\ep_n + \tnorm{\phi_n}_{C^0}}.
  \end{equation}
  Together observations~\eqref{observation_1_} and~\eqref{observation_2_} imply that as $n\to\infty$ the support of the integrand on the left hand side of equation~\eqref{the limit that we would like to show} inhabits a vanishingly thin tubular neighborhood of $\S^4$. Since $\tnorm{\phi_n}_{C^1}\to0$, we can ensure that eventually $|\grad\tp{\v + \phi_n}|\ge 2$ on these tubular neighborhoods of support.

  Therefore, for sufficiently large $n\in\N$, we can use the coarea formula, see, e.g., Theorem 3.2.12 in Federer~\cite{MR257325}, to equate
  \begin{multline}\label{the_first_man}
    \int_{B_{\R^5}[0,3]}\tp{\Gamma_{\ep_n}}'\tp{\v\tp{y} + \phi_n\tp{y} - 1}\varphi_n\tp{y}\;\m{d}\mathcal{H}^5\tp{y} \\
    = \f{1}{\ep_n}\int_{0}^{\ep_n}\gam'\tp{t/\ep_n}\int_{B_{\R^5}[0,3]}\mathds{1}_{\tcb{\v + \phi_n - 1 = t}}\tp{y}\f{\varphi_n\tp{y}}{\tabs{\grad\tp{\v + \phi_n}\tp{y}}}\;\m{d}\mathcal{H}^4\tp{y}\;\m{d}\mathcal{H}^1\tp{t}.
  \end{multline}

  To proceed further, we need to break into the stable local level set coordinates of Theorem~\ref{thm on stable local level set coordinates}. Invoking this result grants us small positive parameters $r_0,\ep_0>0$, a constant $C_0\in\R^+$, and $C^1$ functions
  \begin{equation}
    \Theta_\phi:\S^4\times[-\ep_0,\ep_0]\to B_{\R^5}[0,2]\text{ defined for }\tnorm{\phi}_{C^1_{\m{b}}}\le r_0
    \end{equation}
  such that for all $|t|\le\ep_0$ the map
  \begin{equation}\label{the level set diffeomorphism identity}
    \Theta_\phi\tp{\cdot,t}:\S^4\to\tcb{\v + \phi - 1 = t}
  \end{equation}
  is a diffeomorphism obeying the bi-Lipschitz estimate 
  \begin{equation}
    \tp{1/C_0}\tabs{m - \tilde{m}}\le\tabs{\Theta_\phi\tp{m,t} - \Theta_\phi\tp{\tilde{m},t}}\le C_0\tabs{m - \tilde{m}}.
  \end{equation}

  We use the maps~\eqref{the level set diffeomorphism identity} in conjunction with the area formula (see, e.g., Corollary 3.2.20 in Federer~\cite{MR257325}) followed by a change of variables in the one dimensional integral in order to rewrite (for sufficiently large $n\in\N$) identity~\eqref{the_first_man} as 
  \begin{multline}\label{intermediate_step_}
    \int_{B_{\R^5}[0,3]}\tp{\Gamma_{\ep_n}}'\tp{\v\tp{y} + \phi_n\tp{y} - 1}\varphi_n\tp{y}\;\m{d}\mathcal{H}^5\tp{y} \\
    =\int_{0}^{1}\gam'\tp{\tau}\int_{\S^4}\f{\varphi_n\circ\Theta_{\phi_n}\tp{y,\ep_n\tau}}{\tabs{\grad\tp{\v + \phi_n}\circ\Theta_{\phi_n}\tp{y,\ep_n\tau}}}J\Theta_{\phi_n}(\cdot,\ep_n\tau)\tp{y}\;\m{d}\mathcal{H}^4\tp{y}\;\m{d}\mathcal{H}^1\tp{\tau}.
  \end{multline}
  Here $J\Theta_{\phi}\tp{\cdot,t} = \tabs{\wedge_4D\Theta_\phi\tp{\cdot,t}}$ denotes the Jacobian of the diffeomorphism~\eqref{the level set diffeomorphism identity}.

  The fourth item of Theorem~\ref{thm on stable local level set coordinates} informs us that the map
  \begin{equation}\label{joint_cont}
    B_{C^1_{\m{b}}}[0,r_0]\ni\phi\mapsto\Theta_\phi\in C^1\tp{\S^4\times[-\ep_0,\ep_0];\R^5}
  \end{equation}
  is continuous. Thus we combine the dominated convergence theorem with the continuity of the map~\eqref{joint_cont} and the hypotheses~\eqref{convergence hypotheses} to justifiably pass to the limit $n\to\infty$ on the right hand side of~\eqref{intermediate_step_}. In doing so, we learn that
  \begin{equation}\label{_almost_there_}
    \lim_{n\to\infty}\int_{B_{\R^5}[0,3]}\tp{\Gamma_{\ep_n}}'\tp{\v\tp{y} + \phi_n\tp{y} - 1}\varphi_n\tp{y}\;\m{d}\mathcal{H}^5\tp{y} = \f13\int_{\S^4}\varphi\circ\Theta_0\tp{y,0}J\Theta_0\tp{\cdot,0}\tp{y}\;\m{d}\mathcal{H}^4\tp{y};
  \end{equation}
  Here we have used that the integral of $\gam'$ over $[0,1]$ is unity and the left hand equality in~\eqref{observation_1_}. Finally, the claimed limit of equation~\eqref{the limit that we would like to show} follows from~\eqref{_almost_there_} after another application of the area formula.
\end{proof}

\subsection{Principal operator}\label{SS_principal_operator}

We now analyze more sophisticated forward mapping properties of the map $K_\ep$, whose action we recall is given in the third item of Definition~\ref{defn on spaces and mappings}, and the related map $A_\ep$ from the second item of Proposition~\ref{prop on basic smoothness}.

\begin{prop}[Principal operator's log-Lipschitz estimates]\label{prop on principal operator lipschitz estimate}
  There exists $C,\ep_1\in\R^+$ with $\ep_1\le1/2$ such that for all $\phi,\tilde{\phi}\in\mathcal{B}[\ep_1]$, $0<\ep\le\ep_1$, and $\varphi\in C^0_{\m{b}}\tp{B_{\R^5}[0,3]}$ the following hold.
  \begin{enumerate}
    \item We have the estimate
    \begin{equation}\label{the_claimed_bound_here}
      \tnorm{A_\ep\tp{\phi}\varphi}_{\m{LL}}\le C\tnorm{\varphi}_{C^0}.
    \end{equation}
    \item We have the estimate
    \begin{equation}
      \tnorm{K_\ep\tp{\phi} - K_\ep\tp{\tilde{\phi}}}_{\m{LL}}\le C\tnorm{\phi - \tilde{\phi}}_{C^0}.
    \end{equation}
  \end{enumerate}
\end{prop}
\begin{proof}
  Initially we fix $0<\ep\le1/2$ and $\phi\in\mathcal{B}[1/2]$. From the definition of $\Gamma_\ep$ in equation~\eqref{defn_regularized_heaviside}, we see that $\supp\tp{\Gamma_\ep}'\subseteq[0,\ep]$. Hence, for $x,\tilde{x}\in B_{\R^5}[0,3]$ we deduce the pointwise estimates
  \begin{equation}\label{int_no_1}
    \tabs{\tp{A_\ep\tp{\phi}\varphi}\tp{x}}\le\f{C}{\ep}\tnorm{\varphi}_{C^0}\int_{\tcb{\tabs{\v + \phi - 1}\le\ep}}\f{\mathds{1}_{B_{\R^5}[0,3]}\tp{y}}{\tabs{x - y}^3}\;\m{d}\mathcal{H}^5\tp{y}
  \end{equation}
  and
  \begin{equation}\label{int_no_2}
    \tabs{\tp{A_\ep\tp{\phi}\varphi}\tp{x} - \tp{A_\ep\tp{\phi}\varphi}\tp{\tilde{x}}}\le\f{C}{\ep}\tnorm{\varphi}_{C^0}\int_{\tcb{\tabs{\v + \phi - 1}\le\ep}}\babs{\f{1}{\tabs{x - y}^3} - \f{1}{\tabs{\tilde{x} - y}^3}}\mathds{1}_{B_{\R^5}[0,3]}\tp{y}\;\m{d}\mathcal{H}^5\tp{y}.
  \end{equation}

  To proceed in further bounding of expressions~\eqref{int_no_1} and~\eqref{int_no_2}, we look to Corollary~\ref{bounds for integral operators on thin domains}, whose hypotheses are satisfied thanks to the properties of Hill's solution $\v$ enumerated in equation~\eqref{observation_1_}. Invoking this result gives us a parameter $0<\ep_1\le1/2$ such that whenever $0<\ep\le\ep_1$ and $\phi\in\mathcal{B}[\ep_1]$ we have $\tcb{\tabs{\v + \phi - 1}\le\ep}\Subset B_{\R^5}[0,3]$ and the estimates 
  \begin{multline}\label{int_no_3}
    \int_{\tcb{\tabs{\v + \phi - 1}\le\ep}}\f{1}{\tabs{x - y}^3}\;\m{d}\mathcal{H}^5\tp{y}\le C\ep,\\
    \int_{\tcb{\tabs{\v + \phi - 1}\le\ep}}\babs{\f{1}{\tabs{x - y}^3} - \f{1}{\tabs{\tilde{x} - y}^3}}\;\m{d}\mathcal{H}^5\tp{y}\le C\ep\tabs{x - \tilde{x}}\tp{1 + \tabs{\log\tabs{x - \tilde{x}}}}
  \end{multline} 
  for a constant $C$ independent of $x,\tilde{x}\in B_{\R^5}[0,3]$, $\ep$, and $\phi$. The combination of~\eqref{int_no_1}, \eqref{int_no_2}, and~\eqref{int_no_3} verifies the claimed bound~\eqref{the_claimed_bound_here}.

  The second item is a direct corollary of the first item and the fundamental theorem of calculus identity
  \begin{equation}\label{FTC_IDENTITY}
    K_\ep\tp{\phi} - K_\ep\tp{\tilde{\phi}} = \bp{\int_0^1A_\ep\tp{\tau\phi + \tp{1 - \tau}\tilde{\phi}}\;\m{d}\tau}\tp{\phi - \tilde{\phi}}.
  \end{equation}
  Note that the expression~\eqref{FTC_IDENTITY} is legitimate as a consequence of Proposition~\ref{prop on basic smoothness}.
\end{proof}

As a straightforward consequence of both Propositions~\ref{prop on basic smoothness} and~\ref{prop on principal operator lipschitz estimate}, we find a suitably sharp almost Lipschitz modulus of continuity for the map $K_\ep:\mathcal{B}[\ep_1]\to\mathcal{B}$.

\begin{coro}[Principal operator's uniform continuity estimate]\label{coro on principal log-lip}
  Let $0<\ep_1\le1/2$ be the small parameter from Proposition~\ref{prop on principal operator lipschitz estimate}. There exists a constant $C\in\R^+$ such that for all $\phi,\tilde{\phi}\in\mathcal{B}[\ep_1]$ and all $0<\ep\le\ep_1$ the bound 
  \begin{equation}\label{an_important_continuity_estimate}
    \tnorm{K_\ep\tp{\phi} - K_\ep\tp{\tilde{\phi}}}_{C^1}\le C\tnorm{\phi - \tilde{\phi}}_{C^0}\tp{1 + \tabs{\log\tnorm{\phi - \tilde{\phi}}_{C^0}}}.
  \end{equation}
  holds.
\end{coro}
\begin{proof}
  Let us recall here a standard interpolation estimate, for a proof see, e.g., Lemma A.1 in Radu and Stevenson~\cite{VPD}. In this context, we get that 
  \begin{equation}\label{the interpolation estimate}
    \tnorm{f}_{C^1}\le C\tnorm{f}_{\m{LL}}\bp{1 + \babs{\log\bp{\f{c\tnorm{f}_{\m{LL}^1}}{\tnorm{f}_{\m{LL}}}}}}
  \end{equation}
  for all $f\in\m{LL}^1\tp{B_{\R^5}[0,3]}$, where $C$ and $c$ are positive constants depending only on the domain and the dimension.

  The interpolation estimate~\eqref{the interpolation estimate} is relevant to us in the case that $f = K_\ep\tp{\phi} - K_\ep\tp{\tilde{\phi}}$ for some $0<\ep\le\ep_1$ and $\phi,\tilde{\phi}\in\mathcal{B}[\ep_1]$. Due to Propositions~\ref{prop on basic smoothness} and~\ref{prop on principal operator lipschitz estimate}, the estimates 
  \begin{equation}\label{the_devil_you_know}
    \tnorm{f}_{\m{LL}^1}\le C\text{ and }\tnorm{f}_{\m{LL}}\le C\tnorm{\phi - \tilde{\phi}}_{C^0}
  \end{equation}
  are true for some constant $C$ independent of $\ep$, $\phi$, and $\tilde{\phi}$. The desired logarithmic-Lipschitz estimate~\eqref{an_important_continuity_estimate} then follows immediately from~\eqref{the interpolation estimate} and~\eqref{the_devil_you_know}.
\end{proof}

\subsection{Perturbative operator}\label{SS_peturbative_operator}

We turn to focus on more subtle mapping estimates for the perturbative operator $N_{\ep,k}$ that we remind the reader is defined in the fourth item of Definition~\ref{defn on spaces and mappings}. We begin with a lemma that decomposes the ball $B_{\R^5}[0,3]$ into countably many regions of different behaviors.

\begin{lem}[A decomposition and its properties]\label{lem on a decomposition and properties}
  We define the sequence of sets
  \begin{multline}\label{the sets of the decomposition for N}
    E_0 = \tcb{x\in B_{\R^5}[0,3]\;:\;|x|\le1/2\text{ or }3/2\le|x|}\\\text{ and for }n\in\N^+,\;E_n = \tcb{x\in B_{\R^5}[0,3]\;:\;1/2\le|x|\le3/2,\;2^{-n+1}\le\tabs{\P x}\le 2^{-n+2}},
  \end{multline}
  where $\P$ is the projection operator defined in~\eqref{the_important_projection_operator}. There exist constants $c,\ep_2\in\R^+$ with $\ep_2\le\ep_1$, where $\ep_1$ is from Proposition~\ref{prop on principal operator lipschitz estimate}, such that the following hold for all $\phi\in\mathcal{B}[\ep_2]$ and $t\in[-4/5,6/5]$.
  \begin{enumerate}
    \item We have 
    \begin{equation}
      \forall\;n\in\N,\;E_n\subset B_{\R^5}[0,3]\text{ and }\mathcal{H}^5\bp{B_{\R^5}[0,3]\setminus\bigcup_{n\in\N}E_n} = 0.
    \end{equation}
     \item For all $y\in E_0$ we have that 
     \begin{equation}\label{the_good_estimates_}
      y\in B_{\R^5}[0,1/2]\imp\v\tp{y} + \phi\tp{y} - 1\ge1/c\text{ and }y\not\in B_{\R^5}(0,3/2)\imp\v\tp{y} + \phi(y) - 1\le-1/c. 
     \end{equation}
     \item We have the inclusions
     \begin{equation}\label{a_good_covering}
      \bigcup_{n\ge 1}E_n\subset\tcb{y\in B_{\R^5}[0,3]\;:\;\v\tp{y} + \phi\tp{y} - 1\in[-4/5,6/5]}\Subset B_{\R^5}\tp{0,3}
     \end{equation}
     and the estimates
     \begin{equation}\label{some_good_estimates}
      \min_{\tcb{-4/5\le\v + \phi - 1\le6/5}}\tabs{\grad\tp{\v + \phi}}\ge1/c,\;\mathcal{H}^4\tp{E_n\cap \tcb{\v + \phi - 1 = t}}\le c2^{-4n}
     \end{equation}
     hold.
  \end{enumerate}
\end{lem}
\begin{proof}
  The first item is trivial.

  From the properties of $\v$ defined in~\eqref{Hills_Vortex}, the set 
  \begin{equation}
    \tilde{E} = \tcb{y\in B_{\R^5}[0,3]\;:\;\v\tp{y} - 1\in[-19/27,9/8]}
  \end{equation}
  satisfies
  \begin{equation}\label{properties_of_the_set_tilde_E}
    \tilde{E} = B_{\R^5}[0,3/2]\setminus B_{\R^5}\tp{0,1/2}\supset\bigcup_{n\ge1}E_n.
  \end{equation}
  For all $\phi\in\mathcal{B}[3/40]$, 
  \begin{equation}
    \tilde{E}\subseteq\tcb{\v + \phi - 1\in[-4/5,6/5]}\Subset B_{\R^5}(0,3).
  \end{equation}
  Moreover, we can estimate
  \begin{equation}
    \min_{\tcb{-4/5\le\v + \phi - 1\le6/5}}\tabs{\grad\tp{\v + \phi}}\ge9/80\text{ and }\min_{E_0}\tabs{\v + \phi - 1}\ge679/1080.
  \end{equation}
  Thus, the second item holds so long as we select $\ep_2\le3/40$.

  For the third item, the above computations show that~\eqref{a_good_covering} and the left hand side of~\eqref{some_good_estimates} are true, again for $\ep_2\le3/40$. The only nontrivial assertion is the right hand inequality of~\eqref{some_good_estimates}. For this, we invoke the density estimates of Proposition~\ref{prop on transverse cylinder density estimate} to acquire $0<\ep_2\le\min\tcb{3/40,\ep_1}$ and a suitably uniform constant $c>0$ such that not only the remaining estimate of~\eqref{some_good_estimates} holds, but also all of the assertions of the lemma.
\end{proof}

The decomposition of $B_{\R^5}[0,3]$ into the regions $\tcb{E_n}_{n\in\N}$ determined by Lemma~\ref{lem on a decomposition and properties} is central to the Lebesgue space estimates in the next result.

\begin{prop}[Lebesgue space estimates of the perturbative operator]\label{prop on Lebesgue space estimates on the perturbative operator}
  Let $\ep_2>0$ be the small parameter granted by Lemma~\ref{lem on a decomposition and properties}. There exists a constant $C\in\R^+$ such that for all $1\le p\le\infty$, $0<\ep\le1$, $0\le k\le 1$, and $\phi\in\mathcal{B}[\ep_2]$ the function
  \begin{equation}\label{defn_of_the_density}
    \mathscr{N}_{\ep,k}\tp{\phi}\tp{y} = \Gamma_\ep\tp{\tabs{\P y}^2\tp{\v\tp{y} + \phi\tp{y} - 1} - k} - \Gamma_\ep\tp{\v\tp{y} + \phi\tp{y} - 1},\text{ defined for }y\in B_{\R^5}[0,3],
  \end{equation}
  obeys the bound
  \begin{equation}\label{lebesgue_bounds_WTS}
    \tnorm{\mathscr{N}_{\ep,k}\tp{\phi}}_{L^p\tp{B[0,3]}}\le C\tp{\ep + k}^{1/p}.
  \end{equation}
\end{prop}
\begin{proof}
  For the case $p = \infty$ the claimed bound~\eqref{lebesgue_bounds_WTS} is immediate since $\tabs{\mathscr{N}_{\ep,k}\tp{\phi}\tp{y}}\le 1$ for all $y\in B_{\R^5}[0,3]$. The intermediate cases of $1<p<\infty$ will then of course follow by standard interpolation properties of the Lebesgue space norms, so long as we establish~\eqref{lebesgue_bounds_WTS} in the case $p=1$. This task is the focus of the remainder of the proof.

  The first item of Lemma~\ref{lem on a decomposition and properties} allows us to decompose the region of integration and so the initial estimate
  \begin{equation}\label{series_decomposition}
    \tnorm{\mathscr{N}_{\ep,k}\tp{\phi}}_{L^1\tp{B[0,3]}}\le\sum_{n=0}^\infty\tnorm{\mathscr{N}_{\ep,k}\tp{\phi}}_{L^1\tp{E_n}}
  \end{equation}
  holds. The strategy for the $n=0$ term in the series~\eqref{series_decomposition} is slightly different from the cases $n\ge1$.

  Let $c\in\R^+$ denote the constant of Lemma~\ref{lem on a decomposition and properties} such that the inequalities~\eqref{the_good_estimates_} hold. After enlarging $c$, if necessary, we may assume that $c\ge1$. If $1/c<\ep\le1$, then the pointwise bound $\tabs{\mathscr{N}_{\ep,k}\tp{\phi}}\le1$ gives $\tnorm{\mathscr{N}_{\ep,k}\tp{\phi}}_{L^1\tp{B[0,3]}}\le\mathcal{H}^5\tp{B_{\R^5}[0,3]}\le C\tp{\ep+k}$. We may therefore assume for the remainder of the proof that $0<\ep\le1/c$. Note that if $y\in E_0$, $0<\ep\le1/c$, and $\phi\in\mathcal{B}[\ep_2]$, we have 
  \begin{equation}
    |y|\le1/2\imp\Gamma_\ep\tp{\v\tp{y} + \phi\tp{y} - 1} = 1\text{ and }|y|\ge3/2\imp\Gamma_\ep\tp{\v\tp{y} + \phi\tp{y} - 1} = 0.
  \end{equation}
  In turn, we deduce that 
  \begin{equation}\label{_L1_}
    \tnorm{\mathscr{N}_{\ep,k}\tp{\phi}}_{L^1\tp{E_0}}\le\mathcal{H}^5\tp{E_0^-\tp{\phi,\ep,k}} \\+ \mathcal{H}^5\tp{E_0^+\tp{\phi,\ep,k}}
  \end{equation}
  for the sets
  \begin{multline}\label{_L2_}
    E_0^-\tp{\phi,\ep,k} = \tcb{y\in E_0\;:\;|y|\le1/2,\;\Gamma_\ep\tp{\tabs{\P y}^2\tp{\v\tp{y} + \phi\tp{y} - 1} - k}<1}\\\text{and }E_0^+\tp{\phi,\ep,k} = \tcb{y\in E_0\;:\;|y|\ge3/2,\;\Gamma_\ep\tp{\tabs{\P y}^2\tp{\v\tp{y} + \phi\tp{y} - 1} - k}>0}.
  \end{multline}
  Then, a simple argument based on the inequalities~\eqref{the_good_estimates_} and casework reveals the inclusion
  \begin{equation}\label{_L3_}
    E_0^-\tp{\phi,\ep,k}\cup E_0^+\tp{\phi,\ep,k}\subseteq\tcb{y\in B_{\R^5}[0,3]\;:\;\tabs{\P y}^2\le c\tp{k + \ep}}.
  \end{equation}
  Synthesizing~\eqref{_L1_}, \eqref{_L2_}, and~\eqref{_L3_} yields
  \begin{equation}\label{_0th_}
    \tnorm{\mathscr{N}_{\ep,k}\tp{\phi}}_{L^1\tp{E_0}}\le C\tp{\ep + k}^2\text{ for }0<\ep\le1/c,\;\phi\in\mathcal{B}[\ep_2],\;0\le k\le1.
  \end{equation}

  We now consider the $n^{\text{th}}$ term in the series~\eqref{series_decomposition} for $n\ge 1$. From the third item of Lemma~\ref{lem on a decomposition and properties}, we are justified in equating
  \begin{equation}\label{nth term after coarea}
    \tnorm{\mathscr{N}_{\ep,k}\tp{\phi}}_{L^1\tp{E_n}} = \int_{-4/5}^{6/5}\int_{\tcb{\v + \phi - 1 = t}}\f{\mathds{1}_{I_n\tp{\ep,k}}\tp{t,y}}{\tabs{\grad\tp{\v + \phi}\tp{y}}}\tabs{\Gamma_\ep\tp{\tabs{\P y}^2t - k} - \Gamma_\ep\tp{t}}\;\m{d}\mathcal{H}^4\tp{y}\;\m{d}\mathcal{H}^1\tp{t}.
  \end{equation}
  after using the coarea formula. Notice that in the above we have introduced the sets
  \begin{equation}
    I_n\tp{\ep,k} = \tcb{\tp{t,y}\in[-4/5,6/5]\times E_n\;:\;\tabs{\Gamma_\ep\tp{\tabs{\P y}^2t - k} - \Gamma_\ep\tp{t}}>0}.
  \end{equation}
  If $\tp{t,y}\in I_n\tp{\ep,k}$, then it must be the case that
  \begin{equation}\label{annoying_casework}
   \tsb{t>0\text{ or }\tabs{\P y}^2t - k>0}\text{ and }\tsb{t<\ep\text{ or }\tabs{\P y}^2t - k<\ep}.
  \end{equation}
  In all possible cases of~\eqref{annoying_casework}, we deduce $0<t\le\min\tcb{6/5,2^{2n+4}\tp{\ep+k}}$ and hence
  \begin{equation}\label{nice_inclusion_}
    I_n\tp{\ep,k}\subseteq[0,\min\tcb{6/5,2^{2n + 4}\tp{\ep + k}}]\times E_n.
  \end{equation}

  Upon returning to expression~\eqref{nth term after coarea}, the upper bounds
  \begin{equation}\label{_LABEL_}
    \tnorm{\mathscr{N}_{\ep,k}\tp{\phi}}_{L^1\tp{E_n}}\le C\int_{0}^{\min\tcb{6/5,2^{2n + 4}\tp{\ep + k}}}\mathcal{H}^4\tp{E_n\cap\tcb{\v + \phi - 1 = t}}\;\m{d}\mathcal{H}^1\tp{t}
    \le C2^{-2n}\min\tcb{2^{-2n},\ep + k}
  \end{equation}
  are deduced from~\eqref{nice_inclusion_} and~\eqref{some_good_estimates}.

  The desired estimate~\eqref{lebesgue_bounds_WTS} in the case $p = 1$ now follows from~\eqref{series_decomposition}, \eqref{_0th_}, \eqref{_LABEL_}, and summing over $n\in\N$.
\end{proof}

The Lebesgue space estimates of Proposition~\ref{prop on Lebesgue space estimates on the perturbative operator} directly lead to perturbative $C^1$ estimates on the operator $N_{\ep,k}$.

\begin{coro}[More estimates on the perturbative operator]\label{coro on more estimates on the perturbative operator}
  Let $\ep_2>0$ be the small parameter granted by Lemma~\ref{lem on a decomposition and properties}. There exists a constant $C\in\R^+$ such that for all $0<\ep\le 1$, $0\le k\le 1$, and $\phi\in\mathcal{B}[\ep_2]$ the estimate
  \begin{equation}\label{the small estimate on the per_non}
    \tnorm{N_{\ep,k}\tp{\phi}}_{\mathcal{B}}\le C\tp{\ep + k}^{1/6}.
  \end{equation}
  holds.
\end{coro}
\begin{proof}
  By H\"older's inequality, we have the initial estimate 
  \begin{equation}
    \tnorm{N_{\ep,k}\tp{\phi}}_{\mathcal{B}}\le C\sup_{x\in B[0,3]}\tp{\tnorm{\tabs{x - \cdot}^{-3}}_{L^{6/5}(B(0,3))} + \tnorm{\tabs{x - \cdot}^{-4}}_{L^{6/5}(B(0,3))}}\tnorm{\mathscr{N}_{\ep,k}\tp{\phi}}_{L^6\tp{B[0,3]}}.
  \end{equation}
  Estimate~\eqref{the small estimate on the per_non} now follows from $4\cdot\tp{6/5}<5$ and Proposition~\ref{prop on Lebesgue space estimates on the perturbative operator}.
\end{proof}

\subsection{Swirl operator}\label{SS_swirl_operator}

In this section's final subsection, we examine several operators related to the swirl and helicity of steady vortex rings. The central components of these objects are the comb functions that we introduce next. Let us fix an even function $\chi\in C_\m{c}^\infty\tp{\R}$ satisfying 
\begin{equation}\label{_a_single_bump_}
  \chi\ge0,\;\supp\chi\subset[-1,1],\;\int_{-1}^1\tabs{\chi'\tp{t}}\;\m{d}t = 2,\text{ and }\chi(0) = \max_{\R}\chi = 1.
\end{equation}

Then, for $\N\ni n\ge 20$, we define the elementary spike function
\begin{equation}\label{_elementary_spike_}
  \sig_n\in C^\infty_{\m{c}}\tp{\R}\text{ via }\sig_n\tp{t} = \f{\chi\tp{n\tp{\log n}\tp{\log\log n}^2t}}{n\tp{\log n}\tp{\log\log n}}\text{ for }t\in\R.
\end{equation}
These comprise the comb functions $\Sigma_\rho\in C^\infty_{\m{c}}\tp{\R}$, defined for $\rho,t\in\R$ via
\begin{equation}\label{_the_comb_functions_}
  \Sigma_\rho\tp{t} = \sum_{n=20}^\infty\Gamma_2\tp{\rho - n}\sig_n\tp{\tp{100\cdot\del_\infty}t - \del_n}\text{ where }\del_n = \sum_{m=20}^n\f{4}{m\tp{\log m}\tp{\log\log m}^2}.
\end{equation}
We recall that $\Gamma_2$ is defined in~\eqref{defn_regularized_heaviside} and we define $\del_\infty = \lim_{n\to\infty}\del_n\in\R^+$. The role of the parameter $\rho$ is to smoothly select the number of nonzero terms in the series definition. We remark that there exists a constant $C\in\R^+$ such that
\begin{multline}\label{_remarks_on_the_comb_functions_}
  \bigcup_{\rho\in\R}\supp\Sigma_\rho\subseteq\bsb{0,\f{1}{100}};\;\rho\le20\imp\Sigma_\rho = 0;\;0\le\Sigma_\rho\le\f{1}{50};\;\tabs{\Sigma_\rho\cdot\Sigma_\rho'}\le C;\\\rho\ge22\imp\int_0^{1/100}\Sigma_\rho\tp{t}\;\m{d}t\ge\f{1}{C};\text{ and }\int_{0}^{1/100}\tabs{\Sigma_\rho'\tp{t}}\;\m{d}t\ge\sum_{n=20}^{\tfloor{\rho}-2}\f{2}{n\tp{\log n}\tp{\log\log n}}.
\end{multline}
Furthermore, thanks to support considerations, we have the equality
\begin{equation}\label{_another_comb_function_remark_}
  \tabs{\Sigma_\rho'\tp{t}} = \sum_{n=20}^\infty\Gamma_2\tp{\rho - n}\tp{100\cdot\del_\infty}\tabs{\sig_n'\tp{\tp{100\cdot\del_\infty}t - \del_n}}.
\end{equation}

We shall also demonstrate that for all $0\le\al< 1$
\begin{equation}\label{_more_uniform_estimates_on_the_comb_functions_}
  \sup_{\rho\in\R}\tnorm{\Sigma_\rho}_{\m{LL}\tp{\R}}\le C\text{ and }\sup_{\rho\in\R}\tnorm{\Sigma_\rho}_{C^\al_{\m{b}}\tp{\R}}\le C\tp{1 + \tabs{\log\tabs{1 -\al}}}.
\end{equation}
To see this, we let $\tilde{\sig}_n\in C^\infty_{\m{c}}\tp{\R}$ be defined via $\tilde{\sig}_n\tp{t}=\sig_n\tp{\tp{100\cdot\del_\infty}t-\del_n}$ for all $t\in\R$ and $\N\ni n\ge20$; we also set $\tilde{\sig}_n=0$ for $n<20$. Then
\begin{equation}\label{_the_sum_here_}
  \Sigma_\rho = \sum_{n=20}^\infty\Gamma_2\tp{\rho - n}\tilde{\sig}_n.
\end{equation}
One can directly estimate
\begin{multline}\label{estimates_on_the_spikes}
  \tnorm{\tilde{\sig}_n}_{C^0_{\m{b}}\tp{\R}}\lesssim\f{1}{n\tp{\log n}\tp{\log\log n}},\;\tnorm{\tilde{\sig}_n}_{C^\al_{\m{b}}\tp{\R}}\lesssim n^{\al-1}\tp{\log n}^{\al-1}\tp{\log\log n}^{2\al-1}\lesssim1+\tabs{\log\tp{1-\al}},\\
  \text{and }\tnorm{\tilde{\sig}_n}_{\m{LL}\tp{\R}}\lesssim\f{\log\log n}{\log\tp{n\tp{\log n}\tp{\log\log n}^2}}\lesssim1
\end{multline}
with implicit constants independent of $\N\ni n\ge20$ and $0\le\al<1$.

We now endeavor to estimate the quantity $\tabs{\Sigma_\rho(t_1) - \Sigma_\rho(t_0)}$ for some $t_0,t_1\in\R$ with $|t_0 - t_1|\le 1$. If neither point belongs to the support of a spike $\tilde{\sig}_n$, then the desired estimate is immediate. Otherwise, we may suppose, after swapping labels on $t_0$ and $t_1$, if necessary, that there exists $\N\ni n_0\ge 20$ such that 
\begin{equation}\label{this_thing_gets_a_quote}
  \babs{t_0 - \f{\del_{n_0}}{100\cdot\del_\infty}}\le\f{1}{100\cdot\del_\infty}\f{1}{n_0\tp{\log n_0}\tp{\log\log n_0}^2}
\end{equation} 
so that $t_0\in\supp\tilde{\sig}_{n_0}$.

There are three cases to consider. First, if 
\begin{equation}
  \tsb{n_0=20\text{ and }t_1\le\max\supp\tilde{\sig}_{21}}\text{ or }\tsb{n_0>20\text{ and }\min\supp\tilde{\sig}_{n_0-1}\le t_1\le\max\supp\tilde{\sig}_{n_0+1}}
\end{equation}
then we simply use~\eqref{estimates_on_the_spikes}
\begin{multline}\label{bounded_above_by_this_guy}
  \tabs{\Sigma_\rho\tp{t_1} - \Sigma_\rho\tp{t_0}}\le\sum_{\epsilon = -1}^1\tabs{\Gamma_2\tp{\rho - \tp{n_0 + \epsilon}}}\tabs{\tilde{\sig}_{n_0+\epsilon}\tp{t_1} - \tilde{\sig}_{n_0+\epsilon}\tp{t_0}}\\\lesssim
    \min\tcb{\tp{1 + \tabs{\log\tabs{1 - \al}}}\tabs{t_1 - t_0}^\al,\tabs{t_1 - t_0}\tp{1 + \tabs{\log\tabs{t_1 - t_0}}}}.
\end{multline}

For the second case, we suppose that
\begin{equation}\label{far-away-condition}
  t_1>\max\supp\tilde{\sig}_{n_0 + 1}\text{ or }\tsb{n_0>20\text{ and }t_1<\min\supp\tilde{\sig}_{n_0-1}}
\end{equation}
and that there does not exist $\N\ni n\ge 20$ such that $t_1\in\supp\tilde{\sig}_n$. Then, from~\eqref{this_thing_gets_a_quote}, we know
\begin{equation}
  \tabs{\Sigma_\rho\tp{t_1} - \Sigma_\rho\tp{t_0}} = \tabs{\Gamma_2\tp{\rho - n_0}}\tabs{\tilde{\sig}_{n_0}\tp{t_0}}\lesssim\f{1}{n_0\tp{\log n_0}\tp{\log\log n_0}} \text{ and }|t_1 - t_0|\gtrsim\f{1}{n_0\tp{\log n_0}\tp{\log\log n_0}^2}.
\end{equation}
It follows that
\begin{multline}\label{this_one_is_also_too_cited}
  \tabs{\Sigma_\rho(t_1) - \Sigma_\rho(t_0)}\lesssim\min\Big\{n_0^{\al-1}\tp{\log n_0}^{\al-1}\tp{\log\log n_0}^{2\al - 1}|t_1 - t_0|^\al,\\\f{\log\log n_0}{\log\tp{n_0\tp{\log n_0}\tp{\log\log n_0}^2}}\tabs{t_1 - t_0}\tp{1 + \tabs{\log\tabs{t_1 - t_0}}}\Big\},
\end{multline}
which can be bounded above by the right hand side of~\eqref{bounded_above_by_this_guy}.

The remaining case is that $t_1$ satisfies~\eqref{far-away-condition} and there exists $\N\ni n_1\ge20$ such that $t_1\in\supp\tilde{\sig}_{n_1}$. Then, we are guaranteed that
\begin{equation}
  \tabs{t_1 - t_0}\gtrsim\max\bcb{\f{1}{n_0\tp{\log n_0}\tp{\log\log n_0}^2},\f{1}{n_1\tp{\log n_1}\tp{\log\log n_1}^2}}.
\end{equation}
Whence
\begin{equation}
  \tabs{\Sigma_\rho\tp{t_1} - \Sigma_\rho(t_0)}\le\tabs{\tilde{\sig}_{n_1}\tp{t_1}} + \tabs{\tilde{\sig}_{n_0}\tp{t_0}}\lesssim\max\bcb{\f{1}{n_0\tp{\log n_0}\tp{\log\log n_0}},{\f{1}{n_1\tp{\log n_1}\tp{\log\log n_1}}}}.
\end{equation}
Thus, an estimate of the form~\eqref{this_one_is_also_too_cited} follows as above and is therefore bounded by the right hand side of~\eqref{bounded_above_by_this_guy}. This completes the justification of the bounds~\eqref{_more_uniform_estimates_on_the_comb_functions_}.

Before we can give a definition of our swirl operators incorporating the comb functions~\eqref{_the_comb_functions_}, we require some simple facts about the interior stream lines of Hill's vortex.

\begin{lem}[Bounds on interior streamlines]\label{lem on bounds on interior stream lines}
  For $\phi\in\mathcal{B}[1/2]$ we define the function $\v_\phi\in\mathcal{B}$ via
  \begin{equation}\label{almost_hills_solution}
    \v_\phi\tp{y} = \tabs{\P y}^2\tp{\v\tp{y} + \phi\tp{y} - 1}\text{ for }y\in B_{\R^5}[0,3].
  \end{equation}
  There exists $0<\ep_3\le1/2$ and $C\in\R^+$ with the property that for all $\phi\in\mathcal{B}[\ep_3]$ the following hold.
  \begin{enumerate}
    \item We have the compact inclusion
    \begin{equation}
      F^\phi = \tcb{y\in B_{\R^5}[0,3]\;:\;\v_\phi\tp{y}\in[1/16,5/16]}\Subset B_{\R^5}(0,1).
    \end{equation}
    \item For all $y\in F^\phi$ the estimates 
    \begin{equation}\label{some essential coffee flavorz}
      1/C\le\tabs{\grad\v_\phi\tp{y}}\le C\text{ and }1/C\le\tabs{\P y}\le C
    \end{equation}
    hold.
    \item For all $1/16\le t\le 5/16$ the estimates 
    \begin{equation}\label{do_you_notice_my_labels}
      1/C\le\mathcal{H}^4\tcb{y\in F^\phi\;:\;\v_\phi\tp{y} = t}\le C
    \end{equation}
    hold.
  \end{enumerate}
\end{lem}
\begin{proof}
  All but one of the assertions made in the statement of the lemma are direct consequences of Theorem~\ref{thm on uniform density estimates on families of hypersurfaces}; we shall check that the hypotheses of this result are satisfied with
  \begin{equation}
    U = B_{\R^5}\tp{0,3},\;\Psi\in C^1_{\m{b}}\tp{\Bar{U}}\text{ defined via }\Psi(y) = \v_0\tp{y} = \tabs{\P y}^2\tp{\v\tp{y} - 1}\text{ and }a = \f{1}{16},\;b = \f{5}{16}.
  \end{equation}
  Recall that Hill's solution $\v$ is explicitly given in equation~\eqref{Hills_Vortex}. By direct computation, one observes
  \begin{multline}
    \tcb{\Psi<0} = \tcb{y\in B_{\R^5}(0,3)\;:\;|y|>1,\;\P y\neq 0},\;\tcb{\Psi>0} = \tcb{y\in B_{\R^5}(0,3)\;:\;|y|<1,\;\P y\neq0},\\
    \text{and }\tcb{\Psi = 0} = \tcb{y\in B_{\R^5}(0,3)\;:\;\P y = 0\text{ or }|y|=1}.
  \end{multline}
  For $|y|<1$ we calculate 
  \begin{equation}
    \grad\Psi\tp{y} = 3\tp{\P y\tp{1 - \tabs{y}^2} - \tabs{\P y}^2y}
  \end{equation}
  and hence 
  \begin{equation}
    \grad\Psi\tp{y}=0\text{ if and only if }\tsb{\P y=0}\text{ or }\tsb{y_5 = 0\text{ and }|y|^2 = 1/2}.
  \end{equation}

  From this we deduce that $\Psi$ attains its maximum value of $3/8$ over the latter set above. Armed with these observations, one readily deduces that
  \begin{multline}\label{thanks_to_you}
    F^0=\tcb{y\in B_{\R^5}\tsb{0,3}\;:\;\Psi\tp{y}\in[1/16,5/16]}\Subset B_{\R^5}\tp{0,1},\;\min_{F^0}\tabs{\grad\Psi}>0,\;\min\tcb{|\P y|\;:\;y\in F^0}>0\\
    \text{and }\forall\;t\in[1/16,5/16],\;\tcb{y\in B_{\R^5}[0,3]\;:\;\Psi\tp{y} = t}\neq\es.
  \end{multline}

  For $\phi\in\mathcal{B}[1/2]$, the perturbation of $\Psi$ appearing in $\v_\phi$ is the function $y\mapsto\tabs{\P y}^2\phi\tp{y}$, whose $C^{0,1}$ norm is bounded by $C\tnorm{\phi}_{\mathcal{B}}$. Consequently, after decreasing the perturbation size if necessary, the hypotheses of Theorem~\ref{thm on uniform density estimates on families of hypersurfaces} are satisfied thanks to~\eqref{thanks_to_you}. Upon invoking this result we obtain $0<\ep_3\le1/2$ and $C\in\R^+$ such that the three items of the lemma are satisfied except possibly for the right hand inequality of~\eqref{some essential coffee flavorz}. But, by making $\ep_3$ smaller and $C$ larger, if necessary, this too can be satisfied.
\end{proof}

The previous lemma and simple symmetry considerations ensure that all of the integrals and operators in the forthcoming definition are well-defined.

\begin{defn}[Swirl and helicity operators]\label{defn of swirl and helicity operators}
  Let $0<\ep_3\le1/2$ be the parameter given by Lemma~\ref{lem on bounds on interior stream lines}. Given $\rho\in\R$ and $i\in\tcb{1,2}$ we define the following functions.
  \begin{enumerate}
    \item The swirl operator is the map $M^i_\rho:\mathcal{B}[\ep_3]\to\mathcal{B}$ with the action
    \begin{equation}\label{_the_swirl_operator_}
      M^i_\rho\tp{\phi} = \f{1}{8\pi^2}\int_{B_{\R^5}\tp{0,3}}\f{1}{\tabs{\cdot - y}^3}\f{1}{\tabs{\P y}^2}\tp{\Sigma_\rho\cdot\Sigma_\rho'}\tp{\v_\phi\tp{y} - i/8}\;\m{d}\mathcal{H}^5\tp{y}.
    \end{equation}
    \item The principal partial helicity functional is the map $I^i_\rho:\mathcal{B}[\ep_3]\to\R$ with the action
    \begin{equation}\label{_the_p_par_heli_}
      I^i_\rho\tp{\phi} = \f{30}{\pi}\int_{B_{\R^5}\tp{0,3}}\f{1}{\tabs{\P y}^2}\Sigma_\rho\tp{\v_\phi\tp{y} - i/8}\;\m{d}\mathcal{H}^5\tp{y}.
    \end{equation}
    \item The perturbative partial helicity functional is the map $J^i_\rho:\mathcal{B}[\ep_3]\to\R$ with the action 
    \begin{equation}\label{_the_pe_par_heli_}
      J^i_\rho\tp{\phi} = \f{2}{\pi}\int_{B_{\R^5}\tp{0,3}}\f{1}{\tabs{\P y}^4}\tp{\Sigma_\rho\cdot\Sigma_\rho\cdot\Sigma_\rho'}\tp{\v_{\phi}\tp{y} - i/8}\;\m{d}\mathcal{H}^5\tp{y}.
    \end{equation}
    \item The total meridional vorticity functional is the map $T^i_\rho:\mathcal{B}[\ep_3]\to\R$ with the action
    \begin{equation}\label{total_meridional_vorticty}
      T^i_\rho\tp{\phi} = \f{1}{\pi}\int_{B_{\R^5}\tp{0,3}}\f{1}{\tabs{\P y}^3}\tabs{\Sigma_\rho'\tp{\v_\phi\tp{y} - i/8}}\tabs{\grad\v_\phi\tp{y}}\;\m{d}\mathcal{H}^5\tp{y}.
    \end{equation} 
  \end{enumerate}
\end{defn}

Our next result records simple properties of the maps $M$, $I$, and $J$ which were introduced in Definition~\ref{defn of swirl and helicity operators}. 

\begin{lem}[Simple estimates on the swirl and helicity operators]\label{lem on properties of the M I and J}
  Let $0<\ep_3\le1/2$ be the positive parameter given to us by Lemma~\ref{lem on bounds on interior stream lines}. There exist constants $C,\mathfrak{C}\in\R^+$ such that for all $i\in\tcb{1,2}$, $\rho\in\R$, and $\phi\in\mathcal{B}[\ep_3]$ the following hold.
  \begin{enumerate}
    \item We have the inclusion $M^i_\rho\tp{\phi}\in\m{LL}^1\tp{B_{\R^5}[0,3]}$ and the bound $\tnorm{M^i_\rho\tp{\phi}}_{\m{LL}^1}\le C$.
    \item We have the bounds $\tabs{I^i_\rho\tp{\phi}}\le C$ and $\tabs{J^i_\rho\tp{\phi}}\le C$.
    \item If additionally $\rho\ge22$, then $I^i_\rho\tp{\phi}\ge1/C$ and
    \begin{equation}\label{the weird constant thing that we need to cite later}
      I^1_\rho\tp{\phi}\le\mathfrak{C}I^2_\rho\tp{\phi}.
    \end{equation}
  \end{enumerate}
  Moreover, the induced map
  \begin{equation}\label{_the_induced_maps_}
    \R\times\mathcal{B}[\ep_3]\ni\tp{\rho,\phi}\mapsto\tp{M^i_\rho\tp{\phi},I^i_\rho\tp{\phi},J^i_\rho\tp{\phi}}\in\mathcal{B}\times\R\times\R
  \end{equation}
  is continuous.
\end{lem}
\begin{proof}
  We note that by support considerations, see~\eqref{_remarks_on_the_comb_functions_}, each of the defining integrands of equations~\eqref{_the_swirl_operator_}, \eqref{_the_p_par_heli_}, and~\eqref{_the_pe_par_heli_} are supported on the set
  \begin{equation}
    Q_i(\phi) = \tcb{y\in B_{\R^5}[0,3]\;:\;\v_{\phi}\tp{y}-i/8\in[0,1/100]}.
  \end{equation}
  From Lemma~\ref{lem on bounds on interior stream lines}, we get $Q_i\tp{\phi}\subset F^\phi\Subset B_{\R^5}(0,1)$ and the two sided bounds of~\eqref{some essential coffee flavorz} and \eqref{do_you_notice_my_labels} hold.

  Upon once more using~\eqref{_remarks_on_the_comb_functions_}, we deduce immediately that the integrands defining the $I$ and $J$ functionals are bounded independently of $\rho\in\R$ and $\phi\in\mathcal{B}[\ep_3]$. So the second item trivially follows.

  The first item is similarly simple. The operator $M$ is the convolution of the Newtonian potential with a bounded function of bounded support. The same strategy that appeared in the proof of the first item of Proposition~\ref{prop on basic smoothness} is directly applicable in this situation as well.

  We turn our attention to the third item. Assuming that $\rho\ge 22$, we write $I^i_\rho\tp{\phi}$ by using the coarea formula:
  \begin{equation}\label{co_area_form_of_I}
    I^i_\rho\tp{\phi} = \f{30}{\pi}\int_{i/8}^{i/8+1/100}\Sigma_\rho\tp{t-i/8}\int_{\tcb{\v_\phi = t}}\f{1}{\tabs{\P y}^2}\f{1}{\tabs{\grad\v_\phi\tp{y}}}\;\m{d}\mathcal{H}^4\tp{y}\;\m{d}\mathcal{H}^1\tp{t}.
  \end{equation}
  Then, from~\eqref{_remarks_on_the_comb_functions_} and Lemma~\ref{lem on bounds on interior stream lines} we deduce from~\eqref{co_area_form_of_I} the lower bound
  \begin{equation}
    I^i_\rho\tp{\phi} \ge\f{1}{C}\int_{i/8}^{i/8 + 1/100}\Sigma_\rho\tp{t-i/8}\;\m{d}\mathcal{H}^1\tp{t}\ge\f{1}{C}.
  \end{equation}
  The assertions of the third item are now clear. The continuity assertion follows from the local finiteness of the series defining $\Sigma_\rho$, the uniform support and bounds above, the dominated convergence theorem, and the same Newtonian-potential estimates used for $M^i_\rho$.
\end{proof}

For this subsection's final result, we shall select the parameter $\rho$ appearing in Definition~\ref{defn of swirl and helicity operators} as a continuous function of $\phi$ and several other parameters in such a way that we lock the total meridional vorticity of the corresponding velocity field. In the forthcoming proposition statement, the small positive value $\ep_3$ is that of Lemma~\ref{lem on bounds on interior stream lines}, the constant $\mathfrak{C}$ is that of equation~\eqref{the weird constant thing that we need to cite later}, and the maps $T$ are those of~\eqref{total_meridional_vorticty}.

\begin{prop}[Fixing the total meridional vorticity]\label{prop on fixing the total meridional vorticity}
  There exists $\R\ni\mu_0>0$ and a continuous map
  \begin{equation}\label{a_tortured_continuous_function}
    \pmb{\varrho}:\mathcal{B}[\ep_3]\times[0,2\mathfrak{C}]\times[\mu_0,\infty)\to[22,\infty)
  \end{equation}
  with the property that for all $\phi\in\mathcal{B}[\ep_3]$, $0\le\lambda\le2\mathfrak{C}$, and $\mu\ge\mu_0$ the identity 
  \begin{equation}\label{the_bold_varrho}
    T^1_\rho\tp{\phi} + \lambda T^2_\rho\tp{\phi} = \mu\text{ with }\rho = \pmb{\varrho}\tp{\phi,\lambda,\mu}
  \end{equation}
  is satisfied.
\end{prop}
\begin{proof}
  Let $\bf{P} = \mathcal{B}[\ep_3]\times[0,2\mathfrak{C}]$ and define the function
  \begin{equation}
    \mathsf{f}:\bf{P}\times\R\to\R\text{ via }\mathsf{f}\tp{\tp{\phi,\lambda},\rho} = T^1_\rho\tp{\phi} + \lambda T^2_\rho\tp{\phi}\text{ for all }\tp{\phi,\lambda}\in\bf{P},\;\rho\in\R.
  \end{equation}
  It is clear from the definition of $T$ in~\eqref{total_meridional_vorticty} that $\mathsf{f}$ is a continuous function. By inspection of the identity (that follows from~\eqref{_another_comb_function_remark_})
  \begin{multline}\label{the real expression for the function f}
    \mathsf{f}\tp{\tp{\phi,\lambda},\rho} = \f{1}{\pi}\sum_{n=20}^\infty\Gamma_2\tp{\rho - n}\int_{B_{\R^5}\tp{0,3}}\f{100\cdot\del_\infty}{\tabs{\P y}^3}\tabs{\tp{\sig_n}'\tp{\tp{100\cdot\del_\infty}\tp{\v_\phi\tp{y} - 1/8} - \del_n}}\tabs{\grad\v_\phi\tp{y}}\;\m{d}\mathcal{H}^5\tp{y}\\
    + \f{\lambda}{\pi}\sum_{n=20}^\infty\Gamma_2\tp{\rho - n}\int_{B_{\R^5}\tp{0,3}}\f{100\cdot\del_\infty}{\tabs{\P y}^3}\tabs{\tp{\sig_n}'\tp{\tp{100\cdot\del_\infty}\tp{\v_\phi\tp{y} - 1/4} - \del_n}}\tabs{\grad\v_\phi\tp{y}}\;\m{d}\mathcal{H}^5\tp{y}.
  \end{multline}
  It is also clear that for all $\tp{\phi,\lambda}\in\bf{P}$ the map $\R\ni\rho\mapsto\mathsf{f}\tp{\tp{\phi,\lambda},\rho}\in\R$ is smooth.

  By combining identity~\eqref{the real expression for the function f} with Lemma~\ref{lem on bounds on interior stream lines}, the coarea formula, and the definition of $\sig_n$ from~\eqref{_elementary_spike_} (see also~\eqref{_remarks_on_the_comb_functions_}), we observe the existence of a constant $C\in\R^+$ such that for all $\R\ni\rho\ge22$ and $\tp{\phi,\lambda}\in\bf{P}$ the bounds 
  \begin{equation}\label{two_sided_increasing_bound}
    \f{1+\lambda}{C}\sum_{n=20}^{\tfloor{\rho} - 2}\f{1}{n\tp{\log n}\tp{\log\log n}}\le\mathsf{f}\tp{\tp{\phi,\lambda},\rho}\le C\tp{1 + \lambda}\sum_{n=20}^{\tfloor{\rho}}\f{1}{n\tp{\log n}\tp{\log\log n}}
  \end{equation}
  and (letting $D_2$ denote the derivative in the $\rho$-variable)
  \begin{equation}\label{two_sided_derivative_bound}
    \f{1+\lambda}{C}\f{1}{\rho\tp{\log\rho}\tp{\log\log\rho}} \le D_2\mathsf{f}\tp{\tp{\phi,\lambda},\rho}\le C\tp{1 + \lambda}\f{1}{\rho\tp{\log\rho}\tp{\log\log\rho}}
  \end{equation}
  hold. For the lower bound in~\eqref{two_sided_derivative_bound}, we use the positivity of $\gam'$ on $(0,1)$ and the overlap of two consecutive transition intervals. More precisely, $\min_{0\le s\le1}\tp{\Gamma_2'\tp{s}+\Gamma_2'\tp{s+1}}>0$, and the two corresponding summands are comparable to $\rho^{-1}\tp{\log\rho}^{-1}\tp{\log\log\rho}^{-1}$.

  The estimate~\eqref{two_sided_increasing_bound}, when combined with the intermediate value theorem, gives us the existence of $\R\ni\mu_0>1$ such that 
  \begin{equation}\label{totally_surjective}
    \bigcap_{\tp{\phi,\lambda}\in\bf{P}}\tcb{\mathsf{f}\tp{\tp{\phi,\lambda},\rho}\;:\;\R\ni\rho\ge22}\supset[\mu_0-1,\infty).
  \end{equation}
  Then, the estimate~\eqref{two_sided_derivative_bound} implies that the map $\rho\mapsto \mathsf{f}(\tp{\phi,\lambda},\rho)$ is strictly increasing. Combining this observation with~\eqref{totally_surjective} verifies the next assertion: For all $\tp{\phi,\lambda}\in\bf{P}$ and all $\R\ni\mu\ge\mu_0$ there exists a unique $\R\ni\pmb{\varrho}\tp{\phi,\lambda,\mu}\ge22$ such that $\mathsf{f}\tp{\tp{\phi,\lambda},\pmb{\varrho}\tp{\phi,\lambda,\mu}} = \mu$.

  The proof is complete as soon as we verify the map $\tp{\phi,\lambda,\mu}\mapsto\pmb{\varrho}\tp{\phi,\lambda,\mu}$ to be continuous. But this is a direct consequence of the implicit function theorem with continuous parameter, see, e.g., Proposition 2.1 in Gl\"{o}ckner~\cite{MR2269430}. The hypotheses of this result are satisfied, in particular, due to the positive lower bound in~\eqref{two_sided_derivative_bound}.
\end{proof}

\section{Inversion near Hill's vortex}\label{SO_uni_inv_near_Hill}

Our task in this section is to invert the principal operator $F_\ep$ on a neighborhood of Hill's vortex, uniformly as $\ep\to0$. We first prove nondegeneracy of the limiting linearization by spherical harmonics and then use compactness to obtain the uniform estimates required for a quantitative local surjectivity theorem.


\subsection{Nondegeneracy}\label{SS_nondegen}

Our goal now is to establish that Hill's spherical vortex in the $\R^5$ formulation~\eqref{Hills_Vortex} is nondegenerate in a suitable sense. Roughly speaking, we shall show that the defining functional equation in identity~\eqref{Hills_Integral_Equation} when linearized at Hill's vortex is an operator with a trivial kernel on the symmetry class considered below.

The verification of the aforementioned fact is most comfortably performed with the help of the $\S^4\subset\R^5$ spherical harmonics. For a general reference on the subject, the reader is referred to Atkinson and Han~\cite{MR2934227}; we shall rapidly review here some key elements while fixing notation.

The base Hilbert space is
\begin{equation}\label{space_of_monke}
  L^2\tp{\S^4}\text{ with the inner product }\tbr{f,g}_{L^2\tp{\S^4}} = \int_{\S^4}fg\;\m{d}\mathcal{H}^4.
\end{equation}
For every $n\in\N$ there is a subspace
\begin{equation}\label{spherical_harmonic_space}
  \Y_n\subset L^2\tp{\S^4}\text{ with }\m{dim}\Y_n = \tp{2n + 3}\tp{n+2}\tp{n+1}/6
\end{equation}
called the spherical harmonic space of order $n$. Elements of $\Y_n$ are the restrictions of homogeneous harmonic polynomials in $\R^5$ of degree $n$ to $\S^4$. The spherical harmonic subspaces are mutually orthogonal and collectively complete:
\begin{equation}\label{completeness_of_spherical_harmonics}
  L^2\tp{\S^4} = \bigoplus_{n=0}^\infty\Y_n.
\end{equation}

For our purposes, the most important feature of the spherical harmonic subspaces is that they are invariant under linear integral operators on $\S^4$ with angle-dependent kernels. More precisely, the Funk-Hecke formula holds, see, e.g. Theorem 2.22 in~\cite{MR2934227}. If 
\begin{equation}
  f:\tp{-1,1}\to\R\text{ is measurable and }\int_{-1}^1\tabs{f\tp{t}}\tp{1 - t^2}\;\m{d}\mathcal{H}^1\tp{t}<\infty
\end{equation}
then for each $n\in\N$ and $Y\in\Y_n$ we have the Funk-Hecke formula
\begin{equation}\label{Funk_Hecke}
  \int_{\S^4}f\tp{\xi\cdot\eta}Y\tp{\eta}\;\m{d}\mathcal{H}^4\tp{\eta} = \lambda_n Y\tp{\xi}\text{ for all }\xi\in\S^4
\end{equation}
where the constant $\lambda_n\in\R$ is determined via
\begin{equation}\label{constant_in_Funk_Hecke}
  \lambda_n = 2\pi^2\int_{-1}^1P_n\tp{t}f(t)\tp{1 - t^2}\;\m{d}\mathcal{H}^1\tp{t}.
\end{equation} 

In identity~\eqref{constant_in_Funk_Hecke} the functions $P_n$ are the dimension-five Legendre polynomials. While there are many ways to write explicitly what these functions are, the most convenient for our application is Rodrigues' representation formula:
\begin{equation}\label{legendre_polynomials}
  P_n\tp{t} = \f{\tp{-1}^n}{2^n\tp{n + 1}!}\f{1}{1 - t^2}\bp{\f{\m{d}}{\m{d}t}}^n\ssb{\tp{1 - t^2}^{n+1}}\text{ for }t\in\tp{-1,1}.
\end{equation}
See, e.g., Theorem 2.23 in~\cite{MR2934227}.

We are finally ready for this subsection's first result.
\begin{lem}[A coefficient computation]\label{lem on coefficient computation}
  When specializing to the case
  \begin{equation}\label{kernel_function}
    f(t) = \f{1}{\tp{2-2t}^{3/2}}\text{ for }t\in\tp{-1,1}
  \end{equation}
  in the Funk-Hecke coefficient formula~\eqref{constant_in_Funk_Hecke}, we obtain 
  \begin{equation}\label{diagonalization}
    \lambda_n = \f{8\pi^2}{2n+3}\text{ for every }n\in\N.
  \end{equation}
\end{lem}
\begin{proof}
  We substitute the specific choice of $f$~\eqref{kernel_function} and Rodrigues' formula~\eqref{legendre_polynomials} into the definition of $\lambda_n$~\eqref{constant_in_Funk_Hecke}. Then, we integrate by parts $n$ times, noting that all of the boundary terms vanish. Thus we have 
  \begin{equation}\label{an_intermediate_identity}
    \lambda_n = \f{2\pi^2}{2^n\tp{n+1}!}\int_{-1}^1\tp{1 - t^2}^{n+1}\bp{\f{\m{d}}{\m{d}t}}^n\ssb{\tp{2 - 2t}^{-3/2}}\;\m{d}\mathcal{H}^1\tp{t}.
  \end{equation}

  Next, it is a direct computation to establish that
  \begin{equation}\label{a_simple_identity}
    \bp{\f{\m{d}}{\m{d}t}}^{n}\ssb{\tp{2 - 2t}^{-3/2}} = 2^n\f{\Upgamma\tp{n + 3/2}}{\Upgamma\tp{3/2}}\tp{2 - 2t}^{-3/2 - n}
  \end{equation}
  holds for all $n\in\N$, where $\Upgamma\tp{\cdot}$ denotes the gamma function. Then, we substitute~\eqref{a_simple_identity} into~\eqref{an_intermediate_identity} to get that 
  \begin{equation}\label{insane_identity}
    \lambda_n = \f{2\pi^2\Upgamma\tp{n+3/2}}{\Upgamma\tp{n+2}\Upgamma\tp{3/2}}\f{1}{2^{n+3/2}}\int_{-1}^1\tp{1-t}^{-1/2}\tp{1+t}^{n+1}\;\m{d}\mathcal{H}^1\tp{t}.
  \end{equation}

  The final integral above can be linked to the gamma function via the beta function. It holds that 
  \begin{equation}\label{another_insane_identity}
    \int_{-1}^1\tp{1 - t}^{-1/2}\tp{1 + t}^{n+1}\;\m{d}\mathcal{H}^1\tp{t} = 2^{n+3/2}\int_{0}^1\tp{1-\tau}^{-1/2}\tau^{n+1}\;\m{d}\mathcal{H}^1\tp{\tau} = 2^{n+3/2}\f{\Upgamma\tp{1/2}\Upgamma\tp{n+2}}{\Upgamma\tp{n+5/2}}.
  \end{equation}
  The synthesis of~\eqref{insane_identity} and~\eqref{another_insane_identity} gives the desired formula~\eqref{diagonalization}.
\end{proof}

We are now ready to establish that in our specific formulation of the problem the linearized axisymmetric steady Euler equations at Hill's vortex have a trivial kernel. Similar manifestations of the nondegeneracy of Hill's spherical vortex are employed in the work of Norbury~\cite{MR302044}, Amick and Turner~\cite{MR929976}, and Wan~\cite{MR946444}.

Recall the notation for the subgroup $\O$ of the orthogonal group introduced in the first item of Definition~\ref{defn on spaces and mappings}.

\begin{prop}[Nondegeneracy of Hill's spherical vortex]\label{prop on nondegeneracy of Hills spherical vortex}
  The following are equivalent for $\phi\in C^0\tp{\S^4}$:
  \begin{enumerate}
    \item The identity
  \begin{equation}\label{the_formal_linearization_kernel_condition}
    \phi\tp{x} = \f{5}{8\pi^2}\int_{\S^4}\f{1}{\tabs{x - y}^3}\phi\tp{y}\;\m{d}\mathcal{H}^4\tp{y}
  \end{equation}
  holds for all $x\in\S^4$.
  \item We have the inclusion $\phi\in\Y_1$, where $\Y_1$ is the unit order spherical harmonic subspace as in~\eqref{spherical_harmonic_space}.
  \end{enumerate}
  
  As a consequence, if $\phi\in C^0\tp{\S^4}$ satisfies~\eqref{the_formal_linearization_kernel_condition} and the invariance condition:
  \begin{equation}\label{delete_the_kernel_with_love}
    \phi(Tx) = \phi(x)\text{ for all }x\in\S^4\text{ and }T\in\O,
  \end{equation}
  then $\phi = 0$.
\end{prop}
\begin{proof}
  Let us introduce the bounded linear operator 
  \begin{equation}
    Q:C^0\tp{\S^4}\to C^0\tp{\S^4},\;\tp{Q\varphi}\tp{x} = \f{5}{8\pi^2}\int_{\S^4}\f{1}{\tabs{x - y}^3}\varphi\tp{y}\;\m{d}\mathcal{H}^{4}\tp{y}
  \end{equation}
  defined for $x\in\S^4$ and $\varphi\in C^0\tp{\S^4}$.

  The Funk-Hecke formula, see identity~\eqref{Funk_Hecke}, ensures that the spherical harmonic subspaces $\tcb{\Y_n}_{n\in\N}$ are invariant under the action of the operator $Q$. More precisely, we see that by applying Lemma~\ref{lem on coefficient computation} and the identity 
  \begin{equation}
    \tabs{x - y}^{-3} = \tp{2 - 2x\cdot y}^{-3/2}\text{ for all }x,y\in\S^4
  \end{equation} 
  the equality
  \begin{equation}\label{it_is_now_immediate}
    \tp{QY}\tp{x} = \f{5}{2n+3}Y\tp{x}\text{ for all }x\in\S^4
  \end{equation}
  holds for every $n\in\N$ and $Y\in\Y_n$.

  Armed with identity~\eqref{it_is_now_immediate}, it is now clear that the second item of the asserted equivalence implies the first item. For the opposite direction, let us fix $\phi\in C^0\tp{\S^4}$ a solution to~\eqref{the_formal_linearization_kernel_condition}, i.e. $\phi = Q\phi$. If $n\in\N$ and $Y\in\Y_n$, then
  \begin{equation}
    \tbr{\phi,Y}_{L^2\tp{\S^4}} = \tbr{\phi,QY}_{L^2\tp{\S^4}} = \f{5}{2n + 3}\tbr{\phi,Y}_{L^2\tp{\S^4}}.
  \end{equation}
  Note that we have employed the inner product notation of~\eqref{space_of_monke}, the symmetry of $Q$, and identity~\eqref{it_is_now_immediate}. When $n\neq 1$, we deduce from the above that $\tbr{\phi,Y}_{L^2\tp{\S^4}} = 0$. The completeness of the spherical harmonics then implies that $\phi\in\Y_1$.

  We now consider what happens when $\phi\in C^0\tp{\S^4}$ solves~\eqref{the_formal_linearization_kernel_condition} and obeys the symmetries discussed in equation~\eqref{delete_the_kernel_with_love}. The previously studied equivalence implies that $\phi\in\Y_1$. The validity of the conditions~\eqref{delete_the_kernel_with_love} then permits us to apply a certain structural result, namely Theorem 2.8 in Atkinson and Han~\cite{MR2934227}. The conclusion is that $\phi$ necessarily obeys
  \begin{equation}\label{structural form}
    \phi\tp{x} = \phi\tp{e_5}P_1\tp{x_5} = \phi\tp{e_5}x_5,
  \end{equation}
  where $P_1$ is one of the Legendre polynomials, see~\eqref{legendre_polynomials}.

  Now from~\eqref{structural form} we see that $\phi(x_1,\dots,x_4,x_5) = -\phi\tp{x_1,\dots,x_4,-x_5}$. But from the symmetry conditions~\eqref{delete_the_kernel_with_love} it is also the case that $\phi(x_1,\dots,x_4,x_5) = \phi\tp{x_1,\dots,x_4,-x_5}$. Hence $\phi=0$ and the proof is complete.
\end{proof}

\subsection{Consequences}\label{SS_consequences}

The nondegeneracy of Hill's spherical vortex was established in Proposition~\ref{prop on nondegeneracy of Hills spherical vortex}. We now endeavor to promote this qualitative condition to a quantitative statement regarding the operators $A_{\ep}\tp{\phi}$, which we recall are defined in the second item of Proposition~\ref{prop on basic smoothness}.

\begin{thm}[Quantitative closed range estimates]\label{thm on quantitative closed range estimates}
  Let $0<\ep_1\le1/2$ be the small parameter granted to us by Proposition~\ref{prop on principal operator lipschitz estimate}. There exists $0<\tilde{\ep}_1\le\ep_1$ and $C\in\R^+$ with the property that for all $0<\ep\le\tilde{\ep}_1$, $\phi^0,\phi^1\in\mathcal{B}[\tilde{\ep}_1]$, and $\varphi\in\mathcal{B}^0$
  we have the estimate
  \begin{equation}\label{_all_hail_satan_}
    \tnorm{\varphi}_{\mathcal{B}^0}\le C\bnorm{\varphi + \bsb{\int_0^1A_{\ep}\tp{\tau\phi^0 + \tp{1 - \tau}\phi^1}\;\m{d}\tau}\varphi}_{\mathcal{B}^0}.
  \end{equation}
\end{thm}
\begin{proof}
  We argue via contradiction. If the claimed result were in fact false, then there would necessarily exist sequences 
  \begin{equation}\label{__sequences___}
    \tcb{\ep_n}_{n\in\N}\subset(0,\ep_1],\;\tcb{\phi^0_n}_{n\in\N},\tcb{\phi^1_n}_{n\in\N}\subset\mathcal{B}[\ep_1],\text{ and }\tcb{\varphi_n}_{n\in\N}\subset\mathcal{B}^0
  \end{equation}
  satisfying
  \begin{equation}\label{___absurd_things___}
    \ep_n+\tnorm{\phi_n^0}_{\mathcal{B}}+\tnorm{\phi_n^1}_{\mathcal{B}} + \tnorm{\varphi_n + \Phi_n}_{\mathcal{B}^0}\to0\text{ as }n\to\infty\text{ and }\forall\;n\in\N,\;\tnorm{\varphi_n}_{\mathcal{B}^0} = 1,
  \end{equation}
  where for $n\in\N$ we have denoted
  \begin{equation}\label{the_definition_of_capital_Phi}
    \Phi_n = \bsb{\int_0^1A_{\ep_n}\tp{\tau\phi^0_n + \tp{1 - \tau}\phi^1_n}\;\m{d}\tau}\varphi_n\in C^0\tp{B_{\R^5}[0,3]}.
  \end{equation}

  From the first item of Proposition~\ref{prop on principal operator lipschitz estimate}, we obtain uniform equicontinuous bounds for the sequence $\tcb{\Phi_n}_{n\in\N}\subset\mathcal{B}^0$; more precisely, the estimate
  \begin{equation}\label{Ascoli-Arzela_is_life}
    \sup_{n\in\N}\tnorm{\Phi_n}_{\m{LL}}<\infty
  \end{equation} 
  holds. We now make a call to the theorem of Arzel\`{a}--Ascoli: the sequence $\tcb{\Phi_n}_{n\in\N}\subset\mathcal{B}^0$ is, due to~\eqref{Ascoli-Arzela_is_life}, precompact in $\mathcal{B}^0$. After extraction of a suitable subsequence and relabeling, we are assured the existence of $\Phi\in\mathcal{B}^0$ such that 
  \begin{equation}\label{convergence_from_compactness}
    \tnorm{\Phi_n - \Phi}_{\mathcal{B}^0} \to 0\text{ as }n\to\infty.
  \end{equation}
  Returning with this information to our original convergences of equation~\eqref{___absurd_things___}, we deduce further that 
  \begin{equation}\label{_it_is_deduction_}
    \tnorm{\varphi_n + \Phi}_{\mathcal{B}^0}\to0\text{ as }n\to\infty.
  \end{equation}
  As it is the case that for every $n\in\N$ we have $\tnorm{\varphi_n}_{\mathcal{B}^0} = 1$, the strong limit~\eqref{_it_is_deduction_} implies that $\Phi$ is necessarily nontrivial; in fact $\tnorm{\Phi}_{\mathcal{B}^0} = 1$.

  We now endeavor to pass to the limit $n\to\infty$ in the left hand identity of~\eqref{the_definition_of_capital_Phi}. To assist in doing so, we set the following notation for the $\R^5$ Newtonian potential:
  \begin{equation}\label{R5_Newtonian_Potential}
    \bf{N}:\R^5\setminus\tcb{0}\to\R,\;\bf{N}\tp{x} = -\f{1}{8\pi^2}\f{1}{\tabs{x}^3},\;x\in\R^5\setminus\tcb{0}.
  \end{equation}
  Now, for an arbitrary test function $\psi\in C^0_{\m{c}}\tp{B_{\R^5}(0,3)}$ we take, for every $n\in\N$, the $L^2$-inner product of the definition of $\Phi_n$ from~\eqref{the_definition_of_capital_Phi} with $\psi$. After using symmetry of convolution with $\bf{N}$ and the definition of $A_\ep\tp{\phi}$ from the second item of Proposition~\ref{prop on basic smoothness}, the equality
  \begin{equation}\label{intermediate_integral_identity}
    \tbr{\Phi_n,\psi}_{L^2\tp{B(0,3)}} = 15\int_0^1\;\tbr{\tp{\Gamma_{\ep_n}}'\tp{\v + \tau\phi^0_n + \tp{1 - \tau}\phi_n^1-1}\varphi_n,\bf{N}\ast\psi}_{L^2\tp{B(0,3)}}\m{d}\tau
  \end{equation}
  is unveiled.

  The integrand on the right hand side of identity~\eqref{intermediate_integral_identity} is, for every $n\in\N$, a continuous function in $\tau$; moreover, the same argument via the first item of Proposition~\ref{prop on principal operator lipschitz estimate} that gave us the uniform bounds~\eqref{Ascoli-Arzela_is_life} applies here as well and we deduce that these integrands are uniformly bounded over $n\in\N$. Next, Proposition~\ref{prop on limit identification}, applied with the test functions $\varphi_n\tp{\bf{N}\ast\psi}$, implies pointwise-in-$\tau$ convergence. Thus, the dominated convergence theorem applies, and we can justifiably pass to the limit on either side of identity~\eqref{intermediate_integral_identity}.
  
  The identification of the aforementioned limit on the left is done with the convergence of equation~\eqref{convergence_from_compactness}. For the right side, we instead use the convergence of equations~\eqref{the limit that we would like to show} and~\eqref{_it_is_deduction_}. The result is 
  \begin{equation}
    \tbr{\Phi,\psi}_{L^2\tp{B(0,3)}} = -5\int_{\S^4}\Phi\tp{y}\tp{\bf{N}\ast\psi}\tp{y}\;\m{d}\mathcal{H}^4\tp{y}.
  \end{equation}
  Upon appealing to the arbitrary nature of $\psi$, we see from the above that $\Phi$ solves
  \begin{equation}\label{the_important_integral_equation}
    \Phi(x) = \f{5}{8\pi^2}\int_{\S^4}\f{1}{\tabs{x - y}^3}\Phi\tp{y}\;\m{d}\mathcal{H}^4\tp{y}
  \end{equation}
  for all $x\in B_{\R^5}[0,3]$. Proposition~\ref{prop on nondegeneracy of Hills spherical vortex}, that affirms the nondegeneracy of Hill's vortex, now informs us that $\Phi(y) = 0$ for all $y\in\S^4$. Hence, from~\eqref{the_important_integral_equation}, $\Phi$ vanishes identically in $B_{\R^5}[0,3]$. This is in direct conflict with the fact that $\tnorm{\Phi}_{\mathcal{B}^0} = 1$. Therefore, our contradiction hypotheses, that the statement of the theorem does not hold, must be false.
\end{proof}

Our first consequence of the closed range estimate of Theorem~\ref{thm on quantitative closed range estimates} is pointwise invertibility of the Fr\'echet differential of the maps $F_\ep$ from the fifth item of Definition~\ref{defn on spaces and mappings}. The upshot is that the neighborhood upon which differential invertibility is assured is uniform in $\ep$, despite the lack of uniform control in the operator norms in estimates~\eqref{inverse_operator_norm} and~\eqref{inverse_lipschitz_norm}.

\begin{coro}[Uniform region of differential invertibility]\label{coro on uniform region of differential invertibility}
  Let $0<\tilde{\ep}_1\le1/2$ be the small parameter given to us by Theorem~\ref{thm on quantitative closed range estimates}. There exists $C\in\R^+$ with the property that for all $0<\ep\le\tilde{\ep}_1$ and $\phi,\tilde{\phi}\in\mathcal{B}[\tilde{\ep}_1]$ the derivative $DF_\ep\tp{\phi}:\mathcal{B}\to\mathcal{B}$ is invertible and we have the estimates 
  \begin{equation}\label{inverse_operator_norm}
    \tnorm{\tsb{DF_\ep\tp{\phi}}^{-1}\varphi}_{\mathcal{B}} + \tnorm{DF_\ep\tp{\phi}\varphi}_{\mathcal{B}}\le\tp{C/\ep}\tnorm{\varphi}_{\mathcal{B}}
  \end{equation}
  and
  \begin{equation}\label{inverse_lipschitz_norm}
    \tnorm{\tsb{DF_\ep\tp{\phi}}^{-1}\varphi - \tsb{DF_\ep\tp{\tilde{\phi}}}^{-1}\varphi}_{\mathcal{B}}\le\tp{C/\ep^3}\tnorm{\phi - \tilde{\phi}}_{\mathcal{B}^0}\tnorm{\varphi}_{\mathcal{B}^0}
  \end{equation}
  for all $\varphi\in\mathcal{B}$.
\end{coro}
\begin{proof}
  Let $0<\ep\le\tilde{\ep}_1$ and $\phi\in\mathcal{B}[\tilde{\ep}_1]$. The defining identity of $F_\ep$ in equation~\eqref{defn_principal_part_operator} in conjunction with the estimate~\eqref{the_claimed_bound_here} show that as a linear operator $DF_\ep\tp{\phi}\in\mathcal{L}\tp{\mathcal{B}}$ has the structure of a compact perturbation of the identity and thus its Fredholm index is zero. Therefore, to verify invertibility of $DF_\ep\tp{\phi}$ it is sufficient to observe that $\m{ker}DF_\ep\tp{\phi} = \tcb{0}$.

  Let $\varphi\in\m{ker}DF_\ep\tp{\phi}$. Then, by invoking estimate~\eqref{_all_hail_satan_} with $\phi^0 = \phi^1 = \phi$, in addition to recalling that $\varphi + A_\ep\tp{\phi}\varphi = DF_\ep\tp{\phi}\varphi$, we get that
  \begin{equation}
    \tnorm{\varphi}_{\mathcal{B}^0}\le C\tnorm{\varphi + A_\ep\tp{\phi}\varphi}_{\mathcal{B}^0} = 0.
  \end{equation}
  So the kernel is indeed trivial.

  Let us next verify estimate~\eqref{inverse_operator_norm}. The second item of Proposition~\ref{prop on basic smoothness} allows us to focus only on the term $\tnorm{\tsb{DF_\ep\tp{\phi}}^{-1}\varphi}_{\mathcal{B}}$. Initially, thanks to Theorem~\ref{thm on quantitative closed range estimates} again, the bound
  \begin{equation}\label{_on_the_}
    \tnorm{\tsb{DF_\ep\tp{\phi}}^{-1}\varphi}_{\mathcal{B}^0}\le C\tnorm{\varphi}_{\mathcal{B}^0}
  \end{equation}
  is acquired. On the other hand,
  \begin{equation}
    \tsb{DF_\ep\tp{\phi}}^{-1}\varphi = \varphi - A_\ep\tp{\phi}\tsb{DF_\ep\tp{\phi}}^{-1}\varphi.
  \end{equation}
  We take the norm in $\mathcal{B}$ of the above and use estimate~\eqref{the_bounds_right_here} to find that
  \begin{equation}\label{_job_market_}
    \tnorm{\tsb{DF_\ep\tp{\phi}}^{-1}\varphi}_{\mathcal{B}}\le C\tnorm{\varphi}_{\mathcal{B}} + \tp{C/\ep}\tnorm{\tsb{DF_\ep\tp{\phi}}^{-1}\varphi}_{\mathcal{B}^0}.
  \end{equation}
  Synthesizing~\eqref{_on_the_} and~\eqref{_job_market_} completes the proof of estimate~\eqref{inverse_operator_norm}.

  We now deduce estimate~\eqref{inverse_lipschitz_norm}. Let $\tilde{\phi}\in\mathcal{B}[\tilde{\ep}_1]$ as well. The identity 
  \begin{equation}
    \tsb{DF_\ep\tp{\phi}}^{-1}\varphi - \tsb{DF_\ep\tp{\tilde{\phi}}}^{-1}\varphi = -\tsb{DF_\ep\tp{\phi}}^{-1}\tp{A_\ep\tp{\phi} - A_\ep\tp{\tilde{\phi}}}\tsb{DF_\ep\tp{\tilde{\phi}}}^{-1}\varphi
  \end{equation}
  holds. Take the norm of the above expression in the space $\mathcal{B}$ and combine estimates~\eqref{inverse_operator_norm}, \eqref{the_bounds_right_here}, and~\eqref{_on_the_}.
\end{proof}

As the second consequence of Theorem~\ref{thm on quantitative closed range estimates}, reverse modulus of continuity estimates for the maps $F_\ep$ that are $\ep$-uniform are bestowed upon us.

\begin{coro}[Uniform reverse logarithmic-Lipschitz estimates]\label{coro on uniform reverse log Lipschitz estimates}
  Let $0<\tilde{\ep}_1\le1/2$ be the small parameter granted to us by Theorem~\ref{thm on quantitative closed range estimates}. There exists a constant $C\in\R^+$ with the property that for all $0<\ep\le\tilde{\ep}_1$ and $\phi,\tilde{\phi}\in\mathcal{B}[\tilde{\ep}_1]$ we have the estimate
  \begin{equation}\label{rev_uni_cont_est}
    \tnorm{\phi - \tilde{\phi}}_{\mathcal{B}}\le\tnorm{F_\ep\tp{\phi} - F_\ep\tp{\tilde{\phi}}}_{\mathcal{B}}+C\tnorm{F_\ep\tp{\phi} - F_\ep\tp{\tilde{\phi}}}_{\mathcal{B}^0}\tp{1 + \tabs{\log\tnorm{F_\ep\tp{\phi} - F_\ep\tp{\tilde{\phi}}}_{\mathcal{B}^0}}}.
  \end{equation}
\end{coro}
\begin{proof}
  Let $\phi,\tilde{\phi}\in\mathcal{B}[\tilde{\ep}_1]$. We shall first employ the second item of Proposition~\ref{prop on basic smoothness} and the fundamental theorem of calculus to justifiably equate
  \begin{equation}\label{FTC_is_not_your_theorem}
    F_\ep\tp{\phi} - F_\ep\tp{\tilde{\phi}} = \tp{\phi - \tilde{\phi}} + \bsb{\int_0^1A_\ep\tp{\tau\phi + \tp{1 - \tau}\tilde{\phi}}\;\m{d}\tau}\tp{\phi - \tilde{\phi}}.
  \end{equation}
  We then take the norm in $\mathcal{B}^0$ of~\eqref{FTC_is_not_your_theorem} and employ Theorem~\ref{thm on quantitative closed range estimates} with $\phi^0 = \phi$, $\phi^1 = \tilde{\phi}$, and $\varphi = \phi - \tilde{\phi}$. As a result, the initial estimate
  \begin{equation}\label{as_a_result_there}
    \tnorm{\phi - \tilde{\phi}}_{\mathcal{B}^0}\le C\tnorm{F_\ep\tp{\phi} - F_{\ep}\tp{\tilde{\phi}}}_{\mathcal{B}^0}
  \end{equation}
  is true. On the other hand, directly from the definition~\eqref{defn_principal_part_operator}, we have 
  \begin{equation}\label{FTC_is_not_your_friend}
    \phi - \tilde{\phi} = F_\ep\tp{\phi} - F_\ep\tp{\tilde{\phi}} - \tp{K_\ep\tp{\phi} - K_\ep\tp{\tilde{\phi}}}.
  \end{equation}
  Now we take the norm in $\mathcal{B}$ of~\eqref{FTC_is_not_your_friend} and use the logarithmic-Lipschitz estimate of Corollary~\ref{coro on principal log-lip}. The bound 
  \begin{equation}\label{FTC_is_also_not_food}
    \tnorm{\phi - \tilde{\phi}}_{\mathcal{B}}\le\tnorm{F_\ep\tp{\phi} - F_\ep\tp{\tilde{\phi}}}_{\mathcal{B}} + C\tnorm{\phi - \tilde{\phi}}_{\mathcal{B}^0}\tp{1 + \tabs{\log\tnorm{\phi - \tilde{\phi}}_{\mathcal{B}^0}}}.
  \end{equation}
  is obtained via this process.

  The desired estimate~\eqref{rev_uni_cont_est} is established upon combining~\eqref{as_a_result_there} with~\eqref{FTC_is_also_not_food}.
\end{proof}

To close this subsection, we now combine Corollaries~\ref{coro on uniform region of differential invertibility} and~\ref{coro on uniform reverse log Lipschitz estimates} with the abstract quantitative local surjectivity result developed in Theorem 3.13 of Radu and Stevenson~\cite{VPD}.

\begin{thm}[Uniform existence of local inverse functions]\label{thm on uniform existence of local inverse functions}
  Let $0<\tilde{\ep}_1\le1/2$ be the positive parameter granted to us by Theorem~\ref{thm on quantitative closed range estimates}. There exists $0<\Bar{\ep}_1\le\tilde{\ep}_1$ and $C\in\R^+$ such that for all $0<\ep\le\tilde{\ep}_1$ the following hold.
  \begin{enumerate}
    \item For each $f\in\mathcal{B}[\Bar{\ep}_1]$, there exists a unique $\phi\in\mathcal{B}[\tilde{\ep}_1]$ such that $F_\ep\tp{\phi} = f$. This defines an inverse map $\tp{F_\ep}^{-1}:\mathcal{B}[\Bar{\ep}_1]\to\mathcal{B}[\tilde{\ep}_1]$ mapping $f\mapsto\phi$.
    \item The inverse maps from the previous item obey $\tp{F_\ep}^{-1}\tp{0} = 0$ and the uniform logarithmic-Lipschitz estimates
    \begin{equation}
      \tnorm{\tp{F_\ep}^{-1}\tp{f} - \tp{F_\ep}^{-1}\tp{\tilde{f}}}_{\mathcal{B}}\le\tnorm{f - \tilde{f}}_{\mathcal{B}} + C\tnorm{f - \tilde{f}}_{\mathcal{B}^0}\tp{1 + \tabs{\log\tnorm{f - \tilde{f}}_{\mathcal{B}^0}}}
    \end{equation}
    for all $f,\tilde{f}\in\mathcal{B}[\Bar{\ep}_1]$.
  \end{enumerate}
\end{thm}
\begin{proof}
  The proof is extremely similar to that of Corollary 3.14 in~\cite{VPD}; as such we only sketch the details. As previously mentioned, the main engine for establishing the local inverse maps of the first item is the quantitative local surjectivity theorem, stated precisely in Theorem 3.13 of~\cite{VPD}. To satisfy the hypotheses of this abstract tool, we essentially only need the existence of a right inverse to the Fr\'echet differential and a nonlinear a priori estimate. These are provided, in a suitably $\ep$-uniform fashion, by Corollaries~\ref{coro on uniform region of differential invertibility} and~\ref{coro on uniform reverse log Lipschitz estimates}, respectively. Armed with the existence of the maps $\tp{F_\ep}^{-1}:\mathcal{B}[\Bar{\ep}_1]\to\mathcal{B}[\tilde{\ep}_1]$, the second item is an immediate consequence of the reverse logarithmic-Lipschitz estimate~\eqref{rev_uni_cont_est}.
\end{proof}

\section{Existence of smooth vortex rings}\label{SO_construction_smooth_ring}

We now put the pieces together to construct smooth vortex rings. The inverse maps from Section~\ref{SO_uni_inv_near_Hill} convert the nonlinear problem into a compact fixed point equation, to which we apply Schauder's theorem after a suitable selection of parameters. We then extend the resulting solutions to $\R^5$ and record the properties needed in the Hill--Norbury desingularization limit.

\subsection{Construction of nearby solutions}\label{SS_nearby_construction}

\begin{lem}[Source term]\label{lem on estimates on the source term}
  There exists $C\in\R^+$ such that the following hold for all $0<\ep\le 1$.
  \begin{enumerate}
    \item The function 
    \begin{equation}\label{_defn_of_the_source_term_}
      f_\ep\tp{x} = \f{15}{8\pi^2}\int_{B_{\R^5}\tsb{0,3}}\f{1}{\tabs{x - y}^3}\tp{\Gamma_\ep\tp{\v\tp{y} - 1} - \mathds{1}_{\tp{0,\infty}}\tp{\v\tp{y} - 1}}\;\m{d}\mathcal{H}^5\tp{y}\text{ defined for }x\in B_{\R^5}[0,3]
    \end{equation}
    obeys the inclusions $f_\ep\in\mathcal{B}$ and $f_\ep\in\m{LL}^1\tp{B_{\R^5}[0,3]}$.
    \item The estimates 
    \begin{equation}\label{_estimates_on_the_source_term_}
      \tnorm{f_\ep}_{\m{LL}^1}\le C\text{ and }\tnorm{f_\ep}_{\mathcal{B}}\le C\ep\tp{1 + \tabs{\log\ep}}
    \end{equation}
    hold.
  \end{enumerate}
\end{lem}
\begin{proof}
  An argument similar to that of the proof of the first item of Proposition~\ref{prop on basic smoothness} verifies that $f_\ep\in\m{LL}^1\tp{B_{\R^5}[0,3]}$ -- in fact, we have the left hand estimate of~\eqref{_estimates_on_the_source_term_} with an $\ep$-uniform constant. We omit the routine details. It is similarly trivial to verify that $f_\ep$ is $\mathbb{O}$-invariant (see Definition~\ref{defn on spaces and mappings}) and hence $f_\ep\in\mathcal{B}$.

  Thus, it remains only to check the right hand estimate of equation~\eqref{_estimates_on_the_source_term_}. Initially, we note the pointwise estimate for $x\in B_{\R^5}[0,3]$
  \begin{equation}
    \tabs{f_\ep\tp{x}} + \tabs{\grad f_\ep\tp{x}}\le C\int_{B_{\R^5}[0,3]}\bp{\f{1}{\tabs{x - y}^3} + \f{1}{\tabs{x - y}^4}}\mathds{1}_{\tcb{|\v - 1|\le\ep}}\tp{y}\;\m{d}\mathcal{H}^5\tp{y}
  \end{equation}
  which follows from support considerations. The sought-after $\ep\tp{1 + \tabs{\log\ep}}$ bound is now a direct consequence of the first item of Corollary~\ref{bounds for integral operators on thin domains} when $\ep$ is sufficiently small. For the remaining (uniformly positive) values of $\ep$ in a compact subinterval of $(0,1]$, the same bound follows from the uniform $\m{LL}^1$ estimate already established.
\end{proof}

The statement of the following result involves parameters and operators from across the document: $\ep_2$ and $\ep_3$ are from Lemmas~\ref{lem on a decomposition and properties} and~\ref{lem on bounds on interior stream lines}, respectively; $\mu_0$ is from Proposition~\ref{prop on fixing the total meridional vorticity}; $\mathfrak{C}$ is from Lemma~\ref{lem on properties of the M I and J}; the maps $K_\ep$, $N_{\ep,k}$, and $F_\ep$ are from Definition~\ref{defn on spaces and mappings}, the functions $f_\ep$ are from Lemma~\ref{lem on estimates on the source term}; and the maps $M$, $I$, $J$, and $T$ are from Definition~\ref{defn of swirl and helicity operators}.

\begin{thm}[Existence of fixed points]\label{thm on existence of a fixed point}
  Set $\ep_4 = \min\tcb{\ep_2,\ep_3,1/54}>0$. There exists $\R^+\ni\del_0<1/16$ satisfying $\del_0<1/\mu_0$ with the property that for all $0<\ep\le\del_0$, $0\le k,\del\le\del_0$, $|\eta|\le\del_0$, and $\R\ni\mu\ge1/\del_0$ there exists $\phi\in\mathcal{B}[\ep_4]$, $-2\mathfrak{C}\le\lambda\le2\mathfrak{C}$, and $\R\ni\rho\ge22$ such that the equations
  \begin{equation}\label{the_fixed_point_to_write_home_about}
    \begin{cases}
      \phi + K_\ep\tp{\phi} = N_{\ep,k}\tp{\phi} + f_\ep + \del^2\tp{M^1_\rho\tp{\phi} + \lambda^2M^2_\rho\tp{\phi}},\\
      I^1_\rho\tp{\phi} + \lambda I^2_\rho\tp{\phi} + \del^2\tp{J^1_\rho\tp{\phi} + \lambda^3J^2_\rho\tp{\phi}} = \eta,\\
      T^1_\rho\tp{\phi} + \tabs{\lambda}T^2_\rho\tp{\phi} = \mu
    \end{cases}
  \end{equation}
  are satisfied.
\end{thm}
\begin{proof}
  We begin by making some initial smallness assumptions on the available parameters to ensure that the right hand side of the first equation in system~\eqref{the_fixed_point_to_write_home_about} always sits inside the ball $\mathcal{B}[\Bar{\ep}_1]$ and hence, thanks to Theorem~\ref{thm on uniform existence of local inverse functions}, belongs to the domains of the local inverse maps $\tp{F_\ep}^{-1}$. For any $\phi\in\mathcal{B}[\ep_4]$, $0\le\ep,k,\del\le 1$ with $\ep>0$, $|\lambda|\le 2\mathfrak{C}$, and $\R\ni\rho\ge 22$ the estimate
  \begin{equation}\label{need_to_make_a_small_source}
    \tnorm{N_{\ep,k}\tp{\phi} + f_\ep + \del^2\tp{M^1_\rho\tp{\phi} + \lambda^2 M^2_\rho\tp{\phi}}}_{\mathcal{B}}\le C\tp{\tp{\ep + k}^{1/6} + \ep\tp{1 + \tabs{\log\ep}} + \del^2}
  \end{equation}
  is valid thanks to the conclusions of Corollary~\ref{coro on more estimates on the perturbative operator}, Lemma~\ref{lem on estimates on the source term}, and Lemma~\ref{lem on properties of the M I and J}. The above bound can be made arbitrarily small if we constrain $\ep$, $k$, and $\del$. In particular, there exists $\qoppa_0\in(0,1]$ such that if $0\le\ep,k,\del\le\qoppa_0$, $\ep>0$, and $\phi\in\mathcal{B}[\ep_4]$, then
  \begin{equation}\label{has_to_be_about_this_small}
    \tnorm{N_{\ep,k}\tp{\phi} + f_\ep + \del^2\tp{M^1_\rho\tp{\phi} + \lambda^2 M^2_\rho\tp{\phi}}}_{\mathcal{B}}\le\Bar{\ep}_1/2.
  \end{equation}

  Next, let us recall the continuous map $\pmb{\varrho}$ constructed in Proposition~\ref{prop on fixing the total meridional vorticity}. We shall directly let this map make our choice of $\rho\ge 22$ so that the final equation of system~\eqref{the_fixed_point_to_write_home_about} is trivially satisfied.

  We now introduce for $0<\qoppa\le\min\tcb{1/e,\qoppa_0,\tilde{\ep}_1,1/16,1/\mu_0}$ the sets 
  \begin{equation}
    \mathsf{P}_\qoppa = \tcb{ \tp{\ep,k,\del,\eta,\mu}\in[0,\qoppa]^3\times[-\qoppa,\qoppa]\times[1/\qoppa,\infty)\;:\;\ep>0}\text{ and }\mathsf{K} = \mathcal{B}[\ep_4]\times[-2\mathfrak{C},2\mathfrak{C}]
  \end{equation}
  and the mapping
  \begin{equation}
    \mathsf{P}_\qoppa\times\mathsf{K}\ni\tp{\tp{\ep,k,\del,\eta,\mu},\tp{\phi,\lambda}}\mapsto\mathtt{G}_{\tp{\ep,k,\del,\eta,\mu}}\tp{\phi,\lambda}\in\mathcal{B}\times\R
  \end{equation}
  with the action
  \begin{equation}\label{_map_to_be_fixed}
    \mathtt{G}_{\mathsf{p}}\tp{\phi,\lambda} = \bpm \tp{F_\ep}^{-1}\tp{N_{\ep,k}\tp{\phi} + f_\ep + \del^2\tp{M^1_{\pmb{\varrho}\tp{\phi,|\lambda|,\mu}}\tp{\phi} + \lambda^2 M^2_{\pmb{\varrho}\tp{\phi,|\lambda|,\mu}}\tp{\phi}}} \\ \tp{1/I^2_{\pmb{\varrho}\tp{\phi,|\lambda|,\mu}}\tp{\phi}}\tp{\eta - I^1_{\pmb{\varrho}\tp{\phi,|\lambda|,\mu}}\tp{\phi} - \del^2\tp{J_{\pmb{\varrho}\tp{\phi,|\lambda|,\mu}}^1\tp{\phi} + \lambda^3J_{\pmb{\varrho}\tp{\phi,|\lambda|,\mu}}^2\tp{\phi}}}\epm
  \end{equation}
  for $\mathsf{p} = \tp{\ep,k,\del,\eta,\mu}\in\mathsf{P}_\qoppa$ and $\tp{\phi,\lambda}\in\mathsf{K}$. By the above discussion and the definition of $\qoppa_0$ (see~\eqref{has_to_be_about_this_small}), the maps of~\eqref{_map_to_be_fixed} are well-defined in the sense that the composition with $\tp{F_\ep}^{-1}$ in the first component is valid.

  It is also the case that for each fixed $\mathsf{p}\in\mathsf{P}_\qoppa$ the map
  \begin{equation}\label{the_map_who_is_to_be_fixed}
    \mathsf{K}\ni\tp{\phi,\lambda}\mapsto\mathtt{G}_{\mathsf{p}}\tp{\phi,\lambda}\in\mathcal{B}\times\R
  \end{equation}
  is, thanks to Theorem~\ref{thm on uniform existence of local inverse functions}, Proposition~\ref{prop on basic smoothness}, Lemma~\ref{lem on properties of the M I and J}, Proposition~\ref{prop on fixing the total meridional vorticity}, and the Arzel\`{a}-Ascoli theorem,  both continuous and compact.

  Finally, we need only check now that for $\qoppa$ small enough, the map in~\eqref{the_map_who_is_to_be_fixed} is valued in a subset of $\mathsf{K}$. By combining estimate~\eqref{need_to_make_a_small_source}, the second item of Theorem~\ref{thm on uniform existence of local inverse functions}, and Lemma~\ref{lem on properties of the M I and J}, the estimates
  \begin{equation}\label{__00__}
    \tnorm{\tsb{\mathtt{G}_{\mathsf{p}}\tp{\phi,\lambda}}_1}_{\mathcal{B}}\le C\qoppa^{1/6}\log\tp{1/\qoppa}
  \end{equation}
  and 
  \begin{equation}\label{___00___}
    \tabs{\tsb{\mathtt{G}_{\mathsf{p}}\tp{\phi,\lambda}}_2}\le\mathfrak{C} + C\qoppa,
  \end{equation}
  for all $\tp{\phi,\lambda}\in\mathsf{K}$ and $\mathsf{p}\in\mathsf{P}_\qoppa$, are justified. We are using the notation $\tsb{\cdot}_{i}$ for $i\in\tcb{1,2}$ to denote the first and second components of the map in~\eqref{_map_to_be_fixed}.
  
  So indeed, by inspection of~\eqref{__00__} and~\eqref{___00___}, one is assured that there exists $0<\qoppa_1\le\qoppa_0$ such that for all $\mathsf{p}\in\mathsf{P}_{\qoppa_1}$ the compact and continuous map~\eqref{the_map_who_is_to_be_fixed} has an image that is a subset of the closed and convex set $\mathsf{K}$. After decreasing $\qoppa_1$, if necessary, we may also assume that
  \begin{equation}
    \qoppa_1\le\min\tcb{1/e,\tilde{\ep}_1,1/16,1/\mu_0}.
  \end{equation}
  To conclude, we employ Schauder's fixed point theorem; for a modern reference see, e.g., Theorem 3.2 in Chapter 6 of Granas and Dugundji~\cite{MR1987179}. The map~\eqref{the_map_who_is_to_be_fixed} has at least one fixed point in $\mathsf{K}$. Upon unpacking the equations satisfied by the fixed point we obtain a solution to system~\eqref{the_fixed_point_to_write_home_about}. Taking $\del_0=\qoppa_1/2$ completes the proof.
\end{proof}

\subsection{Instantiation of a solution map}\label{SS_instantiation_soln_map}

In the next result we shall unpack the fixed points produced by Theorem~\ref{thm on existence of a fixed point} and obtain solutions to the $\R^5$-manifestation of Hicks equation~\cite{Hicks1898}. In the following result we use the bracket notation $\tbr{\cdot}:\R^d\to\R$, $\tbr{x} = \tp{1 + |x|^2}^{1/2}$.

\begin{coro}[Solutions to an $\R^5$ Hicks equation]\label{coro on solutions to the Bragg-Hawthorne equation}
  Let $\del_0\in\R^+$ be the small parameter granted to us by Theorem~\ref{thm on existence of a fixed point}. There exists $C\in\R^+$ such that for any $0\le\ep,k,\del\le\del_0$, $|\eta|\le\del_0$, with $\ep,k>0$ and $\R\ni\mu\ge1/\del_0$, any solution $\phi\in\mathcal{B}[\ep_4]$, $|\lambda|\le2\mathfrak{C}$, and $\R\ni\rho\ge22$ to the system of equations~\eqref{the_fixed_point_to_write_home_about} produced by the aforementioned theorem additionally satisfies the following items.
  \begin{enumerate}
    \item The circulation function $S$ and the head derivative function $\Gamma$ defined via
    \begin{equation}\label{circulation_and_head_derivative_functions_}
      S(t) = \del\Sigma_\rho\tp{t - 1/8 + k} + \del\lambda\Sigma_\rho\tp{t - 1/4 + k}\text{ and }\Gamma\tp{t} = 15\cdot\Gamma_\ep\tp{t}\text{ for }t\in\R
    \end{equation}
    satisfy $S,\Gamma\in C^\infty\tp{\R}$, $\supp\Gamma\subseteq[0,\infty)$, $\supp S\subseteq[1/16,5/16]$.
    \item There exists $\Phi\in C^\infty\tp{\R^5}$ such that 
    \begin{equation}\label{extension_property_and_bounds}
      \Phi = \v + \phi\text{ in }B_{\R^5}[0,3],\;\sup\tcb{\Phi(y)\;:\;y\in\R^5\setminus B_{\R^5}(0,3)}<1,\;\forall\;T\in\O,\;y\in\R^5,\;\Phi\tp{T y} = \Phi\tp{y},
    \end{equation}
    and Hicks equation
    \begin{equation}\label{bragg_hawthorne_r5}
      \tabs{\P y}^2\tp{\Delta\Phi}\tp{y} + \tp{S\cdot S'}\tp{\tabs{\P y}^2\tp{\Phi\tp{y} - 1} - k} +  \tabs{\P y}^2\Gamma\tp{\tabs{\P y}^2\tp{\Phi\tp{y} - 1} - k} = 0\text{ for }y\in\R^5
    \end{equation}
    is classically satisfied.
    \item The function $\Phi$ has the helicity
    \begin{equation}\label{known_helicity}
      \del\eta = \f{2}{\pi}\int_{\R^5}\bp{\f{1}{\tabs{\P y}^2}\tp{S\cdot\Gamma}\tp{\tabs{\P y}^2\tp{\Phi\tp{y} - 1} - k} + \f{1}{\tabs{\P y}^4}\tp{S\cdot S\cdot S'}\tp{\tabs{\P y}^2\tp{\Phi\tp{y} - 1} - k}}\;\m{d}\mathcal{H}^5\tp{y}
    \end{equation}
    and the total meridional vorticity
    \begin{equation}\label{known_total_meridional_vorticity}
      \del\mu = \f{1}{\pi}\int_{\R^5}\f{1}{\tabs{\P y}^3}\tabs{S'\tp{\tabs{\P y}^2\tp{\Phi\tp{y} - 1} - k}}\tabs{\grad\tp{\tabs{\P y}^2\tp{\Phi\tp{y} - 1}}}\;\m{d}\mathcal{H}^5\tp{y}.
    \end{equation}
    \item The estimates
    \begin{multline}\label{hill_deviation_estimates}
      \tnorm{\tbr{\cdot}^3\Phi}_{\m{LL}^1\tp{\R^5}}\le C,\;\tnorm{\Phi}_{H^2\tp{\R^5}}\le C,\text{ and }\\\tnorm{\tbr{\cdot}^3\tp{\Phi - \v}}_{L^\infty\tp{\R^5}} + \tnorm{\tbr{\cdot}^4\tp{\grad\Phi - \grad\v}}_{L^\infty\tp{\R^5}}\le C\tp{\tp{\ep + k}^{1/6} + \del^2}\tp{1 + \tabs{\log\tp{\tp{\ep + k}^{1/6} + \del^2}}}
    \end{multline}
    hold.
    \item The qualitative properties
    \begin{multline}\label{_key_qualitative_properties_}
      \tsb{\Phi(y) = \Phi(T y)>0\;\forall\;y\in\R^5\text{ and }T\in\O},\;\tsb{\min\tcb{\Phi\tp{y}\;:\;y\in B_{\R^5}[0,3]}\ge1/54},\\\tsb{\forall\;y\in\R^5\text{ with }y_5>0,\;\pd_5\Phi(y)<0},\text{ and }\supp\Delta\Phi\subset\tcb{y\in\R^5\;:\;|y|\le 2\text{ and }2\tabs{\P y}^2\ge k}.
    \end{multline}
    are satisfied.
  \end{enumerate}
\end{coro}
\begin{proof}
  Throughout the proof we shall let $1\le C<\infty$ denote a constant, whose precise value may change from line to line, but it only depends on the parameters $\ep_4$ and $\mathfrak{C}$.
  
  We fix parameters $0\le\ep,k,\del\le\del_0$, $|\eta|\le\del_0$, with $\ep,k>0$ and $\R\ni\mu\ge1/\del_0$ and let $\phi\in\mathcal{B}[\ep_4]$, $|\lambda|\le 2\mathfrak{C}$, and $\R\ni\rho\ge 22$ be a corresponding solution to system~\eqref{the_fixed_point_to_write_home_about} given by Theorem~\ref{thm on existence of a fixed point}.
  
  Upon defining the functions $S$ and $\Gamma$ as in~\eqref{circulation_and_head_derivative_functions_}, the assertions of the first item are immediate given the properties enumerated in equations~\eqref{reg_vort_func}, \eqref{defn_regularized_heaviside}, \eqref{_the_comb_functions_}, and~\eqref{_remarks_on_the_comb_functions_}.

  Now we shall analyze an auxiliary function. Let us define 
  \begin{equation}
    \Omega\tp{y} = \Gamma\tp{\tabs{\P y}^2\tp{\v\tp{y} + \phi\tp{y} - 1} - k} + \f{1}{\tabs{\P y}^2}\tp{S\cdot S'}\tp{\tabs{\P y}^2\tp{\v\tp{y} + \phi\tp{y} - 1} - k}\text{ for }y\in B_{\R^5}[0,3].
  \end{equation}
 The second term above is understood through its smooth extension by zero across $\tcb{y\in\R^5\;:\;\P y=0}$; indeed, it vanishes in a neighborhood of this set since $k>0$ and $\supp S\subseteq[1/16,5/16]$.
  The second item of Lemma~\ref{lem on bounds on interior stream lines} when combined with the support considerations of~\eqref{_remarks_on_the_comb_functions_} implies that
  \begin{equation}
    1/C\le\tabs{\P y}\le C\text{ for all }y\in B_{\R^5}[0,3]\text{ such that }S\tp{\tabs{\P y}^2\tp{\v\tp{y} + \phi\tp{y} - 1} - k}\neq0,
  \end{equation}
  $\Omega\in\mathcal{B}$, and we have the uniform estimate $\tnorm{\Omega}_{L^\infty\tp{B\tp{0,3}}}\le C$. In fact, further support considerations in conjunction with the second items of Lemmas~\ref{lem on simple observations} and~\ref{lem on a decomposition and properties} show that $\supp\Omega\subseteq B_{\R^5}[0,2]$. We are therefore justified in extending $\Omega$ trivially by zero and regarding it as a continuously differentiable function defined on the whole of $\R^5$.

  To define the desired function $\Phi$, we set
  \begin{equation}\label{definition_of_the_function_Phi}
    \Phi(x) = \f{1}{8\pi^2}\int_{\R^5}\f{1}{\tabs{x - y}^3}\Omega(y)\;\m{d}\mathcal{H}^5\tp{y}\text{ for all }x\in\R^5.
  \end{equation}
  Several properties of $\Phi$ are immediate. First, since $\Omega$ is $\O$-invariant, it is clear that $\Phi(Tx) = \Phi(x)$ for all $x\in\R^5$ and all $T\in\O$. Second, the uniform estimates $\tnorm{\tbr{\cdot}^3\Phi}_{\m{LL}^1\tp{\R^5}}\le C$ and $\norm{\Phi}_{H^2\tp{\R^5}}\le C$ are simple consequences of the uniform bound $\tnorm{\Omega}_{L^\infty\tp{\R^5}}\le C$ and the uniform support constraint $\supp\Omega\subseteq B_{\R^5}[0,2]$.

  To see that for all $y\in B_{\R^5}[0,3]$ we have the equality $\Phi(y) = \v\tp{y} + \phi(y)$, we take the first equation in system~\eqref{the_fixed_point_to_write_home_about} and combine this information with the identity~\eqref{Hills_Integral_Equation} that is satisfied by Hill's solution $\v$. Simple algebra and definition unpacking reveals that for all $x\in B_{\R^5}[0,3]$
  \begin{equation}\label{WE_ARE_HERE}
    \f{1}{8\pi^2}\int_{B_{\R^5}[0,3]}\f{1}{\tabs{x - y}^3}\Omega\tp{y}\;\m{d}\mathcal{H}^5\tp{y} = \v\tp{x} + \phi\tp{x}.
  \end{equation}
  Thus, the combination of~\eqref{WE_ARE_HERE} with $\supp\Omega\subseteq B_{\R^5}[0,2]$ and the definition~\eqref{definition_of_the_function_Phi} reveals $\Phi = \v + \phi$ on $B_{\R^5}[0,3]$, as desired.

  Working towards the global description of $\Phi$ claimed in~\eqref{bragg_hawthorne_r5}, we require the intermediate fact:
  \begin{equation}\label{an_intermediate fact}
    \sup\tcb{\Phi(y)\;:\;y\in\R^5\setminus B_{\R^5}\tp{0,3}}<1.
  \end{equation}
  Fortunately, this is a simple corollary of the maximum principle for harmonic functions. A particular consequence of the second item of Lemma~\ref{lem on a decomposition and properties} is the uniform bound
  \begin{equation}
    \v\tp{y} + \phi\tp{y}\le 1 - 1/C\text{ for all }y\in B_{\R^5}[0,3]\text{ satisfying }|y|\ge3/2.
  \end{equation}
  From the agreement of $\Phi$ and $\v + \phi$ on $B_{\R^5}[0,3]$ we therefore have 
  \begin{equation}
    \max\tcb{\Phi(y)\;:\;y\in\R^5,\;|y| = 3}\le 1 - 1/C.
  \end{equation}
  On the other hand, it was already noted that $\sup\tcb{\tbr{y}^3|\Phi\tp{y}|\;:\;y\in\R^5}\le C$. Thus there is some radius $R>3$ such that 
  \begin{equation}\label{__OOOOO____}
    \sup\tcb{\Phi\tp{y}\;:\;y\in\R^5,\;|y|\ge R}\le\tp{1 - 1/C}/2.
  \end{equation}
  Since $\Delta\Phi\tp{y} = 0$ for all $y\in\R^5$ with $|y|\ge 3$, we conclude that~\eqref{an_intermediate fact} indeed holds as a consequence of~\eqref{__OOOOO____} and an invocation of the maximum principle in the annular region $\tcb{y\in\R^5\;:\;3\le |y|\le R}$.

  The key consequence of inequality~\eqref{an_intermediate fact} is that we must have 
  \begin{equation}
    \Gamma\tp{\tabs{\P y}^2\tp{\Phi\tp{y} - 1} - k} = 0\text{ and }S\tp{\tp{\P y}^2\tp{\Phi\tp{y} - 1} - k}=0\text{ for all }y\in\R^5\setminus B_{\R^5}\tp{0,3}.
  \end{equation}
  Hence we deduce that
  \begin{equation}\label{we_have_now_extended_omega_globally}
    \Omega\tp{y} = \Gamma\tp{\tabs{\P y}^2\tp{\Phi\tp{y} - 1} - k} + \f{1}{\tabs{\P y}^2}\tp{S\cdot S'}\tp{\tabs{\P y}^2\tp{\Phi\tp{y} - 1} - k}\text{ for all }y\in\R^5.
  \end{equation}
  Equality~\eqref{bragg_hawthorne_r5} is justified by combining~\eqref{definition_of_the_function_Phi} and~\eqref{we_have_now_extended_omega_globally}.

  To complete the proof of the second item, we must show that $\Phi\in C^\infty\tp{\R^5}$. This, however, is a consequence of standard interior elliptic regularity estimates for the Laplacian and a trivial induction argument. Indeed, if for some $m\in\N$ we have that $\Phi\in C^m_{\loc}\tp{\R^5}$, then since $\Omega$ is defined in terms of $\Phi$ via~\eqref{we_have_now_extended_omega_globally} we must have $\Omega\in C^m_{\m{c}}\tp{\R^5}$. But then a very crude regularity gain from the formula~\eqref{definition_of_the_function_Phi} implies that actually $\Phi\in C^{m+1}_{\m{loc}}\tp{\R^5}$.

  We also record the following useful consequence of~\eqref{we_have_now_extended_omega_globally}. If $y\in\supp\Omega$, then necessarily
  \begin{equation}
    \tabs{\P y}^2\tp{\Phi\tp{y}-1}-k\ge0.
  \end{equation}
  Since $\supp\Omega\subseteq B_{\R^5}[0,2]$ and $\Phi=\v+\phi$ on this ball, we have $\Phi-1\le 3/2 + 1/54<2$. Hence
  \begin{equation}\label{support_away_from_axis}
    \supp\Omega\subseteq B_{\R^5}[0,2]\cap\tcb{y\in\R^5\;:\;2\tabs{\P y}^2\ge k}.
  \end{equation}
  In particular, since $-\Delta\Phi=\Omega$, the same inclusion holds with $\supp\Delta\Phi$ in place of $\supp\Omega$.

  Let us now consider the third item. Multiplication of the second equation of system~\eqref{the_fixed_point_to_write_home_about} by $\del$ followed by unpacking the various definitions involved and using that $\Gamma = 15$ on the support of $S$ and the two summands defining $S$ have disjoint supports, $\Phi = \v + \phi$ on $B_{\R^5}[0,3]$, and $\Phi<1$ on $\R^5\setminus B_{\R^5}\tp{0,3}$ will reveal~\eqref{known_helicity} to be a valid identity. An identical strategy will verify equation~\eqref{known_total_meridional_vorticity} as well.

  Next, we focus on the estimates asserted in the fourth item. We have already shown the two left hand inequalities in~\eqref{hill_deviation_estimates} to hold. The proof of the remaining two estimates shall proceed in two parts. Firstly, we derive the part of the estimate that holds on the ball $B_{\R^5}\tsb{0,3}$. On this set we recall that we have the equality $\Phi - \v = \phi$ and thus it suffices to estimate the norm of $\phi$. To do this, we can look through the proof of Theorem~\ref{thm on existence of a fixed point}. Specifically, we have estimate~\eqref{need_to_make_a_small_source} and the fixed point identity 
  \begin{equation}
    \phi = \tp{F_\ep}^{-1}\tp{N_{\ep,k}\tp{\phi} + f_\ep + \del^2\tp{M^1_\rho\tp{\phi} + \lambda^2M^2_\rho\tp{\phi}}}.
  \end{equation}
  Combining these with estimate~\eqref{rev_uni_cont_est} gives us the auxiliary bound 
  \begin{equation}\label{_the_initial_estimate_}
    \tnorm{\Phi - \v}_{C^1_{\m{b}}\tp{B[0,3]}}=\tnorm{\phi}_{C^1_{\m{b}}\tp{B[0,3]}}\le C\tp{\tp{\ep + k}^{1/6} + \del^2}\tp{1 + \tabs{\log\tp{\tp{\ep + k}^{1/6} + \del^2}}}.
  \end{equation}

  The deduction of the weighted estimate in the complement of the ball will now follow from a simple promotion argument employing harmonic function theory. The following facts are essential ingredients:
  \begin{equation}
    \Delta\tp{\Phi - \v}\tp{y} = 0\text{ for all }y\in\R^5\setminus B_{\R^5}\tp{0,2}\text{ and }\sup\tcb{\tbr{y}^3\tabs{\Phi\tp{y} - \v\tp{y}}\;:\;y\in\R^5}\le C<\infty.
  \end{equation}
  We are thus in a position to perform a maximum principle argument, similar to the one previously appearing in this proof. For any $\al\in\R^+$, the initial estimate~\eqref{_the_initial_estimate_} gives
  \begin{equation}
    \max\tcb{\Phi(y) - \v\tp{y}-\al|y|^{-3},\v\tp{y} - \Phi\tp{y} - \al|y|^{-3}}\le C\tp{\tp{\ep + k}^{1/6} + \del^2}\tp{1 + \tabs{\log\tp{\tp{\ep + k}^{1/6} + \del^2}}} - \al/27
  \end{equation}
  for all $y\in\pd B_{\R^5}\tp{0,3}$. The choice $\al = 27\cdot C\tp{\tp{\ep + k}^{1/6} + \del^2}\tp{1 + \tabs{\log\tp{\tp{\ep + k}^{1/6} + \del^2}}}$ makes the right hand side nonpositive. Since the two functions inside the maximum are harmonic in $\R^5\setminus B_{\R^5}[0,3]$ and tend to zero at infinity, the maximum principle on the exterior domain gives us the bound
  \begin{equation}\label{_initial_bound_of_3_decay_}
    \tnorm{\tbr{\cdot}^3\tp{\Phi - \v}}_{L^\infty\tp{\R^5}}\le C \tp{\tp{\ep + k}^{1/6} + \del^2}\tp{1 + \tabs{\log\tp{\tp{\ep + k}^{1/6} + \del^2}}}.
  \end{equation}

  The remaining bound in~\eqref{hill_deviation_estimates} is that of the derivatives. For this, we shall use~\eqref{_initial_bound_of_3_decay_} and the mean value property of harmonic functions. For any $x\in\R^5\setminus B_{\R^5}\tp{0,3}$ the formula
  \begin{equation}\label{differentiated_mvf}
    \grad\tp{\Phi - \v}\tp{x} = \f{15}{8\pi^2}\bp{\f{3}{|x|}}^5\int_{\pd B(x,|x|/3)}\f{y - x}{|y - x|}\tp{\Phi - \v}\tp{y}\;\m{d}\mathcal{H}^4\tp{y}.
  \end{equation}
  holds. We then directly estimate the expression~\eqref{differentiated_mvf} and bound the integrand using the already established bound~\eqref{_initial_bound_of_3_decay_}:
  \begin{equation}\label{final_bound_of_4_decay}
    \f{\tabs{\grad\tp{\Phi - \v}\tp{x}}}{\tp{\tp{\ep + k}^{1/6} + \del^2}\tp{1 + \tabs{\log\tp{\tp{\ep + k}^{1/6} + \del^2}}}}\le \f{C}{|x|^5}\f{\mathcal{H}^{4}\tp{\pd B(x,|x|/3)}}{\tbr{x}^3}\le\f{C}{|x|^4}.
  \end{equation}

  Lastly, we consider the fifth item. The $\O$-invariance of $\Phi$ is immediate from Definition~\ref{defn on spaces and mappings}. Another maximum principle argument will show that $\Phi>0$ on all of $\R^5$.
  
  Since $\ep_4\le1/54$ we know immediately from the third item of Lemma~\ref{lem on simple observations} that $\Phi(y)\ge1/54$ for all $y\in B_{\R^5}[0,3]$. On the other hand, we know that $\Delta\Phi(y) = 0$ for all $|y|\ge 2$ and $\Phi(y)\to0$ as $|y|\to\infty$. Therefore, if $\epsilon\in\tp{0,1}$ is arbitrarily fixed, then there exists some $\uprho>3$ such that $\min\tcb{\Phi(y)\;:\;y\in\pd B_{\R^5}[0,\tilde{\uprho}]}\ge-\epsilon$ for all $\tilde{\uprho}\ge\uprho$. By the maximum principle, we deduce that $\Phi(y)\ge-\epsilon$ for all $3\le|y|\le\tilde{\uprho}$. Letting first $\tilde{\uprho}\to\infty$ and then $\epsilon\to0$ shows that $\Phi(y)\ge0$ for all $y\in\R^5$.
  
  If it were the case that for some $y\in\R^5\setminus B_{\R^5}[0,3]$ we had $\Phi(y) = 0$, then we can apply the strong maximum principle in the ball $B_{\R^5}[y,\epsilon]$ for some small $\epsilon>0$ such that $\Delta\Phi = 0$ on this ball. We would then conclude that $\Phi(z) = 0$ for all $|z - y|\le\epsilon$. This is absurd, since $\Phi$ is harmonic (and hence real analytic) and not identically zero on the connected set $\R^5\setminus B_{\R^5}[0,2]$. 

It remains only to establish the monotonicity assertion in~\eqref{_key_qualitative_properties_}. We invoke the one-directional entire-space moving plane result stated as Theorem~$4'$ in Gidas, Ni, and Nirenberg~\cite{MR544879}. Indeed, after dividing~\eqref{bragg_hawthorne_r5} by $\tabs{\P y}^2$ and understanding the swirl term through its smooth extension across $\tcb{y\in\R^5\;:\;\P y=0}$, we obtain an equation of the form~(4.17) in that paper, with $y_5$ as the distinguished variable: the resulting nonlinearity is smooth and independent of $y_5$. Moreover, the representation~\eqref{definition_of_the_function_Phi} and the compact support of $\Omega$ give the asymptotic expansion required there, with positive leading coefficient since $\Phi>0$. It follows that the function $\Phi$ is symmetric about some hyperplane perpendicular to $e_5$ and satisfies $\pd_5\Phi<0$ on one side of this hyperplane. Since $\Phi$ is even in $y_5$ and vanishes at infinity, this hyperplane must be $\tcb{y\in\R^5\;:\;y_5 = 0}$ and hence $\pd_5\Phi(\cdot,y_5)<0$ for $y_5>0$.
\end{proof}

\begin{rmk}
  The qualitative properties of the functions $\Phi$ from the fifth item of Corollary~\ref{coro on solutions to the Bragg-Hawthorne equation} are crucially used in the sequel when we pass to the singular limit $\ep\to0$. We note that similar symmetry and positivity properties have been a recurrent theme in the analysis of steady vortex rings; to name just two examples, see Section 4 in Ni~\cite{MR583638} and Section 3 in Amick and Turner~\cite{MR929976}.
\end{rmk}

\begin{rmk}[Instantiation of a solution map]\label{remark on instantiation of a solution map}
  Let us denote the initial admissible parameter space by 
  \begin{equation}\label{defn_admissible_parameter_space}
    \mathsf{P} = (0,\del_0]^2\times[0,\del_0]\times[-\del_0,\del_0]\times[1/\del_0,\infty),
  \end{equation}
  where $\del_0$ is the small positive parameter granted to us by Theorem~\ref{thm on existence of a fixed point}. Thanks to this theorem, Corollary~\ref{coro on solutions to the Bragg-Hawthorne equation}, and the axiom of choice, there exists some function 
  \begin{equation}\label{choice_of_solutions}
    \mathsf{P}\ni\mathsf{p} = \tp{\ep,k,\del,\eta,\mu}\mapsto\Phi = \Upphi_{\mathsf{p}}\in C^\infty\tp{\R^5}
  \end{equation}
  where $\Phi$ is the map described in the second item of Corollary~\ref{coro on solutions to the Bragg-Hawthorne equation}; in particular, Hicks equation~\eqref{bragg_hawthorne_r5} is satisfied; the helicity and total meridional vorticity identities, equations~\eqref{known_helicity} and~\eqref{known_total_meridional_vorticity}, respectively, hold; and the estimates~\eqref{hill_deviation_estimates} and qualitative properties~\eqref{_key_qualitative_properties_} are true.
\end{rmk}

\section{Euler reconstruction and vortex geometry near Hill's solution}\label{SO_euler_recon_and_vortex_geometry}

We now return from the $\R^5$ formulation to the physical $\R^3$-Euler equations. We first use the Chandrasekhar transformation to reconstruct axisymmetric traveling waves from the solutions of Section~\ref{SO_construction_smooth_ring} and then establish uniform geometric and velocity estimates for reconstructions near Hill's vortex.

\subsection{The Chandrasekhar transformation and Euler reconstruction}\label{SS_chandra_trans_Euler_recon}

Recall the subset of orthogonal transformations $\O$ introduced in the first item of Definition~\ref{defn on spaces and mappings}. For a function $\Phi:\R^5\to\R$ that is $\O$-invariant, meaning $\Phi(Tx) = \Phi(x)$ for all $x\in\R^5$ and $T\in\O$, we define its meridional restriction $\mathfrak{R}\Phi:\Pi\to\R$ by
\begin{equation}\label{defn_meridional_restriction}
  \tp{\mathfrak{R}\Phi}\tp{r,z} = \Phi(re_1 + ze_5)\text{ for all }\tp{r,z}\in\Pi.
\end{equation}
There is nothing particularly special about the choice of $e_1$, replacing this with any $\nu\in\S^4\cap\tcb{e_5}^\perp$ would give the same map $\mathfrak{R}$, thanks to $\O$-invariance.

A computation in cylindrical coordinates shows that every integrable $\O$-invariant function $G:\R^5\to\R$ satisfies
\begin{equation}\label{integral_identity_key}
  \int_{\R^5}G(y)\;\m{d}y = 2\pi^2\int_{0}^\infty\int_{\R}\tp{\mathfrak{R}G}\tp{r,z}r^3\;\m{d}r\;\m{d}z.
\end{equation}

Two more useful properties of the meridional restriction $\mathfrak{R}$ are recorded in the following elementary lemma, whose proof we omit.

\begin{lem}[The Chandrasekhar transformation]\label{lem on the chandrasekhar transformation}
  Let $\Phi\in C^\infty\tp{\R^5}$ be $\O$-invariant. The following hold.
  \begin{enumerate}
    \item We have 
    \begin{equation}\label{chandrasekhar_transformation_identity}
      \bp{r\pd_r\bp{\f{1}{r}\pd_r} + \pd_z^2}\tp{r^2\mathfrak{R}\Phi} = r^2\mathfrak{R}\tp{\Delta_{\R^5}\Phi}
    \end{equation}
    on $\Pi$.
    \item For $k>0$ define the map $\mathscr{Q}_{\Phi,k}\in C^\infty\tp{\R^5}$ by 
    \begin{equation}\label{this_thing_we_care_about}
      \mathscr{Q}_{\Phi,k}\tp{y} = \tabs{\P y}^2\tp{\Phi(y) - 1} - k.
    \end{equation}
    Then
    \begin{equation}\label{THIS_THING_WE_CARE_ABOUT}
      \mathfrak{R}\tabs{\grad_{\R^5}\mathscr{Q}_{\Phi,k}} = \tabs{{\grad_{\Pi}\mathfrak{R}\mathscr{Q}_{\Phi,k}}},
    \end{equation}
    where $\grad_{\Pi} = \tp{\pd_r,\pd_z}$.
  \end{enumerate}
\end{lem}

\begin{rmk}
  To the best of the authors' knowledge, the transformation identity~\eqref{chandrasekhar_transformation_identity} was first noticed by Chandrasekhar~\cite{Chandrasekhar1956}. Its first use in the context of the steady vortex ring problem was due to Ni~\cite{MR583638}.
\end{rmk}

The following result shows that the Chandrasekhar transformation~\eqref{chandrasekhar_transformation_identity} is the key to translating solutions of the $\R^5$-Hicks equation from the second item of Corollary~\ref{coro on solutions to the Bragg-Hawthorne equation} into traveling wave solutions of the $\R^3$-Euler equations.

\begin{prop}[Euler reconstruction]\label{prop on Euler reconstruction}
  Let $k>0$, let $\Phi\in C^\infty\tp{\R^5}$ be $\O$-invariant, and let $S,\Gamma\in C^\infty\tp{\R}$ satisfy~\eqref{_VORTEX_RING_ANSATZ_}; let $H\in C^\infty\tp{\R}$ be such that $H(0) = 0$ and $H' = \Gamma$.
  Assume that for all $y\in\R^5$ we have $\Phi(y)>0$ and
  \begin{equation}\label{another_r5_hicks_eqn}
    \tabs{\P y}^2\Delta\Phi\tp{y} + \tp{SS'}\tp{\tabs{\P y}^2\tp{\Phi(y) - 1} - k} + \tabs{\P y}^2\Gamma\tp{\tabs{\P y}^2\tp{\Phi(y) - 1} - k} = 0.
  \end{equation}
  Assume moreover that $\Phi(y)\to0$ as $|y|\to\infty$ and $\pd_5\Phi(y)<0$ whenever $y_5>0$. Define $\psi,\Psi:\Pi\to\R$ via
  \begin{equation}\label{two_types_of_stokes_stream_functions}
    \psi(r,z) = r^2\tp{\mathfrak{R}\Phi}\tp{r,z}\text{ and }\Psi(r,z) = r^2\tp{\tp{\mathfrak{R}\Phi}\tp{r,z} - 1} - k
  \end{equation}
  where $\mathfrak{R}$ is the meridional restriction defined in~\eqref{defn_meridional_restriction}. Notice in particular that $\mathfrak{R}\mathscr{Q}_{\Phi,k} = \Psi$, so this notation agrees exactly with the modified Stokes stream function introduced in~\eqref{defn_rel_stream_and_head}.

    There exist smooth and axisymmetric functions $u\in C^\infty\tp{\R^3;\R^3}$ and $p\in C^\infty\tp{\R^3}$ whose cylindrical components for $r>0$ are determined via
  \begin{equation}\label{velocity_from_stokes_stream_function}
    u = -\f{\pd_z\psi}{r}e_r + \f{\pd_r\psi}{r}e_z + \f{S\tp{\Psi}}{r}e_\theta\text{ and }p = - H\tp{\Psi} - \f12\tabs{u - 2e_3}^2
  \end{equation}
  such that the following properties are satisfied.
  \begin{enumerate}
    \item The traveling Euler equations
    \begin{equation}\label{_TRAVELING_EULER_EQUATIONS_}
      \tp{u - 2e_3}\cdot\grad u + \grad p = 0\text{ and }\grad\cdot u = 0
    \end{equation}
    are satisfied in $\R^3$.
    \item The vorticity is compactly supported, is given in cylindrical coordinates by
    \begin{equation}\label{_EXPRESSION_FOR_VORTICITY_}
      \grad\times u = \bp{r\Gamma\tp{\Psi} + \f{1}{r}\tp{SS'}\tp{\Psi}}e_\theta +  \f{S'\tp{\Psi}}{r}\tp{-\pd_z\Psi e_r + \pd_r\Psi e_z},
    \end{equation}
    and 
    \begin{equation}\label{_VORTICITY_SUPPORT_CONSIDERATIONS_}
     \supp\tp{\grad\times u} = \overline{\tcb{x\in\R^3\;:\;\Psi\tp{\tp{x_1^2+x_2^2}^{1/2},x_3}>0}}\Subset\tcb{x\in\R^3\;:\;x_1^2+x_2^2\neq0}.
    \end{equation}
    \item Upon recalling the definition of $\mathscr{Q}_{\Phi,k}$ from equation~\eqref{this_thing_we_care_about}, the helicity and total meridional vorticity are given by
    \begin{equation}\label{_HELICITY_ID_}
      \int_{\R^3}u\cdot\tp{\grad\times u} = \f{2}{\pi}\int_{\R^5}\bp{\f{1}{\tabs{\P y}^2}\tp{S\Gamma}\tp{\mathscr{Q}_{\Phi,k}\tp{y}} + \f{1}{\tabs{\P y}^4}\tp{S^2S'}\tp{\mathscr{Q}_{\Phi,k}\tp{y}}}\;\m{d}y
    \end{equation}
    and
    \begin{equation}\label{_TOTAL_MER_VORT_ID_}
      \int_{\R^3}\tabs{\tp{I - e_\theta\otimes e_\theta}\tp{\grad\times u}} = \f{1}{\pi}\int_{\R^5}\f{\tabs{S'\tp{\mathscr{Q}_{\Phi,k}\tp{y}}}}{\tabs{\P y}^3}\tabs{\grad\mathscr{Q}_{\Phi,k}\tp{y}}\;\m{d}y
    \end{equation}
    respectively.
    \item The radial component of $u$ is odd in $x_3$, while its azimuthal and axial components are even in $x_3$; moreover $e_r\cdot u\tp{x}>0$ whenever $x_1^2 + x_2^2>0$ and $x_3>0$.
  \end{enumerate}
\end{prop}
\begin{proof}
  These calculations are standard and straightforward; see, e.g., Ni~\cite{MR583638}, Amick and Fraenkel~\cite{MR816615}, and Amick and Turner~\cite{MR929976}. As such, we only give an abbreviated treatment here.

  We begin by applying the meridional restriction $\mathfrak{R}$ defined in~\eqref{defn_meridional_restriction} to the $\R^5$-Hicks equation~\eqref{another_r5_hicks_eqn}. After using the Chandrasekhar transformation of Lemma~\ref{lem on the chandrasekhar transformation} and the definitions of $\psi$ and $\Psi$ from equation~\eqref{two_types_of_stokes_stream_functions}, we find that
  \begin{equation}\label{r3_hicks_equation}
    \bp{r\pd_r\bp{\f{1}{r}\pd_r} + \pd_z^2}\psi + \tp{SS'}\tp{\Psi} + r^2\Gamma\tp{\Psi} = 0\text{ in }\Pi.
  \end{equation}

  Now let $u$ and $p$ be the smooth functions given by~\eqref{velocity_from_stokes_stream_function}. It is simple to compute that $\grad\cdot u = 0$ and hence the incompressibility constraint of the traveling Euler equations is satisfied.
  
  Working towards checking the traveling Euler momentum equation, we define the modified velocity and pressure $v\in C^\infty\tp{\R^3;\R^3}$ and $q\in C^\infty\tp{\R^3}$ via
  \begin{equation}
    v = u - 2e_3\text{ and }q = p + \f12\tabs{v}^2.
  \end{equation}
  A direct calculation reveals that
  \begin{equation}\label{q_and_v_in_terms_mod_stream_func}
    v = -\f{\pd_z\Psi}{r}e_r + \f{\pd_r\Psi}{r}e_z + \f{S\tp{\Psi}}{r}e_\theta\text{ and }q = - H(\Psi).
  \end{equation}
  From the vector identity
  \begin{equation}
    v\cdot\grad v = \f12\grad\tp{\tabs{v}^2} + \tp{\grad\times v}\times v
  \end{equation}
  we deduce that
  \begin{equation}\label{_EER_1_}
    \tp{u - 2e_3}\cdot\grad u + \grad p = \tp{\grad\times v}\times v + \grad q.
  \end{equation}

  Now, using~\eqref{r3_hicks_equation}, compute that
  \begin{multline}\label{expression_for_vorticty_sub_hicks}
    \grad\times v = \grad\times u = -\f{1}{r}\bp{r\pd_r\bp{\f{1}{r}\pd_r} + \pd_z^2}\psi e_\theta + S'\tp{\Psi}\bp{ -\f{\pd_z\Psi}{r}e_r + \f{\pd_r\Psi}{r}e_z }\\
    = \bp{\f{1}{r}\tp{SS'}\tp{\Psi} + r\Gamma\tp{\Psi}}e_\theta + S'\tp{\Psi}\bp{ -\f{\pd_z\Psi}{r}e_r + \f{\pd_r\Psi}{r}e_z }.
  \end{multline}
  Taking the cross product of the above expression with $v$ yields, after a short computation,
  \begin{equation}\label{_EER_2_}
    \tp{\grad\times v}\times v = \Gamma\tp{\Psi}\grad\Psi = \grad\tp{H(\Psi)} = -\grad q.
  \end{equation}
  The combination of~\eqref{_EER_1_} and~\eqref{_EER_2_} completes the verification that $u$ and $p$ solve the traveling Euler equations~\eqref{_TRAVELING_EULER_EQUATIONS_}. This completes the proof of the first item.

  We have already verified the first identity~\eqref{_EXPRESSION_FOR_VORTICITY_} from the second item in equation~\eqref{expression_for_vorticty_sub_hicks}. It is immediate from the support conditions~\eqref{_VORTEX_RING_ANSATZ_} that
  \begin{equation}
   \supp\tp{\grad\times u}\subseteq\overline{\tcb{x\in\R^3\;:\;\Psi\tp{\tp{x_1^2 + x_2^2}^{1/2},x_3}>0}}.
  \end{equation}
  For the opposite inclusion, we argue by contradiction. Suppose that for some open and connected set
  \begin{equation}
    \es\neq U\subseteq\tcb{\tp{r,z}\in\Pi\;:\;\Psi\tp{r,z}>0}
  \end{equation}
  we had $\grad\times u = 0$ on the solid of revolution generated by $U$. The meridional components in~\eqref{_EXPRESSION_FOR_VORTICITY_} give $S'\tp{\Psi}\grad\Psi = 0$ in $U$. At a point where $\grad\Psi\neq0$, this forces  $S'\tp{\Psi} = 0$, and the azimuthal component in~\eqref{_EXPRESSION_FOR_VORTICITY_} equals $r\Gamma\tp{\Psi}>0$, a contradiction. Hence $\grad\Psi = 0$ throughout $U$ and $\Psi$ is constant in this connected set. The azimuthal identity then reads $\tp{SS'}\tp{\Psi} + r^2\Gamma\tp{\Psi} = 0$ on an open set in which $r$ varies, which is impossible since $\Gamma\tp{\Psi}>0$. Thus the first equality in~\eqref{_VORTICITY_SUPPORT_CONSIDERATIONS_} must hold.

  To see that $\grad\times u$ must be compactly supported, we note that the hypothesis $\Phi\tp{y}\to0$ as $|y|\to\infty$ forces $\Phi\tp{y}<1$ for $y\in\R^5$ outside of some origin centered ball. On the meridional restriction of the complement of this ball the definition of $\Psi$ in~\eqref{two_types_of_stokes_stream_functions} then forces $\Psi<0$.

  To see that $\supp\grad\times u$ is supported away from the symmetry axis, we note that $\Psi(0,z) = -k<0$ for all $z\in\R$. The proof of the second item is now complete.

  Let us turn our attention to the third item; we first shall establish identity~\eqref{_TOTAL_MER_VORT_ID_}. Taking the length of the meridional projection of~\eqref{_EXPRESSION_FOR_VORTICITY_} shows
  \begin{equation}
    \tabs{\tp{I - e_\theta\otimes e_\theta}\tp{\grad\times u}} = \f{1}{r}\tabs{S'(\Psi)}\tabs{\grad\Psi}.
  \end{equation}
  Then, after integrating over $\R^3$, using cylindrical coordinates, and employing the identities~\eqref{integral_identity_key}, \eqref{this_thing_we_care_about}, and~\eqref{THIS_THING_WE_CARE_ABOUT}, we find that~\eqref{_TOTAL_MER_VORT_ID_} holds.

  Next, we prove identity~\eqref{_HELICITY_ID_}. Combining~\eqref{velocity_from_stokes_stream_function} and~\eqref{_EXPRESSION_FOR_VORTICITY_} shows
  \begin{multline}\label{_HELICITY_DENSITY_}
    u\cdot\tp{\grad\times u} = 2e_3\cdot\tp{\grad\times u} + v\cdot\tp{\grad\times v}\\\text{and }v\cdot\tp{\grad\times v} = \f{S\tp{\Psi}}{r^2}\tp{r^2\Gamma\tp{\Psi} + \tp{SS'}\tp{\Psi}} + \f{S'\tp{\Psi}}{r^2}\tabs{\grad\Psi}^2
  \end{multline}
  to be expressions for the helicity density of $u$ and $v$. Since $\grad\Psi\cdot\grad\tp{S\tp{\Psi}} = \tabs{\grad\Psi}^2S'\tp{\Psi}$, we may integrate by parts and use~\eqref{r3_hicks_equation} to equate, after a short computation,
  \begin{equation}\label{_MAGIC_IBP_}
    \int_{\Pi}\f{1}{r^2}S'\tp{\Psi}\tabs{\grad\Psi}^2r\;\m{d}r\;\m{d}z = \int_{\Pi}\tp{r^2\tp{S\Gamma}\tp{\Psi} + \tp{S^2S'}\tp{\Psi}}\f{1}{r}\;\m{d}r\;\m{d}z.
  \end{equation}
  Combining~\eqref{_HELICITY_DENSITY_} with~\eqref{_MAGIC_IBP_} and the fact that
  \begin{equation}
    \int_{\R^3}2e_3\cdot\tp{\grad\times u} = 0
  \end{equation}
  then gives
  \begin{equation}
    \int_{\R^3}u\cdot\tp{\grad\times u} = 4\pi\int_{\Pi}\bp{r\tp{S\Gamma}\tp{\Psi} + \f{1}{r}\tp{S^2S'}\tp{\Psi}}\;\m{d}r\;\m{d}z.
  \end{equation}
  Identity~\eqref{_HELICITY_ID_} now follows from the above after using~\eqref{integral_identity_key}, \eqref{THIS_THING_WE_CARE_ABOUT}, and~\eqref{this_thing_we_care_about}. The third item is now shown.

  The fourth item is trivial and so the proof is complete.
\end{proof}

\begin{rmk}[Exact prescribed invariants]\label{rmk_on_exact_prescribed_invariants}
  Let $\mathsf{p} = \tp{\ep,k,\del,\eta,\mu}\in\mathsf{P}$, where the latter set is the admissible parameter space defined in Remark~\ref{remark on instantiation of a solution map}. Let $\Phi = \Upphi_{\mathsf{p}}\in C^\infty\tp{\R^5}$ be the image of $\mathsf{p}$ under the solution map produced in that same remark. With $S$ and $\Gamma$ defined as in~\eqref{circulation_and_head_derivative_functions_}, Corollary~\ref{coro on solutions to the Bragg-Hawthorne equation} ensures that $\Phi$, $S$, and $\Gamma$ satisfy the hypotheses of Proposition~\ref{prop on Euler reconstruction}. Thus there is a corresponding traveling wave solution with velocity field $u\in C^\infty\tp{\R^3;\R^3}$. 

  By combining the third item of Proposition~\ref{prop on Euler reconstruction} with the third item of Corollary~\ref{coro on solutions to the Bragg-Hawthorne equation}, we see that
  \begin{equation}
    \del\eta = \int_{\R^3}u\cdot\tp{\grad\times u}\text{ and }\del\mu = \int_{\R^3}\tabs{\tp{I - e_\theta\otimes e_\theta}\tp{\grad\times u}}.
  \end{equation}
  Thus, it is in the above sense that the chosen parameters $\del$, $\eta$, and $\mu$ prescribe the helicity and total meridional vorticity of the Euler velocity field $u$.
\end{rmk}

\subsection{Vortex geometry and velocity estimates}\label{SS_vortex_geom_velo_ests}

We begin by discussing the geometry of the vortex core for traveling vortex rings near Hill's spherical vortex.

\begin{lem}[On the vortex core]\label{lem_on_on_the_vortex_core}
  There exist constants $0<\ep_\star,k_\star<1/64$ and $\kappa>1$ with the following property. Suppose $\Phi$, $S$, and $\Gamma$ satisfy the hypotheses of Proposition~\ref{prop on Euler reconstruction}, let $u$ be the resulting Euler velocity field as defined in~\eqref{velocity_from_stokes_stream_function}, and assume additionally that
  \begin{equation}\label{_HYPOTHERMIA_}
    0<k\le k_\star,\;\Phi<1\text{ in }\R^5\setminus B_{\R^5}\tp{0,3},\;\text{ and }\tnorm{\Phi - \v}_{C^1_{\m{b}}\tp{B_{\R^5}[0,3]}}\le\ep_\star,
  \end{equation}
  where $\v$ is Hill's solution as defined in~\eqref{Hills_Vortex}. Set
  \begin{equation}\label{the_positive_set_charged_up}
    \mathfrak{P} = \tcb{\tp{r,z}\in\Pi\;:\;\Psi(r,z)>0}\text{ where, as in~\eqref{two_types_of_stokes_stream_functions}, }\Psi\tp{r,z} = r^2\tp{\tp{\mathfrak{R}\Phi}\tp{r,z} - 1} - k.
  \end{equation}
  The following hold.
  \begin{enumerate}
    \item With
    \begin{equation}\label{_A_SET_FOR_YOU_}
      \mathfrak{K} = \tcb{\tp{r,z}\in\Pi\;:\;\kappa^{-1}\le r\text{ and }r^2 + z^2\le\tp{1 - \kappa^{-1}}^2}
    \end{equation}
    we have
    \begin{equation}\label{assertions_here_for_today}
      \tcb{\tp{r,z}\in\Pi\;:\;\Psi(r,z)\in[1/32,5/16]}\subseteq\mathfrak{K},\;[1/32,5/16]\subseteq\Psi(\Pi),\text{ and }\sup_{\mathfrak{K}}\tabs{\grad\Psi}\le\kappa,
    \end{equation}
    where $\grad = \grad_{\Pi} = \tp{\pd_r,\pd_z}$.
    \item The set $\mathfrak{P}$ is a bounded topological disk compactly contained in $\Pi$ and $\pd\mathfrak{P}$ is a smooth simple closed curve. Consequently, $\supp\tp{\grad\times u}$ is a topological solid torus with smooth boundary and
    \begin{equation}\label{_PUMP_FOIL_}
      \supp\tp{\grad\times u}\subseteq\tcb{x\in\R^3\;:\;|x|<2\text{ and }2\tp{x_1^2+x_2^2}>k}.
    \end{equation}
    \item If $\supp S\subset[1/32,5/16]$, then
    \begin{equation}\label{_swirl_inclusion_here_}
      \supp\tp{e_\theta\cdot u}\subseteq\tcb{x\in\R^3\;:\;\kappa^{-2}\le x_1^2 + x_2^2\text{ and }|x|\le 1 - \kappa^{-1}}.
    \end{equation}
    Moreover, $S\not\equiv0$ implies $e_\theta\cdot u\not\equiv 0$.
  \end{enumerate}
\end{lem}
\begin{proof}
  Consider first the modified Stokes stream function for Hill's solution,
  \begin{equation}
    \Psi_{\m{H}}\tp{r,z} = r^2\tp{\tp{\mathfrak{R}\v}\tp{r,z} - 1}.
  \end{equation}
  One verifies that
  \begin{equation}
    \Psi_{\m{H}}\tp{0,z} = 0\text{ for all }z\in\R\text{ and }\Psi_{\m{H}}\tp{r,z} = 0\text{ whenever }r^2+z^2 = 1.
  \end{equation}
  Moreover, the inequalities
  \begin{equation}
    \Psi_{\m{H}}\tp{r,z}>0\text{ when }r>0,\;r^2+z^2<1\text{ and }\Psi_{\m{H}}\tp{r,z}<0\text{ when }r>0,\;r^2+z^2>1
  \end{equation}
  hold. From these observations paired with the fact that $\max\Psi_{\m{H}} = 3/8$, we deduce the existence of some $\kappa>1$, that determines the set~\eqref{_A_SET_FOR_YOU_}, for which 
  \begin{equation}
    \tcb{\tp{r,z}\in\Pi\;:\;1/32\le\Psi_{\m{H}}\tp{r,z}\le11/32}\Subset\mathfrak{K}.
  \end{equation}

  Now for $\tp{r,z}\in\Pi$ with $r^2+z^2\le 4$ we have
  \begin{equation}
    \Psi - \Psi_{\m{H}} = r^2\mathfrak{R}\tp{\Phi - \v} - k\text{ and hence }\tnorm{\Psi - \Psi_{\m{H}}}_{C^1_{\m{b}}\tp{\Pi\cap B_{\R^2}[0,2]}}\lesssim \ep_\star + k_\star.
  \end{equation}
  Thus, by taking $\ep_\star$ and $k_\star$ small enough and $\kappa$ larger, if required, the first and third assertions of~\eqref{assertions_here_for_today} follow. The intermediate assertion follows by a similar argument, but one also employs the intermediate value theorem.

  Next, we study the vortex core. Let $g_{\m{H}},g:\R^+\to\R$ be defined for $r\in\R^+$ via
  \begin{equation}
    g_{\m{H}}\tp{r} = \Psi_{\m{H}}\tp{r,0} = \begin{cases}
      \f32r^2\tp{1 - r^2}&\text{for }0<r\le 1,\\
      r^2\tp{\f{1}{r^3} - 1}&\text{for }r>1,
    \end{cases}
    \text{ and }g(r) = \Psi(r,0).
  \end{equation}
  Then, smallness of $\ep_\star$ and $k_\star$ guarantees $C^1_{\m{b}}\tp{[0,3/2]}$-closeness of $g$ and $g_{\m{H}}$. Thus, by decreasing these parameters, if necessary, we can enforce that
  \begin{equation}\label{_facts_about_function_g_}
    g'>0\text{ on }(0,1/2],\;g>0\text{ on }[1/2,3^{1/2}/2],\text{ and }g'<0\text{ on }[3^{1/2}/2,3/2].
  \end{equation}
  In fact, the same smallness also gives $\Phi(y)<1$ whenever $y\in\R^5$ satisfies $3/2\le|y|\le3$ while this inequality is assumed for $|y|\ge3$. Consequently, $\Psi(r,z)<0$ whenever $r^2+z^2\ge\tp{3/2}^2$.

  Since $g(0) = -k<0$, there are unique numbers 
  \begin{equation}
    0<r_-<\f12<\f{3^{1/2}}{2}<r_+<\f32
  \end{equation}
  such that $g(r)>0$ for $r\in\R^+$ if and only if $r_-<r<r_+$.

  Then, for every fixed $r>0$, the function $\R\ni z\mapsto \Psi(r,z)\in\R$ is even and strictly decreasing on $\R^+$, since $\pd_z\Psi(r,z) = r^2\mathfrak{R}\tp{\pd_5\Phi}\tp{r,z}<0$ for $z>0$. We also know that $\lim_{z\to+\infty}\Psi(r,z) = -r^2 - k$. From these facts combined with the implicit function theorem, we deduce the existence of a smooth function $Z:\tp{r_-,r_+}\to\R^+$ such that
  \begin{equation}
    \mathfrak{P} = \tcb{\tp{r,z}\in\Pi\;:\;r_-<r<r_+\text{ and }|z|<Z(r)},
  \end{equation}
  where the former set is defined in~\eqref{the_positive_set_charged_up}. Continuity and strict decrease in $z$ force $Z(r)\to0$ as $r\to r_{\pm}$. The non-vanishing derivatives in~\eqref{_facts_about_function_g_} then allow us to solve for $r$ near the two endpoints $r_-$ and $r_+$, so, in fact, the set $\pd\mathfrak{P}$ is a smooth simple closed curve contained within $\Pi$.

  On $\mathfrak{P}$ we have $\mathfrak{R}\Phi - 1\le3/2 + \ep_\star<2$ and hence $\Psi\tp{r,z}\ge0$ on $\tp{r,z}\in\Bar{\mathfrak{P}}$ implies that $2r^2>k$. The vorticity support equality~\eqref{_VORTICITY_SUPPORT_CONSIDERATIONS_} from the second item of Proposition~\ref{prop on Euler reconstruction} now gives the inclusion in~\eqref{_PUMP_FOIL_}, and the solid of revolution of $\Bar{\mathfrak{P}}$ is a topological solid torus with smooth boundary. This gives the second item.

  For the third item, the inclusion~\eqref{_swirl_inclusion_here_} is immediate from the first inclusion of~\eqref{assertions_here_for_today} and the expression for the azimuthal velocity given in~\eqref{velocity_from_stokes_stream_function}. The second inclusion of~\eqref{assertions_here_for_today} readily proves the subsequent stated implication.
\end{proof}

We continue by discussing estimates relevant for the meridional components of the velocity field. Given any $\Phi\in C^1\tp{\R^5}$ that is $\O$-invariant, in the sense that $\Phi(Ty) = \Phi(y)$ for all $y\in\R^5$ and $T\in\O$ (where $\O$ is from the first item of Definition~\ref{defn on spaces and mappings}), we define $\mathcal{U}[\Phi]\in C^0\tp{\R^3;\R^3}$ via
\begin{equation}\label{_MER_VELO_MAP_}
  \mathcal{U}[\Phi] = -\f{\pd_z\tp{r^2\mathfrak{R}\Phi}}{r}e_r + \f{\pd_r\tp{r^2\mathfrak{R}\Phi}}{r}e_z = -r\pd_z\mathfrak{R}\Phi e_r + \tp{2\mathfrak{R}\Phi + r\pd_r\mathfrak{R}\Phi}e_z,
\end{equation}
where $\mathfrak{R}$ is the meridional restriction defined in~\eqref{defn_meridional_restriction}. There are then two key mapping properties of $\mathcal{U}$ which we shall use. The following Banach spaces will play important roles. First,
\begin{equation}\label{yet_another_space}
  \mathsf{H}\tp{\R^5} = \tcb{\Phi\in C^1\tp{\R^5}\;:\;\forall\;T\in\O,\;\forall\;x\in\R^5,\;\Phi\tp{Tx} = \Phi(x),\;\tnorm{\Phi}_{\mathtt{n}}<\infty}
\end{equation}
is defined with the norm
\begin{equation}\label{yet_another_norm}
  \tnorm{\Phi}_{\mathtt{n}} = \sup\tcb{\tbr{x}^3\tabs{\Phi(x)} + \tbr{x}^4\tabs{\grad\Phi\tp{x}}\;:\;x\in\R^5}.
\end{equation}
We also introduce
\begin{equation}
  \tilde{\mathsf{H}}\tp{\R^5} = \tcb{\Phi\in C^1\tp{\R^5}\;:\;\forall\;T\in\O,\;\forall\;y\in\R^5,\;\Phi(Ty) = \Phi(y),\;\tnorm{\Phi}_{\tilde{\mathtt{n}}}<\infty}
\end{equation}
for the norm
\begin{equation}
  \tnorm{\Phi}_{\tilde{\mathtt{n}}} = \tnorm{\tbr{\cdot}^3\Phi}_{\m{LL}^1\tp{\R^5}}.
\end{equation}

The alluded to mapping properties are that $\mathcal{U}$ restricts to bounded linear operators
\begin{equation}\label{mer_map_props}
  \mathcal{U}:\mathsf{H}\tp{\R^5}\to\tp{L^2\cap C^0_{\m{b}}}\tp{\R^3}\text{ and }\mathcal{U}:\tilde{\mathsf{H}}\tp{\R^5}\to\m{LL}\tp{\R^3}.
\end{equation}
The verification of the boundedness asserted in~\eqref{mer_map_props} is elementary and so we omit the proof here. 

The particular consequence of the above that arises in our subsequent application is the following. Suppose that $\Phi_\infty$ and $\tcb{\Phi_n}_{n\in\N}$ are elements of $\tp{\tilde{\mathsf{H}}\cap\mathsf{H}}\tp{\R^5}$ such that
\begin{equation}\label{hypothesis_bounds_and_limits}
  \sup_{n\in\N}\tnorm{\Phi_n}_{\tilde{\mathsf{H}}\tp{\R^5}}<\infty\text{ and }\lim_{n\to\infty}\tnorm{\Phi_n - \Phi_\infty}_{\mathsf{H}\tp{\R^5}} = 0.
\end{equation}
Then, by a simple interpolation argument, for every $0\le\al<1$ it is the case that
\begin{equation}\label{the_conclusion_limit}
  \lim_{n\to\infty}\tnorm{\mathcal{U}[\Phi_n] - \mathcal{U}[\Phi_\infty]}_{\tp{L^2\cap C^\al_{\m{b}}}\tp{\R^3}} = 0.
\end{equation}

\begin{lem}[Estimates on the azimuthal velocity]
\label{lem_on_azimuthal_velocity_estimates}
  Let $\ep_\star$, $k_\star$, and $\kappa$ be as in Lemma~\ref{lem_on_on_the_vortex_core}. There exists $C\in\R^+$ with the following property. Suppose that $\Phi$, $S$, and $\Gamma$ satisfy the hypotheses of both Proposition~\ref{prop on Euler reconstruction} and Lemma~\ref{lem_on_on_the_vortex_core}, and assume additionally that $\supp S\subset[1/32,5/16]$. Let $\Psi$ be the modified Stokes stream function from equation~\eqref{two_types_of_stokes_stream_functions} and let $u$ be the corresponding Euler velocity field as in~\eqref{velocity_from_stokes_stream_function}. Then the following hold.
  \begin{equation}\label{SWIRL_SUPPORT}
    \supp\tp{e_\theta\cdot u}\subseteq\tcb{x\in\R^3\;:\;\kappa^{-2}\le x_1^2+x_2^2,\;|x|\le 1-\kappa^{-1}}.
  \end{equation}
  Moreover, for any $0\le\al<1$, we have the estimates
  \begin{equation}\label{SWIRL_ESTIMATES}
    \tnorm{\tp{e_\theta\otimes e_\theta}u}_{L^2\tp{\R^3}}+\tnorm{\tp{e_\theta\otimes e_\theta}u}_{\m{LL}\tp{\R^3}}\le C\tnorm{S}_{\m{LL}\tp{\R}}\text{ and }\tnorm{\tp{e_\theta\otimes e_\theta}u}_{C^\al_{\m{b}}\tp{\R^3}}\le C\tnorm{S}_{C^\al_{\m{b}}\tp{\R}}.
  \end{equation}
\end{lem}
\begin{proof}
  The support inclusion~\eqref{SWIRL_SUPPORT} is precisely the third item of Lemma~\ref{lem_on_on_the_vortex_core}. In particular,
  $\supp\tp{S\tp{\Psi}} = \supp\tp{e_\theta\cdot u}$
  is contained in a fixed compact subset of $\R^3$ that is separated from the symmetry axis. By the $C^1$-closeness hypothesis~\eqref{_HYPOTHERMIA_} on $\Phi$ and the definition of $\Psi$, there exists a fixed bounded neighborhood $\mathcal{O}$, a positive distance from the symmetry axis, of this compact set on which
  \begin{equation}
    \tnorm{\Psi}_{C^{0,1}_{\m b}(\mathcal{O})}\le C.
  \end{equation}
  Choose $\upchi\in C^\infty_{\m c}(\mathcal{O})$ equal to one on the set on the right hand side of~\eqref{SWIRL_SUPPORT}. Then
  \begin{equation}
    \tp{e_\theta\otimes e_\theta}u=\upchi\f{e_\theta}{r}S(\Psi)\text{ and }
    \tnorm{\upchi r^{-1}e_\theta}_{C^1_{\m b}(\R^3)}\le C.
  \end{equation}
  Composition and product estimates therefore give, for $0\le\al<1$,
  \begin{multline}
    \tnorm{\tp{e_\theta\otimes e_\theta}u}_{\m{LL}(\R^3)}
    \le C\tp{1+\tnorm{\Psi}_{C^{0,1}_{\m b}(\mathcal{O})}}
    \tnorm{S}_{\m{LL}(\R)}
    \le C\tnorm{S}_{\m{LL}(\R)},\\
    \tnorm{\tp{e_\theta\otimes e_\theta}u}_{C^\al_{\m b}(\R^3)}
    \le C\tp{1+\tnorm{\Psi}_{C^{0,1}_{\m b}(\mathcal{O})}}
    \tnorm{S}_{C^\al_{\m b}(\R)}
    \le C\tnorm{S}_{C^\al_{\m b}(\R)}.
  \end{multline}
    The remaining estimate follows from~\eqref{SWIRL_SUPPORT}, since the set therein is bounded and $r^{-1}\le\kappa$ on it:
  \begin{equation}
    \tnorm{\tp{e_\theta\otimes e_\theta}u}_{L^2(\R^3)}
    \le C\tnorm{S}_{C^0_{\m b}(\R)}.
  \end{equation}
\end{proof}

\section{Desingularization of the Hill--Norbury family and main results}\label{SO_desing_of_HNF_and_main_results}

We now identify the singular limits of the smooth solutions constructed above and prove the main results. The key additional ingredient is the local uniqueness theory for the idealized Hill--Norbury rings, which allows us to identify the limits obtained by compactness before applying the reconstruction estimates of Section~\ref{SO_euler_recon_and_vortex_geometry}.

\subsection{On the uniqueness of the Hill--Norbury idealized rings}\label{SS_uniqueness_of_Norbury}

It is a nontrivial fact that Norbury's family of vortex rings is known to be conditionally unique. As this uniqueness theorem is important for our forthcoming analysis, we endeavor now to precisely describe this uniqueness result and a more conveniently formulated corollary.

The reader is encouraged to recall the weak formulation of the idealized vortex ring problem discussed in Section~\ref{SS_sol_HN_and_nondim} and Norbury's small-flux solution family $\tcb{\psi^{\m{N}}_k}_{0<k<k_{\m{N}}}\subset H(\Pi)$. The question of uniqueness for Norbury's family is far more subtle than that of Hill's vortex, the latter of which was shown to be unconditionally unique~\cite{MR816615}. Norbury's original construction only gives a very localized uniqueness result. Later, Amick and Fraenkel~\cite{MR918795} expanded the sense in which Norbury's solutions are unique, although it is not a completely unconditional result. Nevertheless, it is this latter result that is sufficient for our applications.

In~\cite{MR918795} the authors prove the following: There exists $\tilde{M}_0\in\R^+$ such that for all $\tilde{M}_1\ge\tilde{M}_0$, there exists $0<\tilde{k}_1\le 1$ with the property that for all $0<k\le\tilde{k}_1$ and all $\phi_k\in H\tp{\Pi}\setminus\tcb{0}$ satisfying
\begin{equation}
  \tnorm{\phi_k}_{\mathtt{o}}\le\tilde{M}_1\text{ and }\forall\;\tp{r,z}\in\Pi,\;\phi_k\tp{r,z} = \phi_k\tp{r,-z},
\end{equation}
if $\phi_k$ is a solution to the problem~\eqref{_NON_DIM_NORBURY_WEAK_PROBLEM_} with flux $k$, then $\phi_k = \psi^{\m{N}}_k$. In other words, if $\phi_k$ is a $z$-even solution to~\eqref{_NON_DIM_NORBURY_WEAK_PROBLEM_}, that is not too large in the norm $\tnorm{\cdot}_{\mathtt{o}}$, and $k$ is sufficiently small, then $\phi_k$ must be Norbury's original solution $\psi^{\m{N}}_k$.

We note that in the paper~\cite{MR918795}, see the discussion in Section 1.4, there are two completely equivalent formulations of the problem~\eqref{_NON_DIM_NORBURY_WEAK_PROBLEM_}: one for axisymmetric functions on $\R^3$ and another for cylindrically symmetric functions on $\R^5$. This is similar to what we have explored in Proposition~\ref{prop on Euler reconstruction} above. It is this latter $\R^5$-perspective that is most convenient for our application. We shall state a corollary of the main uniqueness result of Amick and Fraenkel~\cite{MR918795}, in which the Banach space $\mathsf{H}\tp{\R^5}$, defined in~\eqref{yet_another_space}, will play a role.

\begin{lem}\label{lem on Norbury uniqueness}
  Given $\R\ni k>0$ we define $\mathcal{N}\tp{k}\subset\mathsf{H}\tp{\R^5}$ to be the collection of $\Phi\in\mathsf{H}\tp{\R^5}\setminus\tcb{0}$ such that for all $\psi\in C^\infty_{\m{c}}\tp{\R^5}$ we have 
  \begin{equation}\label{R5_form_norbury}
    \int_{\R^5}\grad\Phi\tp{y}\cdot\grad\psi\tp{y}\;\m{d}\mathcal{H}^5\tp{y} = 15\int_{\R^5}\mathds{1}_{\mathsf{A}\tp{\Phi,k}}\tp{y}\psi\tp{y}\;\m{d}\mathcal{H}^5\tp{y}
  \end{equation}
  where
  \begin{equation}\label{R5_form_vortex_core}
    \mathsf{A}\tp{\Phi,k} = \tcb{y\in\R^5\;:\;\tabs{\P y}^2\tp{\Phi(y)-1}>k}.
  \end{equation}

  There exists $M_0\in\R^+$ such that for all $\R\ni M_1\ge M_0$, there exists $0<k_1\le 1$ with the property that for all $0<k\le k_1$ we have 
  \begin{equation}\label{uniqueness assertion}
    \m{card}\tcb{\Phi\in\mathcal{N}\tp{k}\;:\;\tnorm{\Phi}_{\mathtt{n}}\le M_1} \le 1,
  \end{equation}
  where $\m{card}(\cdot)$ denotes the cardinality of a set. Moreover, every member of the set in~\eqref{uniqueness assertion} corresponds under $\mathsf{S}$ to Norbury's solution $\psi_k^{\m{N}}$ from Definition~\ref{defn on the Hill Norbury family}.
\end{lem}
\begin{proof}
  As discussed in Section 1.4 in~\cite{MR918795} (see also~\cite{MR583638,MR929976} and Section~\ref{SO_euler_recon_and_vortex_geometry}), there is a transformation between the $\R^3$ and $\R^5$ perspectives of Norbury's and related steady vortex ring problems. Consider the injective linear transformation
  \begin{equation}\label{the_key_transformation}
    \mathsf{S}:\mathsf{H}\tp{\R^5}\to H\tp{\Pi}\cap\tcb{\phi\in H\tp{\Pi}\;:\;\phi\tp{r,z} = \phi\tp{r,-z},\;\forall\;\tp{r,z}\in\Pi}
  \end{equation}
  that is defined via
  \begin{equation}
    \tp{\mathsf{S}\Phi}\tp{r,z} = r^2\tp{\mathfrak{R}\Phi}\tp{r,z},\;\forall\;\Phi\in\mathsf{H}\tp{\R^5},\;\forall\;\tp{r,z}\in\Pi,
  \end{equation}
  where $\mathfrak{R}$ is the meridional restriction introduced in equation~\eqref{defn_meridional_restriction}. The continuity estimate $\tnorm{\mathsf{S}\Phi}_{\mathtt{o}}\lesssim\tnorm{\Phi}_{\mathtt{n}}$ is straightforward to verify given the weighted control in the norm~\eqref{yet_another_norm}.

  Writing $\phi=\mathsf{S}\Phi$, a direct computation following the Chandrasekhar transformation~\eqref{chandrasekhar_transformation_identity} gives
  \begin{equation}
    \bp{r\pd_r\bp{\f1r\pd_r}+\pd_z^2}\phi\tp{r,z}=r^2\mathfrak{R}\tp{\Delta\Phi}\tp{r,z}
  \end{equation}
  in the sense of distributions. Consequently, equation~\eqref{R5_form_norbury} is equivalent by density to the weak identity~\eqref{_NON_DIM_NORBURY_WEAK_PROBLEM_}, since the condition $\tabs{\P y}^2\tp{\Phi(y)-1}>k$ becomes $\phi(r,z)>r^2+k$.
  Thus, for all $k\in\R^+$ and all $\Phi\in\mathcal{N}\tp{k}$ the function $\phi = \mathsf{S}\Phi\in H(\Pi)\setminus\tcb{0}$ is a solution to the problem~\eqref{_NON_DIM_NORBURY_WEAK_PROBLEM_}. After adjusting $M_0$ to absorb the constant in the estimate $\tnorm{\mathsf{S}\Phi}_{\mathtt{o}}\lesssim\tnorm{\Phi}_{\mathtt{n}}$, the results of the lemma therefore follow immediately from the conditional uniqueness result of Amick and Fraenkel~\cite{MR918795} and the injectivity of $\mathsf{S}$.
\end{proof}

\subsection{Flexible desingularization of the Hill--Norbury family}\label{SS_desingularization_HN_fam}

The task of this subsection is to study particular vanishing limits in $\ep>0$ of the solutions $\Phi$ constructed by Corollary~\ref{coro on solutions to the Bragg-Hawthorne equation}. The Hill--Norbury fields $u^k$ fixed in Definition~\ref{defn on the Hill Norbury family} are the intended target.

Lemma~\ref{lem_on_on_the_vortex_core} has introduced new small positive parameters $\ep_\star$ and $k_\star$ such that within their margins we have quite specific information on the vortex core; thus we shall need to restrict our solution map from Remark~\ref{remark on instantiation of a solution map} to ensure that members of its image meet these desired thresholds.


\begin{rmk}[Restriction of the solution map]\label{remark_on_restriction_of_solution_map}
  Recall the small positive parameter $\del_0$ introduced in Theorem~\ref{thm on existence of a fixed point}. As a simple consequence of Corollary~\ref{coro on solutions to the Bragg-Hawthorne equation}, there exists $0<\del_\star\le\min\tcb{\del_0,k_\star}$ such that whenever $0\le\ep,k,\del\le\del_\star$ with $\ep,k>0$, $\eta\in[-\del_\star,\del_\star]$, and $\mu\in[1/\del_\star,\infty)$, the image of $\mathsf{p} = \tp{\ep,k,\del,\eta,\mu}\in\mathsf{P}$, $\Phi = \Upphi_{\mathsf{p}}\in C^\infty\tp{\R^5}$, satisfies
  \begin{equation}
    \tnorm{\Phi - \v}_{C^1_{\m{b}}\tp{B_{\R^5}[0,3]}}\le\ep_\star\text{ and }\Phi<1\text{ in }\R^5\setminus B_{\R^5}\tp{0,3}.
  \end{equation}

  Thus, with the smooth functions $S$ and $\Gamma$ defined in terms of $\mathsf{p}$ as in~\eqref{circulation_and_head_derivative_functions_}, we have that $\Phi$, $S$, $\Gamma$, and $k$ satisfy the hypotheses of both Proposition~\ref{prop on Euler reconstruction} and Lemma~\ref{lem_on_on_the_vortex_core}.

  This justifies introducing the following notation for the restriction of the parameter space~\eqref{defn_admissible_parameter_space}:
  \begin{equation}\label{restricted_parameter_space}
    \mathsf{P}_\star = (0,\del_\star]^2\times[0,\del_\star]\times[-\del_\star,\del_\star]\times[1/\del_\star,\infty).
  \end{equation}
\end{rmk}

In addition to the solution map $\mathsf{p}\mapsto\Upphi_{\mathsf{p}}$ mentioned above, and its restriction to~\eqref{restricted_parameter_space}, our next result also requires the introduction of the following sets of admissible sequences. Let
\begin{multline}\label{set_of_admissible_sequences}
  \mathscr{A}_0 = \tcb{\tp{\tcb{\upeta_n}_{n\in\N},\tcb{\upmu_n}_{n\in\N}}\;:\;\upeta_n = \upmu_n = 0\text{ for all }n\in\N},\\
  \mathscr{A}_1 = \tcb{\tp{\tcb{\upeta_n}_{n\in\N},\tcb{\upmu_n}_{n\in\N}}\;:\;\lim_{n\to\infty}\upeta_n = 0\text{ and }\forall\;n\in\N,\;\upmu_n\ge1,\;|\upeta_n|\le\del_\star^2},
\end{multline}
and $\mathscr{A} = \mathscr{A}_0\cup\mathscr{A}_1$.

\begin{prop}[Flexible desingularization of the Hill--Norbury family]\label{label prop on flexible desingularization of the hill norbury family}
  There exist constants $C\in\R^+$ and $0<k_0\le\del_\star$ (with $\del_\star>0$ from Remark~\ref{remark_on_restriction_of_solution_map}) such that for all $0\le k\le k_0$ and all $\tp{\tcb{\upeta_n}_{n\in\N},\tcb{\upmu_n}_{n\in\N}}\in\mathscr{A}$ there exists a sequence $\tcb{\Phi_n}_{n\in\N}\subset C^\infty\tp{\R^5}$ such that the following hold.
  \begin{enumerate}
    \item The sequence $\tcb{\Phi_n}_{n\in\N}\subset C^\infty\tp{\R^5}$ is in the image of the restricted solution map of Remarks~\ref{remark on instantiation of a solution map} and~\ref{remark_on_restriction_of_solution_map}. More precisely, for each $n\in\N$ there exists $\mathsf{p}_n = \tp{\ep_n,k_n,\del_n,\eta_n,\mu_n}\in\mathsf{P}_\star$ such that $\Phi_n = \Upphi_{\mathsf{p}_n}$.
    \item The parameters introduced in the previous item satisfy
    \begin{equation}\label{_the_parameter_relationships_}
      k_n = \begin{cases}
        k&\text{if }k>0,\\
        2^{-n}\del_\star&\text{if }k=0,
      \end{cases}\;\ep_n = 2^{-n}\del_\star,\;\del_n\eta_n = \upeta_n,\;\del_n\mu_n = \upmu_n
    \end{equation}
    for all $n\in\N$. In particular, $k_n>0$ for all $n\in\N$, $k_n\to k$, $\ep_n\to0$, and $\del_n\to0$ as $n\to\infty$. Additionally, if $\tp{\tcb{\upeta_n}_{n\in\N},\tcb{\upmu_n}_{n\in\N}}\in\mathscr{A}_0$, then $\del_n = 0$ for all $n\in\N$.
    \item We have the uniform bounds
    \begin{equation}\label{key_uniform_bounds}
      \sup_{n\in\N}\tsb{\tnorm{\tbr{\cdot}^3\Phi_n}_{\m{LL}^1\tp{\R^5}} + \tnorm{\Phi_n}_{H^2\tp{\R^5}}}\le C
    \end{equation}
    and the Cauchy condition 
    \begin{equation}\label{_key_cauchy_condition_}
      \limsup_{m,n\to\infty}\ssb{\tnorm{\tbr{\cdot}^3\tp{\Phi_n - \Phi_m}}_{L^\infty\tp{\R^5}} + \tnorm{\tbr{\cdot}^4\tp{\grad\Phi_n - \grad\Phi_m}}_{L^\infty\tp{\R^5}} } = 0
    \end{equation}
    holds.
    \item There exists $\Phi_\infty\in\mathsf{H}\tp{\R^5}\cap\tilde{\mathsf{H}}\tp{\R^5}\cap H^2\tp{\R^5}$ (see~\eqref{yet_another_space}) such that $\Phi_n\to\Phi_\infty$ as $n\to\infty$ in the $\mathsf{H}\tp{\R^5}$ metric. If $k = 0$, then $\Phi_\infty = \v$, where $\v$ is Hill's solution from~\eqref{Hills_Vortex}. If $k>0$, then $\Phi_\infty\in\mathcal{N}(k)$, where the latter set was introduced in Lemma~\ref{lem on Norbury uniqueness}. 
  \end{enumerate}
\end{prop}
\begin{proof}
  Let us fix sequences and parameters as in the statement of the proposition and for now let $0\le k\le\del_\star$. The parameter $k_0>0$ will be selected at a later point in the proof.
  
  We begin by defining the sequence of functions $\tcb{\Phi_n}_{n\in\N}\subset C^\infty\tp{\R^5}$. To enforce the first item, we take $\Phi_n = \Upphi_{\mathsf{p}_n}$ where the definition of the points $\tcb{\mathsf{p}_n}_{n\in\N}\subset\mathsf{P}_\star$ is described next. We take $\mathsf{p}_n = \tp{\ep_n,k_n,\del_n,\eta_n,\mu_n}$. The numbers $k_n$ and $\ep_n$ are defined as in equation~\eqref{_the_parameter_relationships_}.
  
  We shall define $\del_n$, $\eta_n$, and $\mu_n$ so as to satisfy $\del_n\eta_n = \upeta_n$ and $\del_n\mu_n = \upmu_n$. Suppose first that
  \begin{equation*}
    \tp{\tcb{\upeta_n}_{n\in\N},\tcb{\upmu_n}_{n\in\N}}\in\mathscr{A}_1.
  \end{equation*}
  We then take \begin{equation}\label{corrected_delta_selection}
    \del_n = \max\bcb{|\upeta_n|^{1/2},\f{\del_\star}{2^n}},\;
    \eta_n = \f{\upeta_n}{\del_n},\text{ and }\mu_n = \f{\upmu_n}{\del_n}.
  \end{equation}
  Since $|\upeta_n|\le\del_\star^2$ and $\upmu_n\ge1$, we have $0<\del_n\le\del_\star$, $|\eta_n|\le\del_\star$, and $\mu_n\ge1/\del_\star$. Moreover, $\del_n\to0$ since $\upeta_n\to0$ as $n\to\infty$. On the other hand, if $\tp{\tcb{\upeta_n}_{n\in\N},\tcb{\upmu_n}_{n\in\N}}\in\mathscr{A}_0$ we simply take $\del_n = 0$, $\eta_n = 0$, and $\mu_n = 1/\del_\star$ for all $n\in\N$. The first and second items are now satisfied by definition. 

  Let us pause and record the equation~\eqref{bragg_hawthorne_r5} satisfied by the sequence $\tcb{\Phi_n}_{n\in\N}$ in a useful weak formulation. Fix $n\in\N$. We let $S_n,\Gamma_n\in C^\infty\tp{\R}$ denote the circulation and head derivative functions associated via the first item of Corollary~\ref{coro on solutions to the Bragg-Hawthorne equation} to $\Phi_n$. In particular, $\Gamma_n\tp{\cdot} = 15\cdot\Gamma_{\ep_n}\tp{\cdot}$ and $|S_n(\cdot)\cdot S_n'\tp{\cdot}|\le C\del_n^2$. For any $\psi\in C^\infty_{\m{c}}\tp{\R^5}$, we can test against~\eqref{bragg_hawthorne_r5} divided by $\tabs{\P y}^2$ and integrate by parts to derive the identity:
  \begin{multline}\label{the_key_weak_identity}
    \int_{\R^5}\grad\Phi_n\tp{y}\cdot\grad\psi\tp{y}\;\m{d}\mathcal{H}^5\tp{y} = \int_{\R^5}\f{1}{\tabs{\P y}^2}\tp{S_n\cdot S_n'}\tp{\tabs{\P y}^2\tp{\Phi_n\tp{y} - 1} - k_n}\psi\tp{y}\;\m{d}\mathcal{H}^5\tp{y}\\+\int_{\R^5}\Gamma_n\tp{\tabs{\P y}^2\tp{\Phi_n\tp{y} - 1} - k_n}\psi\tp{y}\;\m{d}\mathcal{H}^5\tp{y}.
  \end{multline}
  The main technical task of the remainder of the proof is to identify the limit as $n\to\infty$ of the final term in~\eqref{the_key_weak_identity}. This is subtle since $\Gamma_n$ is becoming a discontinuous function in the limit.

  Now let us prove the third item. Note that the uniform bounds~\eqref{key_uniform_bounds} are satisfied immediately thanks to the fourth item of Corollary~\ref{coro on solutions to the Bragg-Hawthorne equation}. If $k = 0$, then $\ep_n+k_n\to0$ and $\del_n\to0$; consequently, the final estimate in~\eqref{hill_deviation_estimates} gives $\Phi_n\to\v$ in the $\mathsf{H}\tp{\R^5}$ metric. In this case the Cauchy condition~\eqref{_key_cauchy_condition_} and the fourth item are immediate as well.

  Suppose now that $k>0$. Then $k_n=k$ for all $n\in\N$. Let us observe the following reduction. From the proof of Corollary~\ref{coro on solutions to the Bragg-Hawthorne equation}, we recall that $\Delta\Phi_n = 0$ in $\R^5\setminus B_{\R^5}[0,2]$ for each $n\in\N$. This fact, combined with the $\tbr{\cdot}^{-3}$ decay allows us to argue similarly to the derivation of estimates~\eqref{_initial_bound_of_3_decay_} and~\eqref{final_bound_of_4_decay} and show that the estimate
  \begin{equation}\label{_reduction_estimate_}
    \tnorm{\tbr{\cdot}^3\tp{\Phi_n - \Phi_m}}_{L^\infty\tp{\R^5}} + \tnorm{\tbr{\cdot}^4\tp{\grad\Phi_n - \grad\Phi_m}}_{L^\infty\tp{\R^5}}\le C\tnorm{\Phi_n - \Phi_m}_{C^1_{\m{b}}\tp{B_{\R^5}[0,3]}}
  \end{equation}
  is true for every $n,m\in\N$. We therefore learn that it suffices to verify the desired Cauchy condition on the compact domain $B_{\R^5}[0,3]$.

  From the uniform bounds~\eqref{key_uniform_bounds}, the Arzel\`a-Ascoli theorem, and the reflexivity of $H^2\tp{\R^5}$, we find that any subsequence of $\tcb{\Phi_n}_{n\in\N}$ has a further subsequence that is convergent in $C^1_{\m{b}}\tp{B_{\R^5}[0,3]}$ and weakly convergent in $H^2\tp{\R^5}$. Let $\tcb{\Psi_n}_{n\in\N}$ denote any such convergent subsequence of a subsequence and let $\Psi_\infty\in C_{\m{b}}^1\tp{\R^5}$ be its limit. From~\eqref{_reduction_estimate_} and the uniform bounds~\eqref{key_uniform_bounds}, we deduce also that $\tcb{\Psi_n}_{n\in\N}$ is convergent to $\Psi_\infty$ in the space $\mathsf{H}\tp{\R^5}$, defined in~\eqref{yet_another_space}, and we have the inclusions $\Psi_\infty\in\tp{\tilde{\mathsf{H}}\cap H^2}\tp{\R^5}$ with $\tnorm{\Psi_\infty}_{\tilde{\mathsf{H}}\tp{\R^5}} + \tnorm{\Psi_{\infty}}_{H^2\tp{\R^5}}\le C$. 

  By passing to the subsequential limit in the qualitative properties of the fifth item of Corollary~\ref{coro on solutions to the Bragg-Hawthorne equation} we learn also that
  \begin{equation}
    \Psi_\infty\ge0\text{ on }\R^5,\;\Psi_\infty\ge1/54\text{ on }B_{\R^5}[0,3],\;\text{and }\forall\;y\in\R^5\text{ with }y_5>0\text{ we have }\pd_5\Psi_\infty\tp{y}\le0.
  \end{equation}
  In particular $\Psi_\infty$ does not vanish identically. In fact, a straightforward maximum principle argument similar to what is seen in the next step of this proof will show that $\Psi_\infty\tp{y}>0$ for all $y\in\R^5$.
  
  We now claim that actually it is the case that 
    \begin{equation}\label{important_positivity}
      \pd_5\Psi_\infty\tp{y}<0\text{ for all }y\in\R^5\text{ satisfying }y_5>0.
    \end{equation}
  We shall, in proving~\eqref{important_positivity}, follow an argument similar to that of the proof of Theorem 3.4 in Amick and Turner~\cite{MR929976}. As a first step, replace the test function $\psi$ in~\eqref{the_key_weak_identity} with $-\pd_5\psi$ and integrate by parts on the first and last term. The identity
  \begin{multline}\label{the_other_key_weak_identity}
    \int_{\R^5}\grad\pd_5\Psi_n\tp{y}\cdot\grad\psi\tp{y}\;\m{d}\mathcal{H}^5\tp{y} + \int_{\R^5}\f{1}{\tabs{\P y}^2}\tp{S_n\cdot S_n'}\tp{\tabs{\P y}^2\tp{\Psi_n\tp{y} - 1} - k}\pd_5\psi\tp{y}\;\m{d}\mathcal{H}^5\tp{y}\\=\int_{\R^5}\tabs{\P y}^2\pd_5\Psi_n\tp{y}\tp{\Gamma_n}'\tp{\tabs{\P y}^2\tp{\Psi_n\tp{y} - 1} - k}\psi\tp{y}\;\m{d}\mathcal{H}^5\tp{y}
  \end{multline}
  is the result. Now suppose that $\psi\ge0$ and $y\in\supp\psi$ implies that $y_5>0$. Then, from $\Gamma_n'\ge0$ (see~\eqref{reg_vort_func}), and $\pd_5\Psi_n<0$ on $\supp\psi$, we deduce that the right hand side of the above equation is nonpositive. The limit as $n\to\infty$ exists for the first two integrals in~\eqref{the_other_key_weak_identity}; the latter vanishes by the inequality $|S_n\cdot S_n'|\le C\del_n^2$. Therefore we find that
  \begin{equation}
    \int_{\R^5}\grad\pd_5\Psi_\infty\tp{y}\cdot\grad\psi\tp{y}\;\m{d}\mathcal{H}^5\tp{y}\le0.
  \end{equation}
  Thus $\pd_5\Psi_\infty$ is weakly subharmonic on $\tcb{y\in\R^5\;:\;y_5>0}$. By $e_5$-even symmetry, we also know that $\pd_5\Psi_\infty = 0$ on the set $\tcb{y\in\R^5\;:\;y_5 = 0}$.
  
  We are in a position to invoke the strong maximum principle, see, e.g., Theorem 8.19 in Gilbarg and Trudinger~\cite{MR1814364}, on the function $\pd_5\Psi_\infty$. If there were a point $y\in\R^5$ with $y_5>0$ such that $\pd_5\Psi_\infty\tp{y} = 0$ (making $y$ an interior maximizer), then necessarily $\pd_5\Psi_{\infty}(x) = 0$ for all $x\in\R^5$ with $x_5>0$. But this is impossible since $\Psi_\infty$ is positive somewhere in the upper half space and vanishes at infinity. We have shown~\eqref{important_positivity} to be true.

  We are now ready to pass to the limit $n\to\infty$ in the final term in the weak formulation of equation~\eqref{the_key_weak_identity} for an arbitrarily fixed $\psi\in C^\infty_{\m{c}}\tp{\R^5}$. Let us fix $\upepsilon\in\tp{0,1}$. Consider the region
  \begin{equation}\label{the-good-region}
    E\tp{\upepsilon} = \tcb{y\in\R^5\;:\;|y|\le 3,\;|y_5|\ge\upepsilon,\;|\P y|^2\ge\upepsilon}.
  \end{equation}
  We note that $\mathcal{H}^5\tp{B_{\R^5}[0,3]\setminus E\tp{\upepsilon}}\le C\upepsilon$ and hence, by support considerations, 
  \begin{multline}\label{post_support_consideration}
    \babs{\int_{\R^5}\tp{\Gamma_n\tp{\tabs{\P y}^2\tp{\Psi_n\tp{y} - 1} - k} - 15\cdot\mathds{1}_{\tp{0,\infty}}\tp{\tabs{\P y}^2\tp{\Psi_\infty\tp{y} - 1} - k}}\psi\tp{y}\;\m{d}\mathcal{H}^5\tp{y}}\le C\upepsilon\\+\babs{\int_{E\tp{\upepsilon}}\tp{\Gamma_n\tp{\tabs{\P y}^2\tp{\Psi_n\tp{y} - 1} - k} - 15\cdot\mathds{1}_{\tp{0,\infty}}\tp{\tabs{\P y}^2\tp{\Psi_\infty\tp{y} - 1} - k}}\psi\tp{y}\;\m{d}\mathcal{H}^5\tp{y}}.
  \end{multline}

  Now let us write
  \begin{equation}
    \tilde{\Psi}_n\tp{y} = \tabs{\P y}^2\tp{\Psi_n\tp{y} - 1} - k\text{ and } \tilde{\Psi}_\infty\tp{y} = \tabs{\P y}^2\tp{\Psi_\infty\tp{y} - 1} - k
  \end{equation}
  for $y\in\R^5$ and $n\in\N$. From the convergence in $\mathsf{H}\tp{\R^5}$ of the sequence $\tcb{\Psi_n}_{n\in\N}$ and the positivity described in equation~\eqref{important_positivity} we deduce the existence of $\upnu>0$ and $\aleph\in\N$, depending on $\upepsilon$, such that 
  \begin{equation}
    \forall\;\N\ni n\ge\aleph,\;\forall\;y\in E\tp{\upepsilon},\;\min\tcb{\tabs{\pd_5\tilde{\Psi}_\infty\tp{y}},\tabs{\pd_5\tilde{\Psi}_n\tp{y}}}\ge\upnu.
  \end{equation}
  Then, for $\N\ni n\ge\aleph$, we may bound
  \begin{multline}\label{_more_analysis_}
    \babs{\int_{E\tp{\upepsilon}}\tp{\Gamma_n\tp{\tilde{\Psi}_n} - 15\cdot\mathds{1}_{\tp{0,\infty}}\tp{\tilde{\Psi}_\infty}}\psi}\\
    \le\babs{\int_{E\tp{\upepsilon}}\tp{\Gamma_n - 15\cdot\mathds{1}_{\tp{0,\infty}}}\tp{\tilde{\Psi}_n}\psi}+15\babs{\int_{E\tp{\upepsilon}}\tp{\mathds{1}_{\tp{0,\infty}}\tp{\tilde{\Psi}_n} - \mathds{1}_{\tp{0,\infty}}\tp{\tilde{\Psi}_\infty}}\psi}\\
    \le C\tnorm{\psi}_{L^\infty\tp{\R^5}}\tp{\mathcal{H}^5\tp{E\tp{\upepsilon}\cap\tcb{\tabs{\tilde{\Psi}_n}\le\ep_n}} + \mathcal{H}^5\tp{E\tp{\upepsilon}\cap\tcb{\tabs{\tilde{\Psi}_\infty}\le\tabs{\tilde{\Psi}_n - \tilde{\Psi}_\infty}}}}\\
    \le C\tnorm{\psi}_{L^\infty\tp{\R^5}}\upnu^{-1}\tp{\ep_n + \tnorm{\tilde{\Psi}_n - \tilde{\Psi}_\infty}_{L^\infty\tp{\R^5}}}.
  \end{multline}
  
  We combine~\eqref{post_support_consideration} and~\eqref{_more_analysis_} and deduce that 
  \begin{equation}
    \limsup_{n\to\infty}\babs{\int_{\R^5}\tp{\Gamma_n\tp{\tilde{\Psi}_n} - 15\cdot\mathds{1}_{\tp{0,\infty}}\tp{\tilde{\Psi}_\infty}}\psi}\le C\upepsilon.
  \end{equation}
  Upon recalling that $\upepsilon\in\tp{0,1}$ was selected arbitrarily, we use the above to deduce that the limit
  \begin{multline}
    \lim_{n\to\infty}\int_{\R^5}\Gamma_n\tp{\tabs{\P y}^2\tp{\Psi_n\tp{y} - 1} - k}\psi\tp{y}\;\m{d}\mathcal{H}^5\tp{y} = \int_{\R^5}15\cdot\mathds{1}_{\tp{0,\infty}}\tp{\tabs{\P y}^2\tp{\Psi_\infty\tp{y} - 1} - k}\psi\tp{y}\;\m{d}\mathcal{H}^5\tp{y}\\ = 15\int_{\R^5}\mathds{1}_{\mathsf{A}\tp{\Psi_\infty,k}}\tp{y}\psi\tp{y}\;\m{d}\mathcal{H}^5\tp{y}
  \end{multline}
  is valid.

  Hence, by sending $n\to\infty$ along the subsequence in equality~\eqref{the_key_weak_identity}, we deduce that $\Psi_\infty\in\mathcal{N}\tp{k}$, where the latter set was defined in Lemma~\ref{lem on Norbury uniqueness}. The aforementioned lemma, when combined with the uniform bounds~\eqref{key_uniform_bounds} and after increasing $C$, if necessary, so that $C\ge M_0$, provides the existence of some $0<k_0\le\del_\star$ such that if we take $0<k\le k_0$, then
  \begin{equation}\label{_napkin_on_a_sandwich_}
    \Psi_\infty\in\mathcal{N}\tp{k}\cap B_{\mathsf{H}\tp{\R^5}}\tsb{0,C}\text{ and }\m{card}\tp{\mathcal{N}\tp{k}\cap B_{\mathsf{H}\tp{\R^5}}\tsb{0,C}} \le 1.
  \end{equation}

  The right hand inequality in~\eqref{_napkin_on_a_sandwich_} shows that the function $\Psi_\infty$, which is the limit along a further subsequence of an arbitrary subsequence of the original sequence $\tcb{\Phi_n}_{n\in\N}$, is a unique object: let us denote it by $\Phi_\infty$. By the subsequence characterization of convergence, see, e.g., Proposition 1.1 of Chapter I in~\cite{MR1897317}, we conclude that $\lim_{n\to\infty}\Phi_n = \Phi_\infty$ in $\mathsf{H}\tp{\R^5}$. It is now clear that the third and fourth items of the proposition are satisfied.
\end{proof}

\subsection{Proofs of main results}\label{SS_proofs_main_results}

At last we complete the verification of our two main theorems, which are stated in Section~\ref{SS_main_results}.

\begin{proof}[Proof of Theorem~\ref{main_thm_1_swirl_free}]
  Let $k_0>0$ be as in Proposition~\ref{label prop on flexible desingularization of the hill norbury family} and set $k^\star = \min\tcb{\pmb{\mathtt{k}},k_0}$, where $\pmb{\mathtt{k}}$ is from Definition~\ref{defn on the Hill Norbury family}. Fix $0\le k\le k^\star$ and apply the aforementioned proposition with the element of $\mathscr{A}_0$ (see~\eqref{set_of_admissible_sequences}). Let $\tcb{\Phi_n}_{n\in\N}$ and $\tcb{\mathsf{p}_n}_{n\in\N}$ be the resulting sequences, and let $S_n$ and $\Gamma_n$ denote the corresponding circulation and head derivative functions. Since $\del_n=0$, equation~\eqref{circulation_and_head_derivative_functions_} gives $S_n=0$. Remarks~\ref{remark on instantiation of a solution map} and~\ref{remark_on_restriction_of_solution_map} and Proposition~\ref{prop on Euler reconstruction} therefore produce smooth functions $u_n$ and $p_n$ satisfying $u_n=\mathcal{U}[\Phi_n]$, where we recall that $\mathcal{U}[\cdot]$ is defined in~\eqref{_MER_VELO_MAP_}.

  The fourth item of Proposition~\ref{label prop on flexible desingularization of the hill norbury family}, together with Lemma~\ref{lem on Norbury uniqueness} when $k>0$ and Definition~\ref{defn on the Hill Norbury family}, identifies $\mathcal{U}[\Phi_\infty]=u^k$. The third and fourth items of the aforementioned proposition, together with~\eqref{mer_map_props}, \eqref{hypothesis_bounds_and_limits}, and~\eqref{the_conclusion_limit}, give the uniform $\tp{L^2\cap\m{LL}}\tp{\R^3}$ bound and the claimed $\tp{L^2\cap C^\al_{\m b}}\tp{\R^3}$ convergence. All of the remaining assertions follow directly from Proposition~\ref{prop on Euler reconstruction} and Lemma~\ref{lem_on_on_the_vortex_core}, after increasing the uniform constant $C$, if necessary.
\end{proof}

\begin{proof}[Proof of Theorem~\ref{main_thm_2_with_swirl}]
  We may increase the constant $C$ from Theorem~\ref{main_thm_1_swirl_free}, without invalidating that result, so that it dominates all of the constants below, $C\ge\kappa$, and $C^{-1}\le\del_\star^2$. Fix $0\le k\le k^\star$ and sequences as in the statement. Then $\tp{\tcb{h_n}_{n\in\N},\tcb{m_n}_{n\in\N}}\in\mathscr{A}_1$, where the latter set is defined in~\eqref{set_of_admissible_sequences}, and so Proposition~\ref{label prop on flexible desingularization of the hill norbury family} provides sequences $\tcb{\Phi_n}_{n\in\N}$ and $\tcb{\mathsf{p}_n}_{n\in\N}$ with
  \begin{equation}
    \del_n\eta_n=h_n,\;\del_n\mu_n=m_n,\text{ and }\del_n\to0.
  \end{equation}
  Let $\lambda_n$, $\rho_n$, $S_n$, and $\Gamma_n$ be the associated parameters and functions from Corollary~\ref{coro on solutions to the Bragg-Hawthorne equation}, and let $u_n,p_n$ be given by Proposition~\ref{prop on Euler reconstruction}. Remark~\ref{rmk_on_exact_prescribed_invariants}, Proposition~\ref{prop on Euler reconstruction}, and Lemma~\ref{lem_on_on_the_vortex_core} give the asserted invariants, equations, symmetries, signs, and support properties. Moreover, since $\del_n>0$, $\rho_n\ge22$, and the two summands defining $S_n$ have disjoint supports,~\eqref{_remarks_on_the_comb_functions_} gives $S_n\not\equiv0$ and hence $e_\theta\cdot u_n\not\equiv0$.

  From~\eqref{circulation_and_head_derivative_functions_}, $|\lambda_n|\le2\mathfrak{C}$, and~\eqref{_more_uniform_estimates_on_the_comb_functions_}, for every $0\le\al<1$ we have
  \begin{equation}
    \tnorm{S_n}_{C^0_{\m b}(\R)}+\tnorm{S_n}_{\m{LL}(\R)}\le C\del_n\text{ and }\tnorm{S_n}_{C^\al_{\m b}(\R)}\le C\del_n\tp{1+\tabs{\log\tabs{1-\al}}}.
  \end{equation}
  Lemma~\ref{lem_on_azimuthal_velocity_estimates} therefore shows that the azimuthal components are uniformly bounded in $\tp{L^2\cap\m{LL}}\tp{\R^3}$ and converge to zero in $\tp{L^2\cap C^\al_{\m b}}\tp{\R^3}$ for every $0<\al<1$. The meridional components converge to $u^k$ and obey the corresponding uniform bounds exactly as in the proof of Theorem~\ref{main_thm_1_swirl_free}. Since $u_n=\mathcal{U}[\Phi_n]+\tp{e_\theta\otimes e_\theta}u_n$, where we recall that $\mathcal{U}[\cdot]$ is defined in~\eqref{_MER_VELO_MAP_}, the claimed norm conclusions follow. Finally,~\eqref{SWIRL_SUPPORT} and $C\ge\kappa$ give the asserted uniform location of the swirl.
\end{proof}

\appendix
\section{Estimates on level sets}\label{appendix_estimates_on_level_sets}

The purpose of this appendix is to record some quantitative and rather technical information concerning regular level sets, tubular neighborhoods, and the geometry of $C^1$ submanifolds. Some of the ideas implemented below are closely related to classical ones; see, for instance, \cite{MR2977424,MR448362} for similar regular level set and tubular neighborhood constructions and~\cite{MR257325,MR3409135} for the relevant geometric measure theory.

We include the unabridged details here since our application necessitates quite specific results, especially ones that are suitably uniform under small $C^{1}$ perturbations of defining functions of hypersurfaces, together with numerous quantitative estimates adapted to the transverse geometry used within the body of the paper.

\subsection{Level set coordinates}\label{SS_level_set_coordinates}

We now instantiate this subsection's main object of study. Fix an ambient spatial dimension $\N\ni n\ge 2$ and let $U\subset\R^n$ be an open, bounded, and nonempty set. Suppose that we have a function of class $C^1$ $\Psi:\Bar{U}\to\R$ such that 
\begin{equation}\label{this is our base manifold}
    \es\neq M = \tcb{x\in U\;:\;\Psi(x) = 0}\Subset U\text{ and }\min_M\tabs{\grad\Psi(x)}>0.
\end{equation}
We also give notation for the modulus of continuity for $\grad\Psi$ on $\Bar{U}$; more precisely, let 
\begin{equation}\label{this is our modulus of continuity}
    \omega:[0,\infty)\to[0,\infty)\text{ be defined via }\omega\tp{\tau} = \max\tcb{\tabs{\grad\Psi(x) - \grad\Psi(y)}\;:\;x,y\in\Bar{U}\text{ and }\tabs{x - y}\le\tau}.
\end{equation}
We observe that, since $\Bar{U}$ is compact and $\grad\Psi$ is continuous on this set, $\omega$ is nondecreasing and satisfies $\omega\tp{\tau}\to0$ as $0<\tau\to0$.

A simple consequence of the implicit function theorem is that the set $M$ defined in~\eqref{this is our base manifold} is a compact, class $C^1$, $(n-1)$-dimensional submanifold of $\R^n$. By making small $C^1$ perturbations of the function $\Psi$, we expect to obtain small perturbations of the manifold $M$. The next few results allow us to concretely realize these perturbations with a version of graphical coordinates over $M$.

Since the hypersurface $M$ is only assumed to have $C^1$ regularity, the unit normal may fail to be differentiable or even Lipschitz. The next result in particular constructs a suitable regularization of the unit normal.

\begin{lem}[Regularization of the unit normal]\label{lem on regularization of the unit normal}
    There exists an open set $W\Subset U$ with $M\subset W$, a vector field $V\in C^1_{\m{b}}\tp{\Bar{W};\R^n}$, and a constant $0<\ell<\min\tcb{1,\m{dist}\tp{\pd W,\pd U}}$ such that the following hold.
    \begin{enumerate}
        \item For all $x\in W$ we have that $3/4\le\tabs{V(x)}\le 1$.
        \item For all $x,y\in W$ we have that $\tabs{V(x) - V(y)}\le\tp{n/\ell}\tabs{x - y}$.
        \item For all $x\in W$ we have that 
        \begin{equation}\label{youre_on_the_internet}
            \f45\min_{M}\tabs{\grad\Psi}\le\tabs{\grad\Psi(x)}\le\f54\max_M\tabs{\grad\Psi}.
        \end{equation}
        \item For all $x,y\in W$ that satisfy $\tabs{x - y}\le\ell$ we have the bounds 
        \begin{equation}\label{looking_at}
            \f34\min_M\tabs{\grad\Psi}\le V\tp{x}\cdot\grad\Psi(y)\le\f54\max_M\tabs{\grad\Psi}.
        \end{equation}
    \end{enumerate}
\end{lem}
\begin{proof}
      By the latter hypothesis of equation~\eqref{this is our base manifold} we are assured the existence of an open set $M\subset W\Subset U$ such that the stronger estimates
      \begin{equation}
        \f78\min_M\tabs{\grad\Psi}\le\tabs{\grad\Psi(x)}\le\f98\max_M\tabs{\grad\Psi}\text{ for all }x\in W
      \end{equation}
      hold; in particular, the estimates of equation~\eqref{youre_on_the_internet} hold.

  Let
  \begin{equation}\label{get_noted_later}
    0<\ell_0<\min\tcb{1,\m{dist}\tp{\pd W,\pd U},\m{dist}\tp{M,\pd W}}
  \end{equation}
  and set $\tilde{W} = \tcb{x\in U\;:\;\m{dist}\tp{x,W}\le\ell_0}$. The parameter $\ell_0$ can be chosen sufficiently small so that
  \begin{equation}
    \min_{\tilde{W}}\tabs{\grad\Psi}\ge\min_M\tabs{\grad\Psi}/2.
  \end{equation}
  We then define, for some to be determined parameter $0<\ell\le\ell_0$, the vector field
  \begin{equation}\label{initial_definition_of_the_vector_field_V}
    V\in C_{\m{b}}^1\tp{\Bar{W};\R^n}\text{ via }V(x) = \f{\Upgamma\tp{n/2 + 1}}{\tp{\sqrt{\pi}\ell}^n}\int_{B(x,\ell)}\f{\grad\Psi(y)}{\tabs{\grad\Psi(y)}}\;\m{d}y,\text{ for all }x\in\Bar{W}.
  \end{equation}
  
  By arguing directly from formula~\eqref{initial_definition_of_the_vector_field_V}, we may deduce the bound 
  \begin{equation}\label{e_2}
    \sup_{x\in W}\babs{V(x) - \f{\grad\Psi(x)}{\tabs{\grad\Psi(x)}}}\le2\f{\max_{\tilde{W}}\tabs{\grad\Psi}}{\min_{\tilde{W}}\tabs{\grad\Psi}^2}\omega\tp{\ell},
  \end{equation}
  where $\omega$ is the modulus of continuity defined in~\eqref{this is our modulus of continuity}. Since $\omega\tp{\tau}\to0$ as $\tau\to0$, we may take $0<\ell\le\ell_0$ smaller, if necessary, to ensure that
  \begin{equation}\label{e_3}
    \sup_{x\in W}\babs{V(x) - \f{\grad\Psi(x)}{\tabs{\grad\Psi(x)}}}\le\f18\text{ and hence }\forall\;x\in W,\;\f78\le\tabs{V(x)}.
  \end{equation}
  That $\tabs{V\tp{x}}\le 1$ for all $x\in W$ is immediate from its definition~\eqref{initial_definition_of_the_vector_field_V}. This first item is now established.

  The vector field $V$ is also manifestly Lipschitz; indeed, arguing directly from~\eqref{initial_definition_of_the_vector_field_V} we obtain, for all $x,y\in W$,
  \begin{equation}
    \tabs{V(x)-V(y)}\le\f{\mathcal{H}^n\tp{\tp{B(x,\ell)\setminus B(y,\ell)}\cup\tp{B(y,\ell)\setminus B(x,\ell)}}}{\mathcal{H}^n\tp{B(0,\ell)}}\le\f{n}{\ell}\tabs{x-y}.
  \end{equation}
  This is the second item.

  To prove the fourth item, we begin by noting that for any $x,y\in W$ with $|x - y|\le\ell$ we have the identity
  \begin{equation}
    V(x)\cdot\grad\Psi(y) - \tabs{\grad\Psi(x)} = V(x)\cdot\tp{\grad\Psi(y) - \grad\Psi(x)} + \bp{V(x) - \f{\grad\Psi(x)}{\tabs{\grad\Psi(x)}}}\cdot\grad\Psi(x)
  \end{equation}
  and hence, from~\eqref{e_2} and~\eqref{e_3}, the estimate 
  \begin{equation}\label{e_4}
    \tabs{V(x)\cdot\grad\Psi(y) - \tabs{\grad\Psi(x)}}\le\f54\omega\tp{\ell} + 2\bp{\f{\max_{\tilde{W}}\tabs{\grad\Psi}}{\min_{\tilde{W}}\tabs{\grad\Psi}}}^2\omega\tp{\ell}
  \end{equation}
  is valid. By taking $\ell>0$ smaller, if necessary, we can arrange from~\eqref{e_4} that
  \begin{equation}
    \tabs{V(x)\cdot\grad\Psi(y) - \tabs{\grad\Psi(x)}}\le\f18\min_M\tabs{\grad\Psi}
  \end{equation} 
  and so estimate~\eqref{looking_at} holds.
\end{proof}

The next few results will reference the following team of strictly positive parameters that are defined in terms of $M$, $\grad\Psi$, $\omega$, and $\ell$. We set
\begin{equation}\label{parameter_splurge}
    \be = \f{\max\tcb{1,\max_{M}\tabs{\grad\Psi}}}{\min\tcb{1,\min_M\tabs{\grad\Psi}}},\;\al = \f{1}{2^{10}\be},\;\lambda = \al^2\cdot\sup\tcb{0\le\tau\le\ell/n\;:\;\omega\tp{\tau}\le\al^2},\;\ep = \tp{\al\lambda}^2,\;r = \tp{\al\lambda}^3.
\end{equation}
It is helpful to note that from~\eqref{get_noted_later} we have $\ell\le\min\tcb{1,\m{dist}\tp{M,\pd W}}$. We shall also define the following tubular neighborhoods of $M$:
\begin{equation}\label{tub_nei}
    W_\mu = \tcb{x\in\R^n\;:\;\m{dist}\tp{x,M}\le\mu}\Subset U\text{ for any }0<\mu\le\ell.
\end{equation}

Our next result uses the vector field $V$ from Lemma~\ref{lem on regularization of the unit normal} to build a family of functions that serve as a suitable proxy for the signed distance to the hypersurface $M$ and its small perturbations.

\begin{lem}[Regularized distance to the manifold]\label{lem on regularized distance to the manifold}
    Define the complete metric spaces
    \begin{equation}\label{complete_metric_spaces}
        \mathcal{X} = \tcb{g\in C^0_{\m{b}}\tp{W_\ep\times[-\ep,\ep]}\;:\;\tnorm{g}_{C^0_{\m{b}}\tp{W_\ep\times[-\ep,\ep]}}\le\lambda}\text{ and }\mathcal{B} = \tcb{\phi\in C^1_{\m{b}}\tp{\Bar{U}}\;:\;\tnorm{\phi}_{C^{0,1}\tp{\Bar{U}}}\le r}.
    \end{equation}
    There exists a Lipschitz continuous mapping $\bf{f}:\mathcal{B}\to\mathcal{X}$ such that the following hold for all $\phi\in\mathcal{B}$ and all $\tp{x,t}\in W_\ep\times[-\ep,\ep]$.
    \begin{enumerate}
        \item Upon defining
        \begin{equation}\label{occasionally_useful}
             G_\phi\tp{x,t} = x + V(x)\tsb{\bf{f}\tp{\phi}}\tp{x,t},\text{ we have the inclusion }G_\phi(x,t)\in W_{2^{-37}\ell}.
        \end{equation}
        \item We have the equality
        \begin{equation}\label{projection_parametrization_identity}
            \tp{\Psi + \phi}\circ\tp{x + V(x)\tsb{\bf{f}\tp{\phi}}\tp{x,t}} = t.
        \end{equation}
        \item We have the estimate
        \begin{equation}\label{distance to the manifold bound}
            \f{1}{2\be}\tabs{\Psi(x) + \phi(x) - t}\le\tabs{\tsb{\bf{f}\tp{\phi}}\tp{x,t}}\le2\be\tabs{\Psi(x) + \phi(x) - t}.
        \end{equation}
        \item If $\tp{\tilde{x},\tilde{t}}\in W_\ep\times[-\ep,\ep]$ satisfy $\tabs{x - \tilde{x}}\le\lambda/2$, we have the Lipschitz estimate
        \begin{equation}\label{initial_lip_est}
            \tabs{\tsb{\bf{f}\tp{\phi}}\tp{x,t}-\tsb{\bf{f}\tp{\phi}}\tp{\tilde{x},\tilde{t}}}\le6\be\tabs{x - \tilde{x},t - \tilde{t}}.
        \end{equation}
        \item We have that $\bf{f}\tp{\phi}\in C^1_{\m{b}}\tp{W_\ep\times[-\ep,\ep]}$; moreover the mapping
        \begin{equation}
            \mathcal{B}\ni\tilde{\phi}\mapsto \bf{f}\tp{\tilde{\phi}} \in C^1_{\m{b}}\tp{W_\ep\times[-\ep,\ep]}\text{ is continuous}.
        \end{equation}
    \end{enumerate}
\end{lem}
\begin{proof}
    We argue via the contraction mapping principle. Seeking a solution to the level set projection-parametrization identity~\eqref{projection_parametrization_identity}, we introduce the operator $\aleph:\mathcal{B}\times\mathcal{X}\to\mathcal{X}$  with definition
    \begin{equation}\label{e_6}
    \tsb{\aleph\tp{\phi,g}}\tp{x,t} = \f{1}{V(x)\cdot\grad\Psi\tp{x}}\tp{t - \phi\tp{x + g(x,t)V(x)} + \grad\Psi\tp{x}\cdot V(x)g(x,t) - \Psi(x + g(x,t)V(x))}
  \end{equation}
  for $\tp{x,t}\in W_\ep\times[-\ep,\ep]$ and $\tp{\phi,g}\in\mathcal{B}\times\mathcal{X}$. The properties of $V$ from Lemma~\ref{lem on regularization of the unit normal} and the selection of parameters~\eqref{parameter_splurge} are much more conservative than what is required to ensure that $\aleph$ is well-defined, in the sense that all of the compositions in its definition are permissible and it maps into the stated codomain of $\mathcal{X}$. More precisely, it is straightforward to verify that 
  \begin{equation}\label{sup_inf}
    \tnorm{\aleph\tp{\phi,g}}_{C_{\m{b}}^0\tp{W_\ep\times[-\ep,\ep]}}\le\al\lambda,
  \end{equation}
    \begin{equation}
    \tnorm{\aleph\tp{\phi,g} - \aleph\tp{\phi,\tilde{g}}}_{C_{\m{b}}^0\tp{W_\ep\times[-\ep,\ep]}}\le\al\tnorm{g - \tilde{g}}_{C^0_{\m{b}}\tp{W_\ep\times[-\ep,\ep]}},
  \end{equation}
  and
  \begin{equation}
    \tnorm{\aleph\tp{\phi,g} - \aleph\tp{\tilde{\phi},g}}_{C^0_{\m{b}}\tp{W_\ep\times[-\ep,\ep]}}\le\tp{4\be/3}\tnorm{\phi - \tilde{\phi}}_{C^0_{\m{b}}\tp{\Bar{U}}},
  \end{equation}
  for all $\tp{\phi,g},\tp{\tilde{\phi},\tilde{g}}\in\mathcal{B}\times\mathcal{X}$.

  The contraction mapping principle with parameter thus applies, see, e.g., Theorem C.7 in Irwin~\cite{MR1867353}, and we acquire a Lipschitz continuous mapping
  \begin{equation}\label{e_8}
    \bf{f}:\mathcal{B}\to\mathcal{X}\text{ such that }\aleph\tp{\phi,\bf{f}\tp{\phi}} = \bf{f}\tp{\phi},\;\forall\;\phi\in\mathcal{B}.
  \end{equation}
  By unpacking the meaning of this fixed point from the definition~\eqref{e_6}, we deduce that precisely the second item holds. 

  Let us now turn our attention to the third item. From~\eqref{e_6} we deduce, for $\phi\in\mathcal{B}$ and $\tp{x,t}\in W_\ep\times\tsb{-\ep,\ep}$ the estimates
  \begin{equation}
    \tabs{\tsb{\bf{f}\tp{\phi}}\tp{x,t}}\le\f{4\be}{3}\tp{\tabs{\Psi(x) + \phi(x) - t} + \tp{r + \omega\tp{\lambda}}\tabs{\tsb{\bf{f}\tp{\phi}}\tp{x,t}}}
  \end{equation}
  and
  \begin{equation}
    \tabs{\tsb{\bf{f}\tp{\phi}}\tp{x,t}}\ge\f{4}{5\be}\tp{\tabs{\Psi(x) + \phi(x) - t} - \tp{r + \omega\tp{\lambda}}\tabs{\tsb{\bf{f}\tp{\phi}}\tp{x,t}}}
  \end{equation}
  The two sided bounds claimed in~\eqref{distance to the manifold bound} now follow. As a special consequence of these estimates, we can establish the inclusion asserted within the first item. Indeed, if $x_0\in M$ is such that $|x - x_0|\le\ep$, then
  \begin{equation}
    \tabs{x + V(x)\tsb{\bf{f}\tp{\phi}}\tp{x,t} - x_0}\le\ep + 2\be\tabs{\Psi(x) + \phi(x) - t}<2^{-37}\ell
  \end{equation}
  Thus~\eqref{occasionally_useful} holds.

  Next, we justify the fourth item. For any $\phi\in\mathcal{B}$ and $\tp{x,t},\tp{\tilde{x},\tilde{t}}\in W_\ep\times[-\ep,\ep]$ satisfying $|x - \tilde{x}|\le\lambda/2$ the line segment connecting $G_\phi\tp{x,t}$ and $G_\phi\tp{\tilde{x},\tilde{t}}$ (where $G_\phi$ is defined in~\eqref{occasionally_useful}) is contained within the interior of $W_{\ell}$. Thus we can use the fundamental theorem of calculus when examining the difference between~\eqref{projection_parametrization_identity} evaluated at $(x,t)$ and $(\tilde{x},\tilde{t})$. The result of doing so is the identity
  \begin{multline}\label{key_difference_identity}
    t - \tilde{t} = \tp{\tp{x - \tilde{x}} + \tp{V(x) - V(\tilde{x})}\tsb{\bf{f}\tp{\phi}}\tp{x,t}}\cdot\int_0^1\grad\tp{\Psi + \phi}\tp{\tau G_\phi(x,t) + \tp{1 - \tau}G_\phi(\tilde{x},\tilde{t})}\;\m{d}\tau\\
    +\tp{\tsb{\bf{f}\tp{\phi}}\tp{x,t} - \tsb{\bf{f}\tp{\phi}}\tp{\tilde{x},\tilde{t}}}V\tp{\tilde{x}}\cdot\bp{\int_0^1\grad\tp{\Psi + \phi}\tp{\tau G_\phi(x,t) + \tp{1 - \tau}G_\phi(\tilde{x},\tilde{t})}\;\m{d}\tau} - \grad\Psi\tp{G_\phi\tp{x,t}}\\
    +\tp{\tsb{\bf{f}\tp{\phi}}\tp{x,t} - \tsb{\bf{f}\tp{\phi}}(\tilde{x},\tilde{t})}V\tp{\tilde{x}}\cdot\grad\Psi\tp{G_\phi\tp{x,t}}.
  \end{multline}
  In turn, we get the estimates
  \begin{multline}\label{estimates_in_turn}
    \tabs{\tsb{\bf{f}\tp{\phi}}\tp{x,t} - \tsb{\bf{f}\tp{\phi}}(\tilde{x},\tilde{t})}\le\tp{4\be/3}\tp{\omega\tp{\lambda} + r}\tabs{\tsb{\bf{f}\tp{\phi}}\tp{x,t} - \tsb{\bf{f}\tp{\phi}}\tp{\tilde{x},\tilde{t}}} \\+ \be\tp{5/3 + 4r/3}\tp{1 + n\al\lambda/\ell}\tabs{x - \tilde{x}} + \tp{4\be/3}\tabs{t - \tilde{t}}
  \end{multline}
  and hence we get the Lipschitz bound~\eqref{initial_lip_est}.
  
  We may exploit identity~\eqref{key_difference_identity} via a standard argument implementing difference quotients to deduce that $\bf{f}\tp{\phi}$ is actually differentiable at every point in $W_\ep\times[-\ep,\ep]$. The Lipschitz bounds~\eqref{initial_lip_est} directly lead to bounds on the differential. Furthermore, the mapping of $\phi\mapsto D\tsb{\bf{f}\tp{\phi}}$ is observed to be continuous. This completes the proof.
\end{proof}

We now wish to examine further the maps $G_\phi$ introduced in equation~\eqref{occasionally_useful}. The final target is to show that their restriction to $M\times\tcb{t}$ gives a diffeomorphism $M\times\tcb{t}\to W_\mu\cap\tcb{\Psi + \phi = t}$ for small enough $t$, $\phi$, and $\mu$. To proceed, we require the following essential result on almost tangent vectors.

\begin{lem}[Estimates on almost tangent vectors]\label{lem on estimates on almost tangent vectors}
    Suppose that $\tp{x,t},\tp{\tilde{x},\tilde{t}}\in W_\ep\times[-\ep,\ep]$ satisfy $|x - \tilde{x}|\le\lambda/2$. For any $z\in W_\ell$ satisfying $\tabs{z - \tilde{x}}\le\lambda$ and any $\phi\in\mathcal{B}$ we have the estimate
    \begin{equation}\label{the almost tangent vector estimate}
        \tabs{\grad\tp{\Psi+\phi}\tp{z}\cdot\tp{x - \tilde{x}}}\le\tp{\omega\tp{2\lambda} + 2r}\tabs{x - \tilde{x}} + \tabs{\tp{\Psi + \phi}(x) - \tp{\Psi + \phi}(\tilde{x})}.
    \end{equation}
\end{lem}
\begin{proof}
    The argument is very simple. From the fundamental theorem of calculus, we have the equality
    \begin{equation}
       \grad\tp{\Psi + \phi}(z)\cdot\tp{x - \tilde{x}} = \tp{\Psi + \phi}\tp{x} - \tp{\Psi + \phi}\tp{\tilde{x}} + \tp{x - \tilde{x}}\cdot\int_0^1\tp{\grad\tp{\Psi + \phi}\tp{z} - \grad\tp{\Psi + \phi}\tp{\tau x + \tp{1 - \tau}\tilde{x}}}\;\m{d}\tau.
    \end{equation}
    Taking the length and using $\tabs{z-\tp{\tau x+\tp{1-\tau}\tilde{x}}}\le2\lambda$ gives~\eqref{the almost tangent vector estimate}.
\end{proof}

While Lemma~\ref{lem on estimates on almost tangent vectors} may seem trivial, it is the key to the bi-Lipschitz bounds that follow. 

\begin{lem}[Bi-Lipschitz bounds]\label{lem on bilipschitz bounds}
  Regarding the map $G$ defined in the first item of Lemma~\ref{lem on regularized distance to the manifold}, the following hold for all $\phi\in\mathcal{B}$, $\tp{x,t},\tp{\tilde{x},\tilde{t}}\in W_\ep\times[-\ep,\ep]$ that satisfy $\tabs{x - \tilde{x}}\le\lambda/2$.
  \begin{enumerate}
    \item We have the Lipschitz bound
    \begin{equation}\label{LIP_1}
      \f{1}{2\be}\tabs{t - \tilde{t}}\le\tabs{G_\phi\tp{x,t} - G_\phi\tp{\tilde{x},\tilde{t}}}\le8\be\tabs{x - \tilde{x},t - \tilde{t}}.
    \end{equation}
    \item If there exists $y,\tilde{y}\in W_\ep$, $\psi\in\mathcal{B}$ with
    \begin{equation}\label{LIPO}
      x - \tilde{x} = y - \tilde{y},\;\tabs{\tilde{x} - \tilde{y}} + \tabs{x - \tilde{x}}\le\lambda/2,\;\tp{\Psi + \psi}\tp{y} = \tp{\Psi + \psi}\tp{\tilde{y}},
    \end{equation}
    then we have the reverse Lipschitz bound
    \begin{equation}\label{LIP_2}
      \tabs{x - \tilde{x},t - \tilde{t}}\le 8\be^2\tabs{G_\phi\tp{x,t} - G_\phi\tp{\tilde{x},\tilde{t}}}.
    \end{equation}
  \end{enumerate}
\end{lem}
\begin{proof}
    The right-hand estimate of the first item is a direct consequence of the Lipschitz bounds from the second item of Lemma~\ref{lem on regularization of the unit normal} and the fourth item of Lemma~\ref{lem on regularized distance to the manifold}.

    We turn our attention to the left hand inequality of~\eqref{LIP_1}. From identity~\eqref{projection_parametrization_identity}, we acquire that
    \begin{equation}\label{weve_been_here}
      t - \tilde{t} = \tp{\Psi + \phi}\circ G_\phi\tp{x,t}  - \tp{\Psi + \phi}\circ G_\phi\tp{\tilde{x},\tilde{t}}.
    \end{equation}
    As we have noted in the derivation of equation~\eqref{key_difference_identity}, the line segment connecting $G_\phi\tp{x,t}$ and $G_\phi\tp{\tilde{x},\tilde{t}}$ is contained entirely within the interior of $W_\ell$. Hence, the fundamental theorem of calculus takes us from identity~\eqref{weve_been_here} to the bound asserted on the left of~\eqref{LIP_1}.

    The second item is more subtle. We begin with the bound
    \begin{multline}\label{marco}
      \tabs{G_\phi\tp{x,t} - G_\phi\tp{\tilde{x},\tilde{t}}}\ge\tabs{x - \tilde{x}} - \tabs{\tsb{\bf{f}\tp{\phi}}\tp{x,t}}\tabs{V(x) - V(\tilde{x})} - \tabs{\tsb{\bf{f}\tp{\phi}}\tp{x,t} - \tsb{\bf{f}\tp{\phi}}\tp{\tilde{x},\tilde{t}}}\\
      \ge\tp{1 - \al^2}\tabs{x - \tilde{x}} - \tabs{\tsb{\bf{f}\tp{\phi}}\tp{x,t} - \tsb{\bf{f}\tp{\phi}}\tp{\tilde{x},\tilde{t}}}.
    \end{multline}
    To continue, we require a more precise Lipschitz estimate on $\bf{f}\tp{\phi}$. We examine identity~\eqref{key_difference_identity} more carefully. The following more precise version of~\eqref{estimates_in_turn} holds:
    \begin{multline}\label{estimates_in_turn_again}
      \tabs{\tsb{\bf{f}\tp{\phi}}\tp{x,t} - \tsb{\bf{f}\tp{\phi}}(\tilde{x},\tilde{t})}\le\f{4\be}{3}\tp{\omega\tp{\lambda} + r}\tabs{\tsb{\bf{f}\tp{\phi}}\tp{x,t} - \tsb{\bf{f}\tp{\phi}}\tp{\tilde{x},\tilde{t}}} +\f{4\be}{3}\bp{1 + \f{n\al\lambda}{\ell}}\tabs{x - \tilde{x}} \\+ \f{4\be}{3}\tabs{t - \tilde{t}} + \f{4\be}{3}\babs{\tp{x - \tilde{x}}\cdot\int_0^1\grad\tp{\Psi + \phi}\tp{\tau G_\phi\tp{x,t} + \tp{1 - \tau}G_\phi\tp{\tilde{x},\tilde{t}}}\;\m{d}\tau}.
    \end{multline}
    To estimate the final term above we now use the hypotheses~\eqref{LIPO} in conjunction with Lemma~\ref{lem on estimates on almost tangent vectors}. For every $0\le\tau\le1$, the point
    \begin{equation}
      z_\tau=\tau G_\phi\tp{x,t} + \tp{1 - \tau}G_\phi\tp{\tilde{x},\tilde{t}}
    \end{equation}
    belongs to $W_\ell$ and satisfies $\tabs{z_\tau-\tilde{y}}\le\lambda$. We may therefore apply the lemma to $y,\tilde{y},\psi$, and then estimate the difference between $\grad\phi$ and $\grad\psi$. The result is
    \begin{equation}\label{key_almost_tangent}
      \f{4\be}{3}\babs{\tp{x - \tilde{x}}\cdot\int_0^1\grad\tp{\Psi + \phi}\tp{z_\tau}\;\m{d}\tau}\le\f{4\be}{3}\tp{\omega\tp{2\lambda} + 4r}\tabs{x - \tilde{x}}.
    \end{equation}
    Notice that the term corresponding to the final right hand side contribution of~\eqref{the almost tangent vector estimate} vanishes identically. Since $\omega\tp{2\lambda}\le\al^2$, $\lambda\le\al^2$, and $r\le\al^9$, the same absorption in estimates~\eqref{estimates_in_turn_again} and~\eqref{key_almost_tangent} gives us the bound
    \begin{equation}\label{aguas_de}
      \tabs{\tsb{\bf{f}\tp{\phi}}\tp{x,t} - \tsb{\bf{f}\tp{\phi}}\tp{\tilde{x},\tilde{t}}}\le2\al\tabs{x - \tilde{x}} + 2\be\tabs{t - \tilde{t}}
    \end{equation}

    The combination of~\eqref{aguas_de}, \eqref{marco}, and the left-hand bound of~\eqref{LIP_1} gives
    \begin{equation}
      \tabs{t-\tilde{t}}\le2\be\tabs{G_\phi\tp{x,t}-G_\phi\tp{\tilde{x},\tilde{t}}}
    \end{equation}
    and
    \begin{equation}
      \tp{1-2\al-\al^2}\tabs{x-\tilde{x}}\le\tp{1+4\be^2}\tabs{G_\phi\tp{x,t}-G_\phi\tp{\tilde{x},\tilde{t}}}.
    \end{equation}
    Since $\be\ge1$ and $\al\le2^{-10}$, these inequalities yield the reverse Lipschitz inequality~\eqref{LIP_2} with the stated factor $8\be^2$.
\end{proof}

Before we come to the principal theorem of this subsection, we record some properties of the operators $G$ introduced in the first item of Lemma~\ref{lem on regularized distance to the manifold} that suggest that they behave like a nearest point projection.

\begin{lem}[Projection properties]\label{lem on projection properties}
  Define the class $C^1$ map
  \begin{equation}\label{proj_defn}
    \Pi:W_\ep\to W_{2^{-37}\ell}\text{ via }\Pi\tp{x} = x + V(x)\tsb{\bf{f}\tp{0}}\tp{x,0}\text{ for all }x\in W_\ep.
  \end{equation}
  The following hold.
  \begin{enumerate}
    \item For all $x\in W_\ep$ we have that $\Psi\circ\Pi\tp{x} = 0$ and hence $\Pi\tp{x}\in M$.
    \item If $x\in M$, then $\Pi\tp{x} = x$.
    \item If $x,\tilde{x}\in W_\ep$ satisfy $\tabs{x - \tilde{x}}\le\lambda/2$, then we have the estimate
    \begin{equation}\label{the projection estimate_}
      \tabs{x - \Pi\tp{\tilde{x}}}\le 7\be\tabs{x - \tilde{x}} + 2\be\tabs{\Psi(x)}.
    \end{equation}
  \end{enumerate}
\end{lem}
\begin{proof}
  The first item is immediate from the second item of Lemma~\ref{lem on regularized distance to the manifold} and the definition of $M$~\eqref{this is our base manifold}. 

  For the second item, we note that if $x\in M$, then $\Psi(x) = 0$ and so the third item of Lemma~\ref{lem on regularized distance to the manifold} implies that $\tsb{\bf{f}\tp{0}}\tp{x,0} = 0$. Hence $\Pi(x) = x$ directly from the definition~\eqref{proj_defn}.

  Finally, to prove the third item we begin by writing
  \begin{equation}
    x - \Pi\tp{\tilde{x}} = x - \tilde{x} - V\tp{\tilde{x}}\tp{\tsb{\bf{f}\tp{0}}\tp{\tilde{x},0} - \tsb{\bf{f}\tp{0}}\tp{x,0}} - V(\tilde{x})\tsb{\bf{f}\tp{0}}\tp{x,0}.
  \end{equation}
  We then measure the length and use the estimates from the third and fourth items of Lemma~\ref{lem on regularized distance to the manifold} to deduce the bound~\eqref{the projection estimate_}.
\end{proof}

We are now ready for the main result of this subsection.

\begin{thm}[Stable local level set coordinates]\label{thm on stable local level set coordinates}
  Consider the class $C^1$ maps
  \begin{equation}\label{The coordinate maps}
    \Theta_\phi:M\times[-\ep,\ep]\to W_{2^{-37}\ell}\text{ given via the restriction }\Theta_\phi\tp{m,t} = G_\phi\tp{m,t}
  \end{equation}
  and defined for $\phi\in\mathcal{B}$, $\tp{m,t}\in M\times[-\ep,\ep]$, where $G$ is from the first item of Lemma~\ref{lem on regularized distance to the manifold}. The following hold.
  \begin{enumerate}
    \item The maps $\Theta_\phi$ are injective; in fact, for all $\phi\in\mathcal{B}$ and all $\tp{m,t},\tp{\tilde{m},\tilde{t}}\in M\times[-\ep,\ep]$ we have the two-sided Lipschitz bound
    \begin{equation}\label{two_sided_lipschitz}
      \f{1}{8\be^2}\tabs{m - \tilde{m},t - \tilde{t}}\le\tabs{\Theta_\phi\tp{m,t} - \Theta_\phi\tp{\tilde{m},\tilde{t}}}\le 8\be\tabs{m - \tilde{m},t - \tilde{t}}.
    \end{equation}
    \item The maps $\Theta_\phi$ are uniformly locally surjective in the sense that 
    \begin{equation}\label{local_surjectivity}
      \forall\;\phi\in\mathcal{B},\;\Theta_\phi\tp{M\times[-\ep,\ep]}\supseteq W_{\tp{\al\lambda}^3}.
    \end{equation}
    \item Whenever $\phi\in\mathcal{B}$ and $t\in\R$ satisfy $\tnorm{\phi}_{C^{0,1}\tp{\Bar{U}}}\le\al\tp{\al\lambda}^3$ and $|t|\le\al\tp{\al\lambda}^3$ we have the inclusion
    \begin{equation}\label{level_set_covering}
      \tcb{x\in W_\ep\;:\;\Psi(x) + \phi(x) = t}\subset W_{\tp{\al\lambda}^3}.
    \end{equation}
    \item The mapping $\mathcal{B}\ni\tilde{\phi}\mapsto\Theta_{\tilde{\phi}}\in C^1_{\m{b}}\tp{M\times[-\ep,\ep];\R^n}$ is continuous.
  \end{enumerate}
\end{thm}
\begin{proof}
  We begin by remarking that the fourth item is a trivial consequence of the fifth item of Lemma~\ref{lem on regularized distance to the manifold}.

  Let us continue with the proof of the first item. The two-sided Lipschitz bound~\eqref{two_sided_lipschitz} was already established under the additional assumption that $\tabs{m - \tilde{m}}\le\lambda/2$, see Lemma~\ref{lem on bilipschitz bounds}. So we consider the case that $|m - \tilde{m}|>\lambda/2$. For the upper bound, we note that the definition of $\Theta_\phi$ and the uniform bound~\eqref{sup_inf} give: 
  \begin{equation}
    \tabs{\Theta_\phi\tp{m,t} - \Theta_\phi\tp{\tilde{m},\tilde{t}}}\le\tabs{m - \tilde{m}} + 2\al\lambda\le\tp{1 + 4\al}\tabs{m - \tilde{m}}.
  \end{equation}
  The same strategy gives the lower bound 
  \begin{equation}
    \tabs{\Theta_\phi\tp{m,t} - \Theta_\phi\tp{\tilde{m},\tilde{t}}}\ge\tabs{m - \tilde{m}} - 2\al\lambda\ge\tp{1 - 4\tp{\al + \al^2}}\tabs{m - \tilde{m}} + \tabs{t - \tilde{t}}.
  \end{equation}

  The second item requires a somewhat delicate argument. Let $\phi\in\mathcal{B}$ and let $y\in W_{\tp{\al\lambda}^3}$. By definition~\eqref{tub_nei}, there exists $x_0\in M$ such that $\tabs{x_0 - y}\le\tp{\al\lambda}^3$. Let us set $t = \tp{\Psi + \phi}\tp{y}$. We get the estimate 
  \begin{equation}\label{estimate_on_t}
    \tabs{t}\le r + \tabs{\Psi(y) - \Psi(x_0)}\le\tp{1 + 5\be/4}\tp{\al\lambda}^3<\ep.
  \end{equation}
  Now define the ball and the map
  \begin{equation}\label{map_to_be_fixed}
    B = \tcb{x\in\R^n\;:\;\tabs{x - x_0}\le\tp{\al\lambda}^2}\subset W_{\tp{\al\lambda}^2}\text{ and }F:B\to\R^n,\;F(x) = \Pi\tp{y - V\tp{\Pi\tp{x}}\tsb{\bf{f}\tp{\phi}}\tp{\Pi\tp{x},t}}
  \end{equation}
  for all $x\in B$. Recall that the map $\Pi$ appearing above is the projection operator introduced in Lemma~\ref{lem on projection properties}. To check that $F$ is indeed well-defined, we need that for all $x\in B$ the argument of the outermost $\Pi$ belongs to its domain $W_\ep$. By using the third item of Lemma~\ref{lem on regularized distance to the manifold} and the intermediate bound of~\eqref{estimate_on_t}, we get the estimate
  \begin{equation}
    \tabs{y - V\tp{\Pi\tp{x}}\tsb{\bf{f}\tp{\phi}}\tp{\Pi\tp{x},t} - x_0}\le\tp{\al\lambda}^3 + 2\be\tp{\tp{\al\lambda}^3 + |t|}\le\tp{\al\lambda}^3\tp{1 + 2\be\tp{2 + 5\be/4}}<\ep.
  \end{equation}
  So indeed $F$ is well-defined; this map is also continuous. We claim now that $F(B)\subseteq B$. To estimate the difference between $F(x)$ and $x_0$ for an arbitrary $x\in B$, we employ first the third item of Lemma~\ref{lem on projection properties}, thus
  \begin{equation}
    \tabs{x_0 - F(x)}\le 7\be\tabs{x_0 - y + V(\Pi(x))\tsb{\bf{f}\tp{\phi}}\tp{\Pi(x),t}}.
  \end{equation}
  Next, we use again the third item of Lemma~\ref{lem on regularized distance to the manifold} and the intermediate bound of~\eqref{estimate_on_t}, therefore 
  \begin{equation}
    \tabs{x_0 - F(x)}\le 7\be\tp{\al\lambda}^3 + 14\be^2\tp{1 + 5\be/4}\tp{\al\lambda}^3\le2^{-20}\tp{\al\lambda}^2.
  \end{equation}

  Brouwer's fixed point theorem thus applies for the map $F$. We therefore deduce the existence of $x\in B$ such that $x = F(x)$. From the definition of $F$~\eqref{map_to_be_fixed} and the first item of Lemma~\ref{lem on projection properties}, we see that the fixed point $x$ satisfies the inclusion $x\in M$. In turn, the second item of the aforementioned result implies that $x = \Pi(x)$. Thus, upon unpacking the fixed point identity $x = F(x)$ we learn that 
  \begin{equation}\label{step_on}
    \Pi(x) = \Pi\tp{y - V(x)\tsb{\bf{f}\tp{\phi}}\tp{x,t}}.
  \end{equation}
  We would like to conclude that the arguments of $\Pi$ in~\eqref{step_on} are also equal. The second item of Lemma~\ref{lem on bilipschitz bounds} is the key to doing so. Let us set
  \begin{equation}
    w = x,\;\tilde{w} = y - V(x)\tsb{\bf{f}\tp{\phi}}\tp{x,t},\;z = x + V(x)\tsb{\bf{f}\tp{\phi}}\tp{x,t},\;\tilde{z} = y.
  \end{equation}
  Then 
  \begin{equation}
    w - \tilde{w} = z - \tilde{z},\;\tabs{\tilde{w} - \tilde{z}} + \tabs{w - \tilde{w}}\le\lambda/2,\;\tp{\Psi + \phi}\tp{z} = t = \tp{\Psi + \phi}\tp{\tilde{z}}.
  \end{equation}
  Thus the reverse Lipschitz bound~\eqref{LIP_2} applies and we find that 
  \begin{equation}
    \tabs{x - y + V(x)\tsb{\bf{f}\tp{\phi}}\tp{x,t}} = \tabs{w - \tilde{w}}\le8\be^2\tabs{\Pi(w) - \Pi(\tilde{w})} = 0\text{ and hence }y = \Theta_\phi\tp{x,t}.
  \end{equation}

  Finally, we turn our attention to the third item. We first claim that for all $x\in W_\ep$ we have that 
  \begin{equation}\label{kitty_}
    \m{dist}\tp{x,M}\le2\be\tabs{\Psi(x)}.
  \end{equation}
  To prove~\eqref{kitty_}, we begin by letting $x\in W_\ep$ and selecting $x_0\in M$ such that $\tabs{x - x_0} = \m{dist}\tp{x,M}$. The class $C^1$ function $M\ni y\mapsto\tabs{x - y}^2\in\R$ has a minimum at $y = x_0$. Thus the differential vanishes for all directions belonging to the tangent space at $x_0$. Hence:
  \begin{equation}
    \tp{x - x_0}\cdot T = 0\text{ for all }T\in\R^n\text{ such that }T\cdot\grad\Psi(x_0) = 0
  \end{equation}
  and so we get equality in Cauchy-Schwarz:
  \begin{equation}\label{reverse_cauchy_schwarz}
    \tabs{\tp{x - x_0}\cdot\grad\Psi(x_0)} = \tabs{x - x_0}\tabs{\grad\Psi(x_0)}.
  \end{equation}
  Now, from the fundamental theorem of calculus and the fact that $\Psi(x_0) = 0$, we get the bound
  \begin{equation}\label{FTC_AGAIN}
    \tabs{\Psi(x)} \ge \tabs{\grad\Psi(x_0)\cdot\tp{x - x_0}} - \omega\tp{\lambda}\tabs{x - x_0}.
  \end{equation}
  Inequality~\eqref{kitty_} follows by combining~\eqref{reverse_cauchy_schwarz} with~\eqref{FTC_AGAIN}.

  Now if $x\in W_\ep$ is such that $\Psi(x) + \phi(x) = t$ with $\tnorm{\phi}_{C^{0,1}\tp{\Bar{U}}}\le\al\tp{\al\lambda}^3$ and $|t|\le\al\tp{\al\lambda}^3$, then~\eqref{kitty_} directly informs us that $\m{dist}\tp{x,M}\le\tp{\al\lambda}^3$. The proof is complete.
\end{proof}


We close this subsection with a technical estimate regarding the maps $\Theta_\phi$ that are defined in~\eqref{The coordinate maps}; the following result gives conditions under which these functions remain locally injective even when composed with a linear projection operator.

\begin{lem}[More reverse Lipschitz bounds]\label{lem on reverse Lipschitz bound with projection}
  Suppose that $\hat{m}\in M$ and $\nu\in\S^{n-1}$ satisfy $\nu\cdot\grad\Psi\tp{\hat{m}}\neq 0$. Choose $0<\kappa\le\lambda/4$ so small that
  \begin{equation}\label{definition of the parameter kappa}
    \omega\tp{2\kappa}\le\f{1}{2^{10}\be}\f{\tabs{\nu\cdot\grad\Psi\tp{\hat{m}}}}{\max\tcb{1,\tabs{\grad\Psi\tp{\hat{m}}}}}.
  \end{equation}

  Then for all $m,\tilde{m}\in M$, $\phi\in\mathcal{B}$, and $t\in[-\ep,\ep]$ satisfying
  \begin{equation}\label{wave_is_on_its_way}
    \max\tcb{|m - \hat{m}|,|\tilde{m} - \hat{m}|}\le\kappa,\;\max\tcb{\tnorm{\phi}_{C^{0,1}\tp{\Bar{U}}},|t|}\le\min\bcb{\f{1}{2^{10}}\f{\tabs{\nu\cdot\grad\Psi\tp{\hat{m}}}}{\tabs{\grad\Psi\tp{\hat{m}}}}\f{\ell}{n\be},\f{\kappa}{8\be}}
  \end{equation}
  the estimate
  \begin{equation}\label{rev_lip_proj}
    \tabs{m - \tilde{m}}\le4\f{\tabs{\grad\Psi\tp{\hat{m}}}}{\tabs{\nu\cdot\grad\Psi\tp{\hat{m}}}}\tabs{\P\tp{\Theta_\phi\tp{m,t} - \Theta_\phi\tp{\tilde{m},t}}}.
  \end{equation}
  holds with $\P = I - \nu\otimes\nu$.
\end{lem}
\begin{proof}
  We begin by deriving the claimed estimate~\eqref{rev_lip_proj} in the case that $t = 0$ and $\phi = 0$. In this case, the operator $\Theta$ collapses to the identity on $M$, see Lemma~\ref{lem on projection properties}. Fix $m,\tilde{m}\in \Bar{B(\hat{m},\kappa)}\cap M$. By arguing as in the proof of Lemma~\ref{lem on estimates on almost tangent vectors}, we estimate
  \begin{equation}\label{the_almost_tangent_estimate_redux}
    \tabs{\tp{m - \tilde{m}}\cdot\grad\Psi\tp{\hat{m}}}\le\omega\tp{\kappa}\tabs{m - \tilde{m}}.
  \end{equation}
  Combining~\eqref{the_almost_tangent_estimate_redux} with the decomposition $m - \tilde{m} = \nu\tp{m - \tilde{m}}\cdot\nu + \P\tp{m - \tilde{m}}$ shows that
  \begin{equation}
    \tabs{\nu\cdot\grad\Psi(\hat{m})}\tabs{\nu\cdot\tp{m - \tilde{m}}}\le\omega\tp{\kappa}\tabs{m - \tilde{m}} + \tabs{\P\tp{m - \tilde{m}}\cdot\grad\Psi\tp{\hat{m}}}.
  \end{equation}
  Hence from~\eqref{definition of the parameter kappa}, the triangle inequality, and
  $1+\tabs{\grad\Psi\tp{\hat{m}}}/\tabs{\nu\cdot\grad\Psi\tp{\hat{m}}}\le2\tabs{\grad\Psi\tp{\hat{m}}}/\tabs{\nu\cdot\grad\Psi\tp{\hat{m}}}$, we obtain
  \begin{equation}\label{estimate in the special case}
    \tabs{m - \tilde{m}}\le3\f{\tabs{\grad\Psi\tp{\hat{m}}}}{\tabs{\nu\cdot\grad\Psi\tp{\hat{m}}}}\tabs{\P\tp{m - \tilde{m}}}.
  \end{equation}
  
  Now suppose that $\phi\in\mathcal{B}$ and $t\in[-\ep,\ep]$. From the definition of $\Theta$~\eqref{The coordinate maps}, it follows that 
  \begin{multline}\label{__FiNaL__}
    \tabs{\P\tp{m - \tilde{m}}}\le\tabs{\P\tp{\Theta_\phi\tp{m,t} - \Theta_\phi\tp{\tilde{m},t}}} + \tabs{\P\tp{V(m)\tsb{\bf{f}\tp{\phi}}\tp{m,t} - V(\tilde{m})\tsb{\bf{f}\tp{\phi}}\tp{\tilde{m},t}}}\\
    \le\tabs{\P\tp{\Theta_\phi\tp{m,t} - \Theta_\phi\tp{\tilde{m},t}}} + \tp{n/\ell}\tabs{m - \tilde{m}}\tabs{\tsb{\bf{f}\tp{\phi}}\tp{m,t}} + \tabs{\tsb{\bf{f}\tp{\phi}}\tp{m,t} - \tsb{\bf{f}\tp{\phi}}\tp{\tilde{m},t}}.
  \end{multline}
  Then, from the third item of Lemma~\ref{lem on regularized distance to the manifold} and the hypotheses~\eqref{wave_is_on_its_way}, we find
  \begin{equation}\label{_not_done_yet_}
    \f{n}{\ell}\tabs{\tsb{\bf{f}\tp{\phi}}\tp{m,t}}\le2\f{n\be}{\ell}\tabs{\phi(m) - t}\le\f{1}{24}\f{\tabs{\nu\cdot\grad\Psi(\hat{m})}}{\tabs{\grad\Psi\tp{\hat{m}}}}.
  \end{equation}

  To estimate the difference $\tsb{\bf{f}\tp{\phi}}\tp{m,t} - \tsb{\bf{f}\tp{\phi}}\tp{\tilde{m},t}$, we note first that $\tabs{m-\tilde{m}}\le2\kappa\le\lambda/2$. We may therefore examine identity~\eqref{key_difference_identity} as in the derivation of estimate~\eqref{estimates_in_turn_again}. One finds that
  \begin{equation}\label{e_1__}
    \tabs{\tsb{\bf{f}\tp{\phi}}\tp{m,t} - \tsb{\bf{f}\tp{\phi}}\tp{\tilde{m},t}}\le4\be\tnorm{\phi}_{C^{0,1}\tp{\Bar{U}}}\tabs{m - \tilde{m}} + 2\be\babs{\tp{m - \tilde{m}}\cdot\int_0^1\grad\Psi\tp{\tau\Theta_\phi\tp{m,t} + \tp{1 - \tau}\Theta_\phi\tp{\tilde{m},t}}\;\m{d}\tau}.
  \end{equation}
  The final term above can then be estimated as
  \begin{equation}\label{e_2__}
    \babs{\tp{m - \tilde{m}}\cdot\int_0^1\grad\Psi\tp{\tau\Theta_\phi\tp{m,t} + \tp{1 - \tau}\Theta_\phi\tp{\tilde{m},t}}\;\m{d}\tau}\le\tabs{m - \tilde{m}}\omega\tp{2\be\tp{|t| + \tnorm{\phi}_{C^{0,1}\tp{\Bar{U}}}}}.
  \end{equation}
  Indeed, subtract the same integral along the segment joining $m$ and $\tilde m$, whose pairing with $m-\tilde m$ equals $\Psi(m)-\Psi(\tilde m)=0$, and then use the third item of Lemma~\ref{lem on regularized distance to the manifold}.

  The hypotheses imply
  \begin{equation}
    2\be\tp{|t|+\tnorm{\phi}_{C^{0,1}\tp{\Bar U}}}\le\kappa/2.
  \end{equation}
  Moreover, the two coefficients on the right-hand side of~\eqref{e_1__} are bounded respectively by
  \begin{equation}
    2^{-8}\f{\tabs{\nu\cdot\grad\Psi\tp{\hat m}}}{\tabs{\grad\Psi\tp{\hat m}}}
    \text{ and }
    2^{-9}\f{\tabs{\nu\cdot\grad\Psi\tp{\hat m}}}{\tabs{\grad\Psi\tp{\hat m}}}.
  \end{equation}
  Thus~\eqref{definition of the parameter kappa} and~\eqref{wave_is_on_its_way} give
  \begin{equation}\label{good_lippy_bound}
    \tabs{\tsb{\bf{f}\tp{\phi}}\tp{m,t} - \tsb{\bf{f}\tp{\phi}}\tp{\tilde{m},t}}\le\f{1}{24}\f{\tabs{\nu\cdot\grad\Psi(\hat{m})}}{\tabs{\grad\Psi\tp{\hat{m}}}}\tabs{m - \tilde{m}}.
  \end{equation}

  Combining estimates~\eqref{__FiNaL__}, \eqref{_not_done_yet_}, and~\eqref{good_lippy_bound} gives
  \begin{equation}
    \tabs{\P\tp{m-\tilde m}}\le\tabs{\P\tp{\Theta_\phi\tp{m,t}-\Theta_\phi\tp{\tilde m,t}}}+\f1{12}\f{\tabs{\nu\cdot\grad\Psi\tp{\hat m}}}{\tabs{\grad\Psi\tp{\hat m}}}\tabs{m-\tilde m}.
  \end{equation}
  Inserting this estimate into~\eqref{estimate in the special case} yields
  \begin{equation}
    \tabs{m-\tilde m}\le3\f{\tabs{\grad\Psi\tp{\hat m}}}{\tabs{\nu\cdot\grad\Psi\tp{\hat m}}}\tabs{\P\tp{\Theta_\phi\tp{m,t}-\Theta_\phi\tp{\tilde m,t}}}+\f14\tabs{m-\tilde m},
  \end{equation}
  and the absorption of the final term proves~\eqref{rev_lip_proj}.
\end{proof}

\subsection{Density estimates}\label{SS_density_estimates}

Let us record a well-known result in geometric measure theory, whose simple proof we include for the readers' convenience.
\begin{lem}[Density estimates on a fixed manifold]\label{lem on density estimates on a fixed manifold}
  Suppose that $1\le k<n$ and $M\subset \R^n$ is a compact $k$-dimensional manifold of class $C^1$. Then, there exists a constant $C\in\R^+$ such that for all $x\in\R^n$ and all $r\in\R^+$ we have the estimate
  \begin{equation}\label{density_estimate_1}
    \mathcal{H}^k\tp{M\cap B(x,r)}\le Cr^k
  \end{equation}
  where $\mathcal{H}^k$ denotes the $k$-dimensional Hausdorff measure.
\end{lem}
\begin{proof}
  By the local graph representation for class $C^1$ submanifolds, see, e.g., Section 3.1 in Federer~~\cite{MR257325}, the manifold $M$ is locally diffeomorphic to its tangent spaces. More precisely, for each $m\in M$, upon letting $A_m = \tcb{m + v\;:\;v\in T_mM}$, there exists $0<\rho_m<1$ and a $C^1$ diffeomorphism
  \begin{equation}
    f_m:B(m,\rho_m)\to f_m\tp{B(m,\rho_m)}\subset\R^n
  \end{equation}
  satisfying
  \begin{equation}
    \max\tcb{\tnorm{\grad f_m}_{L^\infty\tp{B(m,\rho_m)}},\tnorm{\grad f_m^{-1}}_{L^\infty\tp{f_m(B(m,\rho_m))}}}\le 2
  \end{equation}
  and
  \begin{equation}\label{flattening_identity}
    B(m,\rho_m/2)\cap M = B(m,\rho_m/2)\cap f_m^{-1}\tp{A_m}.
  \end{equation}

  By compactness of $M$, there exists $N\in\N^+$ and $\tcb{m_j}_{j=1}^N\subset M$ such that 
  \begin{equation}
    M \subset \bigcup_{j=1}^N B(m_j,\rho_{j})\text{ where we have defined }\rho_j = \f14\rho_{m_j}.
  \end{equation}

  Now suppose that $x\in\R^n$ and $0<r\in\R^+$. Thanks to identity~\eqref{flattening_identity} and Theorem 2.8 in Evans and Gariepy~\cite{MR3409135} we have the estimate
  \begin{equation}\label{eq_1}
    \mathcal{H}^k\tp{M\cap B(x,r)}\le\sum_{j=1}^N\mathcal{H}^k\tp{M\cap B(m_j,\rho_j)\cap B(x,r) }\le 2^k\sum_{j=1}^N\mathcal{H}^k\tp{A_{j}\cap f_j\tp{B(m_j,\rho_j)\cap B(x,r)}},
  \end{equation}
  in which we have abbreviated $f_j = f_{m_j}$ and $A_{j} = A_{m_j}$. The affine subspace $A_j$ is isometrically isomorphic to $\R^k$. Therefore, by the coincidence of Hausdorff and Lebesgue measures on $\R^k$ and the isodiametric inequality (see, e.g., Theorems 2.4 and 2.5 in~\cite{MR3409135}) we see that for each $j\in\tcb{1,\dots,N}$ it holds that
  \begin{equation}\label{eq_2}
    \m{diam}\tp{A_j\cap f_j\tp{B(m_j,\rho_j)\cap B(x,r)}}\le 4r\text{ and }\mathcal{H}^k\tp{A_{j}\cap f_j\tp{B(m_j,\rho_j)\cap B(x,r)}}\le\f{\tp{2\sqrt{\pi}r}^k}{\Upgamma\tp{k/2 + 1}},
  \end{equation}
  where $\Upgamma$ denotes Euler's gamma function. The desired estimate~\eqref{density_estimate_1} now follows by combining~\eqref{eq_1} with~\eqref{eq_2} and taking $C = N\tp{4\sqrt{\pi}}^k/\Upgamma\tp{k/2 + 1}$.
\end{proof}

We require a version of Lemma~\ref{lem on density estimates on a fixed manifold} in which the constant in the density estimate is stable under small perturbations of the underlying manifold. We shall specialize to the case of compact hypersurfaces. The stable local level set coordinates of Theorem~\ref{thm on stable local level set coordinates} readily dispose of this task. In the next series of results, we shall use the following notation: for $0\le q<\infty$ we introduce the function
\begin{equation}\label{weird_functions}
  \lambda^q:\R^n\setminus\tcb{0}\to\R\text{ via the action }\lambda^q(z) = \begin{cases}
    \log\tabs{z}&\text{if }q = 0,\\
    \tabs{z}^{-q}&\text{if }q>0.
  \end{cases}
\end{equation}

\begin{thm}[Uniform density estimates on families of hypersurfaces]\label{thm on uniform density estimates on families of hypersurfaces}
  Let $\R\ni q\ge0$, $\N\ni n\ge 2$, $\es\neq U\subset\R^n$ be open and bounded, $\Psi\in C^1_{\m{b}}\tp{\Bar{U}}$, and $a,b\in\R$ with $a\le b$. Assume that 
  \begin{equation}\label{hypotheses_are_here}
    E = \tcb{x\in\Bar{U}\;:\;\Psi(x)\in[a,b]},\;E\Subset U,\;\min_E\tabs{\grad\Psi}>0.
  \end{equation}
  Then, there exist constants $C,\ep\in\R^+$ such that the following hold for all $a - \ep \le t\le b+\ep$, $\phi\in C^1\tp{\Bar{U}}$ satisfying $\tnorm{\phi}_{C^{0,1}\tp{\Bar{U}}}\le \ep$, and $\rho\in\R^+$.
  \begin{enumerate}
    \item We have the inclusions
    \begin{equation}\label{INC_1}
      M(\phi,t) = \tcb{x\in U\;:\;\Psi(x) + \phi(x) = t}\subseteq\tcb{x\in U\;:\;a-2\ep\le\Psi(x)\le b+2\ep}\Subset U.
    \end{equation}
    \item If $M\tp{\phi,t}\neq\es$, we have the uniform positivity estimate
    \begin{equation}\label{POS_1}
      \min\tcb{\tabs{\grad\Psi(x) + \grad\phi(x)}\;:\;x\in M(\phi,t)}\ge1/C.
    \end{equation}
    Moreover, if $M\tp{0,s}\neq\es$ for all $a\le s\le b$, then, for $t\in[a,b]$,
    \begin{equation}\label{POS_2}
      \mathcal{H}^{n-1}\tp{M\tp{\phi,t}}\ge1/C.
    \end{equation}
    \item We have the uniform density estimate
    \begin{equation}\label{EST_1}
      \sup_{x\in\R^n}\mathcal{H}^{n-1}\tp{M\tp{\phi,t}\cap B(x,\rho)}\le C\rho^{n-1}.
    \end{equation}
    \item If $0\le q<n-1$, we have the estimate
    \begin{equation}\label{EST_2}
      \sup_{x\in U}\int_{M\tp{\phi,t}\cap B(x,\rho)}\tabs{\lambda^q\tp{x - z}}\;\m{d}\mathcal{H}^{n-1}\tp{z}\le C
      \begin{cases}
      \rho^{n- q -1}&\text{if }q>0,\\
      \rho^{n-1}\tp{1+\tabs{\log\rho}}&\text{if }q = 0.
      \end{cases}
    \end{equation}
    \item If $n-1\le q$, we have the estimate
    \begin{equation}\label{EST_3}
      \sup_{x\in U}\int_{M(\phi,t)\setminus B(x,\rho)}\tabs{\lambda^q\tp{x - z}}\;\m{d}\mathcal{H}^{n-1}\tp{z}\le C\begin{cases}
      \rho^{n-q-1}&\text{if }q>n-1,\\
      1 + \tabs{\log\tabs{\rho}}&\text{if }q = n-1.
      \end{cases}
    \end{equation}
  \end{enumerate}
\end{thm}
\begin{proof}
  By using that $E\Subset U$ and that $\min_{E}\tabs{\grad\Psi}>0$, we may select $\hat{a},\hat{b}\in\R$ such that $\hat{a}<a$, $b<\hat{b}$, and
  \begin{equation}\label{fat_set}
    \hat{E} = \tcb{x\in U\;:\;\Psi(x)\in[\hat{a},\hat{b}]}\Subset U,\text{ and }\min\tcb{\tabs{\grad\Psi(x)}\;:\;x\in\hat{E}}>0.
  \end{equation}

  Set $\mathcal{K}=\Psi\tp{\hat E}$, which is a compact subset of $[\hat a,\hat b]$. Now, for each $s\in\mathcal{K}$ consider the map $\Psi_s\in C_{\m{b}}^1\tp{\Bar{U}}$ defined via $\Psi_s(x) = \Psi(x) - s$; set also $M_s = \tcb{x\in U\;:\;\Psi_s(x) = 0}$. By the definition of $\mathcal K$, each $M_s$ is nonempty. The hypotheses~\eqref{hypotheses_are_here} ensure that Theorem~\ref{thm on stable local level set coordinates} is applicable to each $\Psi_s$. By carefully unpacking the conclusions of this result, we obtain positive constants $C_s, r_s, \ep_s, \tilde{\ep}_s, \hat{\ep}_s\in\R^+$ such that
  \begin{equation}\label{AI_SLOP}
    \ep_s<\tilde{\ep}_s<\hat{\ep}_s<\m{dist}(\hat{E},\pd U),\;\ep_s + r_s\le\min\tcb{\tabs{a -\hat{a}},\tabs{b - \hat{b}}},\;r_s\le\min_{\hat{E}}\tabs{\grad\Psi}/2
  \end{equation}  
  and the following hold.
  \begin{enumerate}
    \item For each $\phi\in C^1_{\m{b}}\tp{\Bar{U}}$ with $\tnorm{\phi}_{C^{0,1}\tp{\Bar{U}}}\le r_s$ there exists a map $\Theta_{\phi}^s:M_s\times[-\hat{\ep}_s,\hat{\ep}_s]\to U$ that is a $C^1$ diffeomorphism onto its image. Moreover, we have the bi-Lipschitz bound
    \begin{equation}\label{bi_lip}
      \tp{1/C_s}\tabs{m - \tilde{m},\tau - \tilde{\tau}}\le\tabs{\Theta^s_\phi\tp{m,\tau} - \Theta^s_\phi\tp{\tilde{m},\tilde{\tau}}}\le C_s\tabs{m - \tilde{m},\tau - \tilde{\tau}}
    \end{equation}
    and the level-set parametrization identity
    \begin{equation}\label{lev_par}
      \tp{\Psi_s + \phi}\circ\Theta_\phi^s\tp{m,\tau} = \tau
    \end{equation}
    for all $\tp{m,\tau},\tp{\tilde{m},\tilde{\tau}}\in M_s\times[-\hat{\ep}_s,\hat{\ep}_s]$.
    \item For each $\tau\in\R$ and $\phi\in C^1_{\m{b}}\tp{\Bar{U}}$ satisfying $|\tau|\le\ep_s$ and $\tnorm{\phi}_{C^{0,1}\tp{\Bar{U}}}\le r_s$ we have that 
    \begin{multline}\label{some_inclusions}
      \tcb{x\in U\;:\;|\Psi_s(x)|\le2\ep_s}\subseteq\tcb{x\in U\;:\;\m{dist}\tp{x,M_s}\le\hat{\ep}_s}\text{ and}\\
      \tcb{x\in U\;:\;\m{dist}\tp{x,M_s}\le\hat{\ep}_s,\;\Psi_s(x) + \phi(x) = \tau}\subseteq\tcb{x\in U\;:\;\m{dist}(x,M_s)\le\tilde{\ep}_s}\subseteq\Theta^s_\phi\tp{M_s\times[-\ep_s,\hat{\ep}_s]}\Subset U.
    \end{multline}
  \end{enumerate}

  The intervals $\tcb{t\in\R:\tabs{t-s}<\ep_s}$ form an open cover of $\mathcal K$. Hence there exist a finite set $\es\neq\digamma\subset\mathcal K$ and a number $\ep_0>0$ such that
  \begin{equation}\label{the_covering_condition}
    \tcb{t\in\R\;:\;\m{dist}\tp{t,\mathcal K}\le\ep_0}\subset\bigcup_{s\in\digamma}\tcb{t\in\R\;:\;|s - t|<\ep_s}.
  \end{equation}
  We then instantiate the positive parameter
  \begin{equation}
    \ep = \min\bcb{\ep_0,\f{a-\hat a}{3},\f{\hat b-b}{3},\min\tcb{\min\tcb{\ep_s,r_s}/3\;:\;s\in\digamma}}.
  \end{equation}

  We now prove the first item. If $x\in M(\phi,t)$ for some $\phi\in C^1_{\m{b}}\tp{\Bar{U}}$ with $\tnorm{\phi}_{C^{0,1}\tp{\Bar{U}}}\le\ep$ and $t\in[a-\ep,b + \ep]$ then it is clear that $a-2\ep\le\Psi(x)\le b + 2\ep$ and hence $x\in\hat{E}\Subset U$.

  Whenever $M\tp{\phi,t}\neq\es$, the second item's estimate~\eqref{POS_1} is immediate from the rightmost bound of equation~\eqref{AI_SLOP}.

  To justify the third item, we begin by letting $x\in\R^n$ and $\rho\in\R^+$ while fixing $t$ and $\phi$ as above. It suffices to consider the case that $M(\phi,t)\cap B(x,\rho)\neq\es$.
  
  Fix $x_\ast\in M\tp{\phi,t}\cap B(x,\rho)$. The first item gives $x_\ast\in\hat E$, and hence
  $\m{dist}\tp{t,\mathcal K}\le\tabs{t-\Psi\tp{x_\ast}}=\tabs{\phi\tp{x_\ast}}\le\ep_0$. From the covering condition~\eqref{the_covering_condition}, we deduce the existence of $s\in\digamma$ such that $|s - t|<\ep_s$. Then, $\tau = t - s$ satisfies $|\tau|\le\ep_s$ and hence we can invoke the inclusion~\eqref{some_inclusions}, the identity~\eqref{lev_par}, and the estimate~\eqref{bi_lip} to deduce
  \begin{equation}
    M\tp{\phi,t}\cap B(x,\rho) = \tcb{\Theta_\phi^s\tp{m,\tau}\;:\;m\in M_s(x,\rho)}\text{ where }M_s\tp{x,\rho} = \tcb{m\in M_s\;:\;\Theta_\phi^s\tp{m,\tau}\in B(x,\rho)}.
  \end{equation}
  Now we use Theorem 2.8 in Evans and Gariepy~\cite{MR3409135} in conjunction with the Lipschitz estimate~\eqref{bi_lip} to arrive at the intermediate inequality
  \begin{equation}\label{com_1}
    \mathcal{H}^{n-1}\tp{M\tp{\phi,t}\cap B(x,\rho)}\le\tp{C_s}^{n-1}\mathcal{H}^{n-1}\tp{M_s\tp{x,\rho}}.
  \end{equation}
  
  Next, we again employ the bounds~\eqref{bi_lip} to deduce that 
  \begin{equation}
    m\in M_s\tp{x,\rho}\imp m\in B(\hat{m},\tp{2C_s}\rho)\cap M_s,
  \end{equation}
  where $\hat{m}\in M_s\tp{x,\rho}$ is arbitrarily fixed. In turn, Lemma~\ref{lem on density estimates on a fixed manifold} gives us the existence of a finite constant $\tilde{C}_s$, depending only on the fixed compact hypersurface $M_s$, such that 
  \begin{equation}\label{com_2}
    \mathcal{H}^{n-1}\tp{M_s\tp{x,\rho}}\le\mathcal{H}^{n-1}\tp{M_s\cap B(\hat{m},2C_s\rho)}\le\tilde{C}_s\tp{2 C_s\rho}^{n-1}.
  \end{equation}
  Combining~\eqref{com_1} and~\eqref{com_2} gives the desired density estimate~\eqref{EST_1} with the finite constant 
  \begin{equation}
    C = 2^{n-1}\max_{s\in\digamma}\tilde{C}_s\tp{C_s}^{2n - 2}.
  \end{equation}

  The strategy for the claimed estimate~\eqref{POS_2} is similar, yet simpler. Assume the additional nonemptiness hypothesis and fix $t\in[a,b]$ and $\phi$ as above. Then $M\tp{0,t}\neq\es$, so $t\in\mathcal K$. From condition~\eqref{the_covering_condition}, we are given $s\in\digamma$ such that $|s - t|<\ep_s$; then, $\tau = t - s$ again obeys $|\tau|\le\ep_s$ thus allowing us to employ equations~\eqref{lev_par} and~\eqref{some_inclusions} and deduce 
  \begin{equation}
    M\tp{\phi,t} = \tcb{\Theta^s_\phi\tp{m,\tau}\;:\;m\in M_s}.
  \end{equation}
  Now we use yet again Theorem 2.8 in~\cite{MR3409135} in conjunction with the Lipschitz estimate~\eqref{bi_lip} to estimate
  \begin{equation}
    \mathcal{H}^{n-1}\tp{M_s}\le\tp{C_s}^{n-1}\mathcal{H}^{n-1}\tp{M\tp{\phi,t}}
  \end{equation}
  and hence 
  \begin{equation}
    \mathcal{H}^{n-1}\tp{M\tp{\phi,t}}\ge\min\tcb{\tp{C_s}^{1-n}\mathcal{H}^{n-1}\tp{M_s}\;:\;s\in\digamma}>0.
  \end{equation}
  Note that $M_s\neq\es$ and hence $\mathcal{H}^{n-1}\tp{M_s}>0$ for all $s\in\digamma$, since $\digamma\subset\mathcal K$.

  To prove the fourth item, we use the third item and dyadic decomposition. For $\ell\in\Z$ and $x\in\R^n$ we let 
  \begin{equation}\label{annulus_region}
    A_\ell\tp{x,\rho} = B(x,2^{-\ell}\rho)\setminus B(x,2^{-\ell-1}\rho),
  \end{equation}
  and expand the integral as a sum over the regions~\eqref{annulus_region}. Upon doing so, we have the estimate
  \begin{equation}
    \int_{M\tp{\phi,t}\cap B(x,\rho)}\f{1}{\tabs{x - z}^q}\;\m{d}\mathcal{H}^{n-1}\tp{z}\le\rho^{-q}\sum_{\ell=0}^\infty2^{q\tp{\ell+1}}\mathcal{H}^{n-1}\tp{M\tp{\phi,t}\cap B(x,2^{-\ell}\rho)}\le\f{2^qC}{1 - 2^{q-n+1}}\rho^{n-q-1}.
  \end{equation}
  This implies the desired bound~\eqref{EST_2} in the case that $q>0$. The case $q = 0$ follows by a nearly identical argument.

  The fifth item is also very similar. If $M\tp{\phi,t}=\es$, there is nothing to prove. Otherwise, let $x\in U$ and set 
  \begin{equation}
    K = \min\tcb{k\in\N\;:\;M(\phi,t)\subset B(x,2^k\rho)}.
  \end{equation}
  By dyadic decomposition again, we have
  \begin{multline}\label{slighty}
    \int_{M\tp{\phi,t}\setminus B(x,\rho)}\f{1}{\tabs{x - z}^q}\;\m{d}\mathcal{H}^{n-1}\tp{z}\le\rho^{-q}\sum_{\ell=1}^K 2^{q\tp{-\ell + 1}}\mathcal{H}^{n-1}\tp{M\tp{\phi,t}\cap B(x,2^{\ell}\rho)}\\\le\begin{cases}
      \f{2^{n-1}C}{1 - 2^{n-1-q}}\rho^{n-q-1}&\text{if }n-1<q,\\
      2^qCK&\text{if }n-1 = q.
    \end{cases}
  \end{multline}
  We obtain the sought-after estimate~\eqref{EST_3} from~\eqref{slighty} by estimating  $K \le 1 + \tabs{\log_2\tp{\m{diam}U/\rho}}$.
\end{proof}

We now use the density estimates of Theorem~\ref{thm on uniform density estimates on families of hypersurfaces} to establish a number of useful potential bounds over thin domains.

\begin{coro}[Integral bounds on thin domains]\label{bounds for integral operators on thin domains}
  Let $\R\ni q\ge0$ satisfy $q<n$, $\N\ni n\ge2$, $\es\neq U\subset\R^n$ be open and bounded, $\Psi\in C^1_{\m{b}}\tp{\Bar{U}}$, and $a,b\in\R$ with $a\le b$. Assume that~\eqref{hypotheses_are_here} holds and let $\ep>0$ be the small parameter guaranteed to exist by Theorem~\ref{thm on uniform density estimates on families of hypersurfaces}. There exists a constant $C\in\R^+$ such that for all $t\in[a,b]$, $\phi\in C^1_{\m{b}}\tp{\Bar{U}}$ with $\tnorm{\phi}_{C^{0,1}\tp{\Bar{U}}}\le\ep$, $0<\del\le\ep$, and $x,z\in U$ the following items hold, where we implement the shorthand
  \begin{equation}\label{short_hand}
    \tcb{\tabs{\Psi + \phi - t}\le\del} = \tcb{x\in U\;:\;\tabs{\Psi(x) + \phi(x) - t}\le\del}.
  \end{equation}
  \begin{enumerate}
    \item We have the bounds
    \begin{equation}\label{thin_est_1}
      \f{1}{\del}\int_{\tcb{\tabs{\Psi + \phi - t}\le\del}}\tabs{\lambda^q\tp{x - y}}\;\m{d}\mathcal{H}^n\tp{y}\le C\begin{cases}
        1&\text{if }0\le q<n-1,\\
        1+\tabs{\log\del}&\text{if }q = n-1,\\
        \del^{n-1-q}&\text{if }n-1<q<n.
      \end{cases}
    \end{equation}
    \item If $0\le q<n-1$, we have the uniform continuity estimates 
    \begin{equation}\label{thin_est_2}
      \f{1}{\del}\int_{\tcb{\tabs{\Psi + \phi - t}\le\del}}\tabs{\lambda^q\tp{x - y} - \lambda^q\tp{z - y}}\;\m{d}\mathcal{H}^n\tp{y}\le C\begin{cases}
        \tabs{x - z}&\text{if }0\le q<n-2,\\
        \tabs{x - z}\tp{1 + \tabs{\log\tabs{x - z}}}&\text{if }q = n-2,\\
        \tabs{x - z}^{n-1-q}&\text{if }n-2<q<n-1.
      \end{cases}
    \end{equation}
  \end{enumerate}
\end{coro}
\begin{proof}
  To prove the first item, we begin by splitting the integral and employing the coarea formula; for a reference to the latter, the reader is referred to, e.g., Theorem 3.2.12 in Federer~\cite{MR257325}. The equality
  \begin{multline}\label{coarea_equality}
    \f{1}{\del}\int_{\tcb{\tabs{\Psi + \phi - t}\le\del}}\tabs{\lambda^q\tp{x - y}}\;\m{d}\mathcal{H}^n\tp{y} =\f{1}{\del}\int_{\tcb{\tabs{\Psi + \phi - t}\le\del}\cap B(x,\del)}\tabs{\lambda^q\tp{x - y}}\;\m{d}\mathcal{H}^n\tp{y}\\+ \f{1}{\del}\int_{-\del}^\del\int_{M\tp{\phi,t+s}\setminus B(x,\del)}\f{\tabs{\lambda^q\tp{x - y}}}{\tabs{\grad\tp{\Psi + \phi}\tp{y}}}\;\m{d}\mathcal{H}^{n-1}\tp{y}\;\m{d}\mathcal{H}^1\tp{s}
  \end{multline}
  is obtained; note that we are tacitly using the first and second items of Theorem~\ref{thm on uniform density estimates on families of hypersurfaces}.
  
  The first term on the right of~\eqref{coarea_equality} is estimated using radial spherical coordinates (e.g. 3.2.13 in~\cite{MR257325}):
  \begin{multline}\label{EE_1}
    \f{1}{\del}\int_{\tcb{\tabs{\Psi + \phi - t}\le\del}\cap B(x,\del)}\tabs{\lambda^q\tp{x - y}}\;\m{d}\mathcal{H}^n\tp{y}\le \f{2\pi^{n/2}}{\Upgamma\tp{n/2}}\f{1}{\del}\int_0^\del s^{n-1}\tabs{\lambda^q\tp{s}}\;\m{d}\mathcal{H}^1\tp{s}\\\le C\begin{cases}
    \del^{n-1-q}&\text{if }0<q<n,\\
    \del^{n-1}\tp{1 + \tabs{\log\del}}&\text{if }q = 0.
    \end{cases}
  \end{multline}
  The second term on the right hand side of~\eqref{coarea_equality} is instead estimated via the third, fourth, and fifth items of Theorem~\ref{thm on uniform density estimates on families of hypersurfaces}:
  \begin{equation}\label{EE_2}
    \f{1}{\del}\int_{-\del}^\del\int_{M\tp{\phi,t+s}\setminus B(x,\del)}\f{\tabs{\lambda^q\tp{x - y}}}{\tabs{\grad\tp{\Psi + \phi}\tp{y}}}\;\m{d}\mathcal{H}^{n-1}\tp{y}\;\m{d}\mathcal{H}^1\tp{s}\le C\begin{cases}
      1 &\text{if }0\le q<n-1,\\
      1 + \tabs{\log\del}&\text{if }q = n-1,\\
      \del^{n-1 -q}&\text{if }n-1<q<n.
    \end{cases}
  \end{equation}
  The proof of the first item is completed by combining inequalities~\eqref{EE_1} and~\eqref{EE_2}.

  We now prove the second item. It is sufficient to address the case that $x\neq z$. Define the region
  \begin{equation}\label{the_important_disk}
    D(x,z) = \tcb{y\in\R^n\;:\;\tabs{y - (x+z)/2}\le2\tabs{x - z}}.
  \end{equation}
  The following simple observations are useful. The inclusion
  \begin{equation}\label{observation_1}
    D(x,z)\subset B(x,5|x - z|/2)\cap B(z,5|x - z|/2)
  \end{equation}
  holds. The estimates
  \begin{equation}\label{observation_2}
    \forall\;y\in\R^n\setminus D(x,z),\;\forall\;\tau\in[0,1],\;\tabs{\tau x + \tp{1 - \tau}z - y}\ge 3\tabs{y - (x + z)/2}/4.
  \end{equation}
  hold.

  We now split the integral on the left hand side of~\eqref{thin_est_2} according to the region~\eqref{the_important_disk}. When $y\in D(x,z)$, we do not exploit the difference and instead just use the coarea formula again, the second and fourth items of Theorem~\ref{thm on uniform density estimates on families of hypersurfaces}, and observation~\eqref{observation_1}. In doing so, we learn that
  \begin{equation}\label{part_1}
    \f{1}{\del}\int_{\tcb{\tabs{\Psi + \phi - t}\le\del}\cap D(x,z)}\tabs{\lambda^q\tp{x - y} - \lambda^q\tp{z - y}}\;\m{d}\mathcal{H}^n\tp{y}\le C\begin{cases}
    \tabs{x - z}^{n-1}\tp{1 + \tabs{\log\tabs{x - z}}}&\text{if }q = 0,\\
    \tabs{x - z}^{n-1-q}&\text{if }0<q<n-1.
    \end{cases}
  \end{equation}
  
  On the other hand, when $y\in \R^n\setminus D(x,z)$, we do exploit the difference between $\lambda^q\tp{x-y}$ and $\lambda^q\tp{z-y}$ via the fundamental theorem of calculus. Thanks to observation~\eqref{observation_2}, this results in the bound
  \begin{equation}\label{diff_bou}
    \tabs{\lambda^q\tp{x - y} - \lambda^q\tp{z - y}}\le C\tabs{x - z}\tabs{\lambda^{q+1}\tp{y - (x+z)/2}}.
  \end{equation}
  Upon combining~\eqref{diff_bou} with the coarea formula and the second, fourth, and fifth items of Theorem~\ref{thm on uniform density estimates on families of hypersurfaces}, one finds that 
  \begin{equation}\label{part_2}
    \f{1}{\del}\int_{\tcb{\tabs{\Psi + \phi - t}\le\del}\setminus D(x,z)}\tabs{\lambda^q\tp{x - y} - \lambda^q\tp{z - y}}\;\m{d}\mathcal{H}^{n}\tp{y}\le C\begin{cases}
      |x - z|&\text{if }q<n-2,\\
      |x - z|\tp{1 + \tabs{\log\tabs{x - z}}}&\text{if }q = n-2,\\
      \tabs{x - z}^{n-1-q}&\text{if }n-2<q<n-1.
    \end{cases}
  \end{equation}
  Summing~\eqref{part_1} with~\eqref{part_2} completes the proof of the second item.
\end{proof}

Our final density result uniformly estimates the measure of the intersection of certain cylindrical regions with families of hypersurfaces.
\begin{prop}[Transverse cylindrical density estimates]\label{prop on transverse cylinder density estimate}
  Let $\N\ni n\ge 2$, $\es\neq U\subset\R^n$ be open and bounded, $\Psi\in C^1_{\m{b}}\tp{\Bar{U}}$, $a,b\in\R$ with $a\le b$, and assume that~\eqref{hypotheses_are_here} holds. Fix $p\in\R^n$ and $\nu\in\S^{n-1}$ and assume additionally that for some $\sig_0\in\R$ we have $p + \sig_0\nu\in E$ and
  \begin{equation}\label{transverse intersection condition}
    \forall\;\sig\in\R,\;p + \sig\nu\in E\imp\nu\cdot\grad\Psi\tp{p + \sig\nu}\neq0.
  \end{equation}
  Then, there exist constants $C,\ep\in\R^+$ such that for all $a-\ep\le t\le b + \ep$, $\phi\in C^1\tp{\Bar{U}}$ satisfying $\tnorm{\phi}_{C^{0,1}\tp{\Bar{U}}}\le\ep$, and $\rho\in\R^+$ we have the estimate
  \begin{equation}\label{cylindrical_density_estimate}
    \mathcal{H}^{n-1}\tp{M(\phi,t)\cap\tcb{x\in\R^n\;:\;\tabs{\P\tp{x - p}}\le\rho}}\le C\rho^{n-1}
  \end{equation}
  where $M(\phi,t)\Subset U$ is the hypersurface defined on the left hand side of equation~\eqref{INC_1} and $\P = I - \nu\otimes\nu$.
\end{prop}
\begin{proof}
  We begin by arguing in a similar fashion to the first part of the proof of Theorem~\ref{thm on uniform density estimates on families of hypersurfaces}. Let $\hat{a}$, $\hat{b}$, and $\hat{E}$ be as in~\eqref{fat_set}, but satisfying additionally
  \begin{equation}\label{additional_condition}
    \forall\;\sig\in\R\text{ such that }p + \sig \nu\in\hat{E},\text{ we have }\nu\cdot\grad\Psi\tp{p + \sig\nu}\neq 0.
  \end{equation}
  After shrinking the enlargement in~\eqref{fat_set}, if necessary, this condition can be arranged by continuity and compactness. Indeed, otherwise there would exist a sequence $x_j=p+\sig_j\nu\in\Bar U$ such that $\m{dist}\tp{\Psi\tp{x_j},[a,b]}\to0$ and $\nu\cdot\grad\Psi\tp{x_j}=0$. A convergent subsequence would have a limit in $E$, contradicting~\eqref{transverse intersection condition}.

  Set $\mathcal{K}=\Psi\tp{\hat E}$. For each $s\in\mathcal{K}$ we again define the map $\Psi_s\in C^1_{\m{b}}\tp{\Bar{U}}$ via $\Psi_s(x) = \Psi(x) - s$ and the nonempty set $M_s=\tcb{x\in U\;:\;\Psi_s\tp{x} = 0}$. Again, each $\Psi_s$ is within the scope of Theorem~\ref{thm on stable local level set coordinates}; invoking this result gives us positive constants $C_s,r_s,\ep_s,\tilde{\ep}_s,\hat{\ep}_s\in\R^+$ (with $C_s\ge 1$) and a function $\Theta^s_\phi$ such that the inequalities of equation~\eqref{AI_SLOP} and the two items listed in that part of the proof of Theorem~\ref{thm on uniform density estimates on families of hypersurfaces} are true. By also appealing to the third item of Lemma~\ref{lem on regularized distance to the manifold} in the process, we can enlarge the constant $C_s$, if necessary, to arrange that 
  \begin{equation}\label{we_need_this}
    \tabs{\Theta_\phi^s\tp{m,\tau} - m}\le C_s\tabs{\Psi_s(m) + \phi\tp{m} - \tau}\text{ for all }\tp{m,\tau}\in M_s\times[-\hat{\ep}_s,\hat{\ep}_s].
  \end{equation}

  Before we appeal to compactness in some way to reduce to a finite number of model level sets, we first need to import the conclusions of Lemma~\ref{lem on reverse Lipschitz bound with projection}. 
  
  We now note that for any $s\in \mathcal{K}$, the set 
  \begin{equation}
    L_s = \tcb{x\in M_s\;:\;\exists\;\sig\in\R\;:\;p + \sig\nu = x}
  \end{equation}
  is necessarily finite: $L_s$ is evidently compact and each member of this set is isolated. Indeed, for each $x\in L_s$ the condition~\eqref{additional_condition} guarantees that $\nu\cdot\grad\Psi_s\tp{x}\neq0$. The implicit function theorem can then be applied to completely describe the zero set $M_s$ in an open neighborhood of $x$ as a graph $y\mapsto y + \nu\varphi(y)$ for $y\cdot\nu = 0$ and $\varphi$ a continuously differentiable scalar function. Such a local graph can only have one intersection with the line $\tcb{p + \sig\nu\;:\;\sig\in\R}$.

  With $s\in\mathcal{K}$ fixed, if $L_s\neq\es$, we can apply Lemma~\ref{lem on reverse Lipschitz bound with projection} (whose hypotheses are satisfied thanks again to~\eqref{additional_condition}) a finite number of times, once for each member of $L_s$. After executing this procedure, we acquire positive parameters $\kappa_s,\del_s,\mathsf{C}_s\in\R^+$, with $\del_s<\min\tcb{\ep_s,r_s}$, such that for all $\hat{m}\in L_s$, $m,\tilde{m}\in M_s$, $\phi\in C^1_{\m{b}}\tp{\Bar{U}}$, and $\tau\in\R$ obeying the inequalities 
  \begin{equation}
    \max\tcb{\tabs{m - \hat{m}},\tabs{\tilde{m} - \hat{m}}}\le\kappa_s,\;\max\tcb{\tnorm{\phi}_{C^{0,1}\tp{\Bar{U}}},\tabs{\tau}}\le\del_s
  \end{equation} 
  we have the reverse Lipschitz estimate
  \begin{equation}\label{REV_rev_LIP}
    \tabs{m - \tilde{m}}\le\mathsf{C}_s\tabs{\P\tp{\Theta^s_\phi\tp{m,\tau} - \Theta^s_\phi\tp{\tilde{m},\tau}}},
  \end{equation}
  where $\P = I - \nu\otimes\nu$. If $L_s=\es$, choose any $\kappa_s,\mathsf C_s\in\R^+$ and $0<\del_s<\min\tcb{\ep_s,r_s}$; the preceding assertion is then vacuous.

  We shall also require the following simple auxiliary fact. For each $s\in\mathcal K$, there exists $0<\hat{\kappa}_s\le\kappa_s/2$ with the property that 
  \begin{equation}\label{simple_auxiliary_fact}
    M_s\cap\tcb{x\in\R^n\;:\;\tabs{\P\tp{x - p}}\le2\hat{\kappa}_s}\subset\bigcup_{\hat{m}\in L_s}B(\hat{m},\kappa_s).
  \end{equation}
  If $L_s=\es$, this follows from compactness, since $\min_{m\in M_s}\tabs{\P\tp{m-p}}>0$. If $L_s\neq\es$, this assertion is readily proved via contradiction. If no such $\hat{\kappa}_s$ existed, then we could construct a sequence $\tcb{x_\mu}_{\mu = 0}^\infty\subset M_s$ satisfying
  \begin{equation}\label{the condition on the}
    \lim_{\mu\to\infty}\tabs{\P\tp{x_\mu - p}} = 0\text{ and }\forall\;\hat{m}\in L_s,\;\inf_{\mu\in\N}\tabs{x_\mu - \hat{m}}\ge\kappa_s.
  \end{equation}
  Thanks to the compactness of $M_s$, we can extract a convergent subsequence to some limit point $x_\infty\in M_s$. The condition on the left hand side of~\eqref{the condition on the} ensures that for some $\sig\in\R$ we must have $x_\infty = p + \sig\nu$. But then $x_\infty\in L_s$, which is impossible by the rightmost condition of~\eqref{the condition on the}. So the sought-after $0<\hat{\kappa}_s\le\kappa_s/2$ indeed exists.

  The corresponding open intervals cover the compact set $\mathcal K$. Hence there exist a finite set $\es\neq\digamma\subset\mathcal K$ and a number $\ep_0>0$ such that
  \begin{equation}\label{yet_another_compactness}
    \tcb{t\in\R\;:\;\m{dist}\tp{t,\mathcal K}\le\ep_0}\subset\bigcup_{s\in\digamma}\bcb{t\in\R\;:\;\tabs{t - s}<\min\bcb{\del_s,\f{\hat{\kappa}_s}{2C_s}}}.
  \end{equation}
  Let $\ep_{\m d}>0$ be an admissible perturbation size in Theorem~\ref{thm on uniform density estimates on families of hypersurfaces}, applied with $q=0$. By using this finite set $\digamma$, we may then instantiate the finite and positive parameters
  \begin{multline}
    \ep=\min\bcb{\ep_0,\ep_{\m d},\f{a-\hat a}{3},\f{\hat b-b}{3},
      \min\bcb{\min\bcb{\del_s,\f{\hat{\kappa}_s}{2C_s}}\;:\;s\in\digamma}},\;
    \varrho=\min\tcb{\hat{\kappa}_s\;:\;s\in\digamma},\\\text{and }
    \mathsf{C}=\max\tcb{\max\tcb{\mathsf{C}_s,C_s}\;:\;s\in\digamma}.
  \end{multline}

  We are now ready to prove the density estimate~\eqref{cylindrical_density_estimate}. Let $t\in[a - \ep,b +\ep]$, $\phi\in C^1_{\m{b}}\tp{\Bar{U}}$ satisfy $\tnorm{\phi}_{C^{0,1}\tp{\Bar{U}}}\le\ep$, and $\rho\in\R^+$. Consider first the case that $\rho\ge\varrho$. Then, we crudely bound the Hausdorff measure on the left hand side of~\eqref{cylindrical_density_estimate} as
  \begin{equation}
    \mathcal{H}^{n-1}\tp{M\tp{\phi,t}\cap\tcb{x\in\R^n\;:\;\tabs{\P\tp{x - p}}\le\rho}}\le\mathcal{H}^{n-1}\tp{M(\phi,t)}\le\tp{\Bar{C}/\varrho^{n-1}}\rho^{n-1},
  \end{equation}
  where $\Bar{C}\in\R^+$ is the finite positive constant obtained from the right hand side of the third item of Theorem~\ref{thm on uniform density estimates on families of hypersurfaces} (evaluated when the radius is taken to be $\m{diam}\tp{U}$). 

  Now suppose that $0<\rho<\varrho$. If the set on the left-hand side of~\eqref{cylindrical_density_estimate} is empty, there is nothing to prove. Otherwise, choose $x_\ast$ in this set. Then $x_\ast\in\hat E$ and $\m{dist}\tp{t,\mathcal K}\le\tabs{\phi\tp{x_\ast}}\le\ep_0$. From the finite covering~\eqref{yet_another_compactness}, there exists $s\in\digamma$ with
    $|t - s|<\min\tcb{\del_s,\hat\kappa_s/(2C_s)}<\ep_s$. In turn, we may invoke the inclusion~\eqref{some_inclusions}, the identity~\eqref{lev_par}, and the estimate~\eqref{bi_lip} to deduce that
  \begin{equation}
    M\tp{\phi,t}\cap\tcb{x\in\R^n\;:\;\tabs{\P\tp{x - p}}\le\rho} = \tcb{\Theta^s_\phi\tp{m,t - s}\;:\;m\in \tilde{M}_s\tp{\rho}}
  \end{equation} 
  where
  \begin{equation}
    \tilde{M}_s\tp{\rho} = \tcb{m\in M_s\;:\;\tabs{\P\tp{\Theta_\phi^s\tp{m,t - s} - p}}\le\rho}.
  \end{equation}
  By Theorem 2.8 in Evans and Gariepy~\cite{MR3409135} and the Lipschitz bounds~\eqref{bi_lip}, we thus have
  \begin{equation}
    \mathcal{H}^{n-1}\tp{M\tp{\phi,t}\cap\tcb{x\in\R^n\;:\;\tabs{\P\tp{x - p}}\le\rho}}\le\mathsf{C}^{n-1}\mathcal{H}^{n-1}\tp{\tilde{M}_s\tp{\rho}}.
  \end{equation}
  So it suffices to measure $\tilde{M}_s\tp{\rho}$. Notice that if $m\in\tilde{M}_s\tp{\rho}$, then from estimate~\eqref{we_need_this}, we get 
  \begin{equation}
    \tabs{\P\tp{m - p}}\le\rho + \tabs{\Theta_\phi^s\tp{m,t-s} - m}\le\rho+C_s\tp{\tnorm{\phi}_{C^{0,1}\tp{\Bar U}}+\tabs{t-s}}<2\hat{\kappa}_s.
  \end{equation}
  In turn, appealing to the inclusion~\eqref{simple_auxiliary_fact}, we learn that 
  \begin{equation}
    \tilde{M}_s\tp{\rho} \subset\bigcup_{\hat{m}\in L_s}B\tp{\hat{m},\kappa_s}\cap\tilde{M}_s\tp{\rho}.
  \end{equation}
  For every $\hat m\in L_s$ for which $B\tp{\hat m,\kappa_s}\cap\tilde M_s\tp{\rho}$ is nonempty, choose
  $m_{\hat m}\in B\tp{\hat m,\kappa_s}\cap\tilde M_s\tp{\rho}$. If $m$ belongs to the set $B(\hat{m},\kappa_s)\cap\tilde{M}_s\tp{\rho}$ as well, then~\eqref{REV_rev_LIP} gives
  \begin{equation}
    \tabs{m-m_{\hat m}}\le\mathsf C_s\tabs{\P\tp{\Theta_\phi^s\tp{m,t-s}-\Theta_\phi^s\tp{m_{\hat m},t-s}}}\le2\mathsf C_s\rho.
  \end{equation}
  Consequently,
  \begin{equation}
    B\tp{\hat m,\kappa_s}\cap\tilde M_s\tp{\rho}\subset B\tp{m_{\hat m},2\mathsf C_s\rho}\cap M_s.
  \end{equation}
  Synthesizing these observations, we deduce that
  \begin{equation}
    \mathcal{H}^{n-1}\tp{\tilde{M}_s\tp{\rho}}\le\sum_{\substack{\hat{m}\in L_s\\B(\hat m,\kappa_s)\cap\tilde M_s(\rho)\neq\es}}\mathcal{H}^{n-1}\tp{B\tp{m_{\hat m},2\mathsf{C}_s\rho}\cap M_s}.
  \end{equation}
  The desired density estimate~\eqref{cylindrical_density_estimate} now follows from Lemma~\ref{lem on density estimates on a fixed manifold}, since $\digamma$ is finite and each $L_s$ is finite.
\end{proof}
\section*{Statements and declarations}

\paragraph{Competing interests}
The authors have no competing interests to declare that are relevant to
the content of this article.

\paragraph{Data availability}
Data sharing is not applicable to this article as no datasets were generated
or analyzed

\paragraph{Use of generative AI}
The authors acknowledge the use of generative AI tools to assist with proofreading and the preparation of figures. All mathematical ideas, arguments, proofs, and conclusions in this work were developed and verified by the authors.

\bibliographystyle{abbrv}
\bibliography{bib.bib}
\end{document}